\documentclass[a4paper,11pt]{article}

\usepackage{latexsym}
\usepackage{amsmath,amsfonts,amssymb,amsthm}
\newtheorem{thm}{Theorem}[section]

\newtheorem{cor}[thm]{Corollary}
\newtheorem{prop}[thm]{Proposition}

\newtheorem{rem}[thm]{Remark}
\newtheorem{defin}[thm]{Definition}
\input epsf

\title{Recurrence matrices}
\author{Roland Bacher, roland.bacher@ujf-grenoble.fr
}
\date{}
\begin{document}
\maketitle



\section{Introduction}

The subject of this paper are {\it $p-$recurrence matrices} 
(or recurrence matrices for short) over a fixed 
ground field ${\mathbf K}$. A recurrence 
matrix is an element of the product
$$\prod_{l=0}^\infty {\mathbf K}^{p^l\times q^l}$$
(where ${\mathbf K}^{p^l\times q^l}$ denotes the vector-space
of $p^l\times q^l$ matrices with
coefficients in ${\mathbf K}$) satisfying finiteness
conditions which are suitable for computations. Concretely,
a recurrence matrix has a finite description involving finitely many 
elements in ${\mathbf K}$, the set $\hbox{Rec}_{p\times q}({\mathbf
  K})$ of all recurrence matrices in 
$\prod_{l=0}^\infty {\mathbf K}^{p^l\times q^l}$
is a vector space and the obvious product 
$AB\in \prod_{l=0}^\infty {\mathbf K}^{p^l\times q^l}$ of two
recurrence matrices $A\in \hbox{Rec}_{p\times r}({\mathbf
  K}),B\in  \hbox{Rec}_{r\times q}({\mathbf
  K})$ is again a recurrence matrix.
The subset $\hbox{Rec}_{p\times p}({\mathbf K})\subset 
\prod_{l=0}^\infty {\mathbf K}^{p^l\times p^l}$ 
of recurrence matrices of ``square-size'' is thus an algebra. 
We show how to do computations in this algebra
and describe a few features of it.

The set of all invertible elements in the algebra
$\hbox{Rec}_{p\times p}({\mathbf K})$ forms a group
$\hbox{GL}_{p-rec}({\mathbf K})$ containing interesting subgroups.
Indeed, a recurrence matrix $\hbox{GL}_{p-rec}({\mathbf K})$
is closely linked to a
finite-dimensional matrix-representation of the free
monoid on $p^2$ elements. Recurrence matrices for which
the image of this representation is a finite monoid
are in bijection with (suitably defined) 
``automatic functions'' associated to finite-state
automata. This implies easily that 
$\hbox{GL}_{p-rec}({\mathbf K})$ contains all
``automata groups'' or
``$p-$self-similar groups'' (formed by bijective
finite-state transducers acting by automorphisms
on the infinite plane rooted $p-$regular tree).

In particular, the group
$\hbox{GL}_{2-rec}({\mathbf K})$ contains a
famous group of Grigorchuk, see \cite{Gri} (and all similarly 
defined groups) in a natural way.

Part of the present paper is contained in a
condensed form in \cite{Bac} whose main result was the initial
motivation for developping the theory.

\section{Monoids}\label{xectmonoids}

{\bf Definition} A {\it monoid} (or {\it a semi-group with identity})
is a set $\mathcal M$ endowed with an
associative product $\mathcal M\times \mathcal M\longrightarrow
\mathcal M$ admitting a two-sided identity.

\begin{rem} It is enough to require that a monoid admits a left
  identity $\epsilon$ and a right identity $f$ for its associative
product since then $ef=e=f$. This shows moreover unicity of the
identity.
\end{rem}

In a {\it commutative} monoid, the product is commutative. A
{\it morphism} of monoids 
$\mu:\mathcal M\longrightarrow \tilde{\mathcal M}$
is an application such that $\mu(e)=\tilde e$ and
$\mu(ab)=\mu(a)\mu(b)$ where $e$ is the identity of $\mathcal M$,
$\tilde e$ the identity of $\tilde{\mathcal M}$ and where
$a,b\in\mathcal M$ are arbitrary. A {\it submonoid} of a monoid
$\mathcal M$ is a subset which is closed under the product and
contains the identity of $\mathcal M$. A {\it (linear) representation}
of a monoid $\mathcal M$ is a morphism $\rho:\mathcal M\longrightarrow
\hbox{End}(V)$ from $\mathcal M$ into the monoid (with product given
by composition) of endomorphisms of a vector space $V$.

Given a subset $\mathcal S\subset M$ of a monoid, the submonoid 
$\mathcal M(\mathcal S)$ generated
by $\mathcal S$ is the smallest submonoid of $\mathcal M$ containing 
$\mathcal S$. A monoid $\mathcal M$ is finitely generated if $\mathcal
M=\mathcal M(\mathcal G)$ for some finite subset $\mathcal G\subset
\mathcal M$, called a {\it generating set}.

The {\it free monoid} $\mathcal M_{\mathcal A}$
over a set (also called alphabet)
$\mathcal A$ is the set of all finite words with letters in 
$\mathcal A$. We call $\mathcal A$ the free generating set of
$\mathcal A$. Two free monoids over $\mathcal A$, respectively
$\tilde{\mathcal A}$, are isomorphic if and only if $\mathcal A$ and
$\tilde{\mathcal A}$ are equipotent. In particular, there exists,
up to isomorphism, a unique free monoid whose free generating set
has a given cardinal number. One can thus speak of the free monoid on
$n$ letters for a natural integer $n\in\mathbb N$.

There is a natural notion of a monoid presented by generators and 
relations. A free monoid has no relations and an arbitrary 
finitely generated monoid ${\mathcal Q}$ can always be given in the form
$${\mathcal Q}=\langle {\mathcal G}:{\mathcal R}\rangle$$
where ${\mathcal G}$ is a finite set of generators and ${\mathcal R}
\subset {\mathcal G}^*\times {\mathcal G}^*$ a (perhaps infinite)
set of relations of the form 
$L=R$ with $L,R\in{\mathcal G}^*$. 
The quotient monoid ${\mathcal Q}$ of the free monoid 
${\mathcal G}^*$ by the relations ${\mathcal R}$ is the set of 
equivalence classes of words in ${\mathcal G}^*$ by the equivalence 
relation generated by $UL_iV\sim UR_iV$ for $U,V\in {\mathcal G}^*$
and $(L_i,R_i)\in {\mathcal R}$.

\begin{rem} Given a quotient monoid $\mathcal Q=\langle \mathcal
  G:\mathcal R\rangle$, the free monoid $\mathcal F_\mathcal G$
on $\mathcal G$ surjects
onto $\mathcal Q$ . The set $\pi^{-1}(\tilde e)\subset \mathcal Q$
of preimages of the identity $\tilde e
\in\mathcal Q$ is of course a submonoid of $\mathcal F_\mathcal G$.

However, unlike in the case of groups, the kernel 
$\pi^{-1}(\tilde e)\subset \mathcal F_\mathcal G$ does not 
characterize $\mathcal Q$. An example is given by the monoid $
\mathcal Q=\{\tilde e,a\}$
with identity $\tilde e$ and product $aa=a$. The free monoid on $1$ generator
surjects onto $\mathcal Q$ but $\pi^{-1}(\tilde
e)=\{\emptyset\}\in\mathcal F_1$ reduces to the trivial submonoid 
in the free monoid $\mathcal F_1$ on $1$ generator. 
\end{rem}

\begin{rem} The composition law $\mathcal M\times \mathcal M
\longrightarrow \mathcal M$ of a monoid $\mathcal M$ endows
the free vector space generated by $\mathcal M$ with an associative
algebra-structure, denoted $\mathbf K[\mathcal M]$ and called the 
{\it monoid-algebra} of $\mathcal M$. The monoid algebra 
$\mathbf K[\mathcal M]$ of a free monoid is simply the
polynomial algebra on free generators of $\mathcal M$, considered
as non-commuting variables.

A monoid $\mathcal M$ has the {\it finite-factorisation property}
if every element $m\in \mathcal M$ has only a finite
number of distinct factorisations $m=m_1 m_2$ with $m_1,m_2\in\mathcal 
M$. If a monoid $\mathcal M$ has the finite-factorisation property
(eg. if $\mathcal M$ is a finitely generated quotient monoid 
such that $L_i$ and $R_i$ are of equal length for every relation 
$(L_i,M_i)\in\mathcal R$), then the algebraic dual $\mathbf
K[[\mathcal M]]$ (which consists
of all formal sums of elements in $\mathcal M$) is also an algebra 
(for the obvious product) and contains the monoid-algebra  
$\mathbf K[\mathcal M]$. In the case where $\mathcal M$ is a 
free monoid, the algebra $\mathbf K[[\mathcal M]]$ is the algebra
of formal power-series with unknowns the free generators of $\mathcal
M$, considered as non-commuting variables. 
\end{rem}

\section{The category $K^{\mathcal M}$}

Given two natural integers $p,q\in{\mathbb N}$ and a natural integer
$l\in{\mathbb N}$, 
we define ${\mathcal M}_{p\times
  q}^l$ as the set of all $p^lq^l$
pairs of words $(U,W)$ of common length $l$ with $U=u_1\dots
u_l\in\{0,\dots,p-1\}^l$ and $W=w_1\dots w_l\in\{0,\dots,q-1\}^l$. 
The set ${\mathcal
  M}_{p\times q}^0$ contains by convention only 
$(\emptyset,\emptyset)$. We denote by ${\mathcal
  M}_{p\times q}=\bigcup_{l\in{\mathbb N}} {\mathcal M}_{p\times q}^l$
the union of all finite sets ${\mathcal M}_{p\times q}^l$. 
The common length
$l=l(U)=l(W)\in{\mathbb N}$ is the {\it length} $l(U,W)$ of a
{\it word} $(U,W)\in{\mathcal M}_{p\times q}^l\subset
{\mathcal M}_{p\times q}$. We write
${\mathcal M}_{p\times q}^{\leq l}$ for the obvious
set of all $1+pq+\dots+p^lq^l=
\frac{(pq)^{l+1}-1}{pq-1}$ words of length at most $l$ 
in ${\mathcal M}_{p\times q}$. 

The concatenation $(U,W)(U',W')=(UU',WW')$ turns 
${\mathcal M}_{p\times q}=({\mathcal M}_{p\times q}^1)^*$ 
into a free monoid on the set $\mathcal M_{p\times q}^1$ of all 
$pq$ words with length $1$ in  ${\mathcal M}_{p\times q}$.

Equivalently, ${\mathcal M}_{p\times q}$ can be described as the
submonoid of the product-monoid ${\mathcal M}_p\times
{\mathcal M}_q$ (where ${\mathcal M}_r$ stands for the free monoid
on $r$ generators) whose elements $(U,W)$ are 
all pairs $U\in{\mathcal M}_p,W\in{\mathcal M}_q$ of the same length.

\begin{rem} \label{remmonoid}
Since the free cyclic monoid $\{0\}^*$ contains 
a unique word of length $l$ for every $l\in{\mathbb N}$, we have an
obvious isomorphism ${\mathcal M}_{p\times 1}\sim {\mathcal M}_p
=\{0,\dots,p-1\}^*$. Moreover, since the free monoid $\mathcal M_0$ 
on the empty alphabet is the trivial monoid on one element,
the monoid ${\mathcal M}_{p\times q}$
is reduced to the empty word $(\emptyset,\emptyset)$ if
$pq=0$. 
\end{rem}

\begin{rem} Many constructions involving the monoid ${\mathcal
    M}_{p\times q}$ have generalizations to 
${\mathcal M}_{p_1\times p_2\times\cdots \times p_k}$ 
defined in the obvious way. In particular, using 
Remark \ref{remmonoid}, we write often simply ${\mathcal M}_p$
instead of ${\mathcal M}_{p\times 1}$ or ${\mathcal M}_{1\times p}$.
\end{rem}

In the sequel, we consider a fixed commutative field ${\mathbf K}$
(many results continue to hold for commutative rings). 
We denote by ${\mathbf K}^{{\mathcal M}_{p\times q}}$ the vector space 
of all functions from ${\mathcal M}_{p\times q}$ into ${\mathbf K}$.
We denote by $A[\mathcal S]$ the restriction of $A\in{\mathbf K}^{{
\mathcal M}_{p\times q}}$ to a subset ${\mathcal S}\subset 
{\mathcal M}_{p\times q}$. In particular, we write 
$A[U,W]$ for the evaluation of $A$ on a word $(U,W)\in
{\mathcal M}_{p\times q}$. For $A\in{\mathbf K}^{{
\mathcal M}_{p\times q}}$ and $l\in {\mathbb N}$ we consider
the restriction $A[{\mathcal M}_{p\times q}^l]$ of $A$ to
${\mathcal M}_{p\times q}^l$ as a matrix of size $p^l\times q^l$
with coefficients $A[U,W]$ indexed by all words 
$(U,W)=(u_1\dots u_l,w_1\dots w_l)\in {\mathcal M}_{p\times q}^l$.

For $A\in {\mathbf K}^{{\mathcal M}_{p\times r}}$ and 
$B\in {\mathbf K}^{{\mathcal M}_{r\times q}}$, we define the {\it
matrix product} or {\it product} 
$AB\in {\mathbf K}^{{\mathcal M}_{p\times q}}$ 
of $A$ and $B$ by
$$(AB)[U,W]=\sum_{V\in\{0,\dots,r-1\}^l}A[U,V]B[V,W]$$
for $(U,W)\in{\mathcal M}_{p\times q}^l$ a word of length $l$.
The matrix product is obviously bilinear
and associative. We get thus a category ${\mathbf K}^{\mathcal M}$
(see for instance \cite{SaMaLa} for definitions)
as follows: An {\it object} of ${\mathbf K}^{\mathcal M}$ is a vector space
of the form ${\mathbf K}^{{\mathcal M}_p}={\mathbf K}^{{\mathcal M}_{
p\times 1}}$ for $p\in{\mathbb N}$. A {\it morphism} (or {\it arrow})
is given by $A\in{\mathbf K}^{{\mathcal M}_{q\times p}}$ 
and defines a linear application from  ${\mathbf K}^{{\mathcal M}_p}$ 
to ${\mathbf  K}^{{\mathcal M}_q}$
by matrix-multiplication.

The matrix-product turns the set ${\mathbf K}^{{\mathcal M}_{p\times p}}$
of endomorphisms of an object ${\mathbf K}^{{\mathcal M}_p}$
into an algebra.

\begin{rem} \label{variationKM}
The category ${\mathbf K}^{{\mathcal M}}$ has also the
following slightly different realization: Associate to an object
${\mathbf K}^{{\mathcal M}_p}$
corresponding to the natural integer $p\in{\mathbb N}$ 
the ${\mathbb N}-$graded vector space $\mathcal{FS}_p
=\bigoplus_{l=0}^\infty
{\mathbf K}^{p^l}$, identified with the subspace of 
${\mathbf K}^{{\mathcal M}_p}$ of all functions with finite support.  
Morphisms are linear applications $\mathcal{FS}_p\longrightarrow 
\mathcal{FS}_q$ preserving the grading (and are given by a product 
$\prod_{l=0}^\infty {\mathbf K}^{q^l\times p^l}$
of linear maps
${\mathbf K}^{p^l}\longrightarrow {\mathbf K}^{q^l}$).
\end{rem}

\begin{rem} The vector-spaces $\mathcal{FS}_{p\times q}\subset
{\mathbf K}^{{\mathcal M}_{p\times q}}$
can be identified with the vector-spaces 
${\mathbf K}[X_{0,0},\dots,X_{p-1,q-1}]\subset
{\mathbf K}[[X_{0,0},\dots,X_{p-1,q-1}]]$ of polynomials and formal
power-series in $pq$ non-commuting variables 
$X_{u,w}, (u,w)\in\mathcal M_{p\times q}^1$. The vector-space
${\mathbf K}^{{\mathcal M}_{p\times q}}$ can also be considered
as the algebraic dual of $\mathcal{FS}_{p\times q}$.
\end{rem}

\begin{rem} A vector $X\in{\mathbf K}^{{\mathcal M}_{p\times q}}$ 
can be given
as a projective limit by considering the projection
$${\mathbf K}^{{\mathcal M}_{p\times q}^{\leq l+1}}\longrightarrow
 {\mathbf K}^{{\mathcal M}_{p\times q}^{\leq l}}$$
obtained by restricting the function $X[{\mathcal M}_{p\times q}^{\leq
  l+1}]$
to the subset ${\mathcal M}_{p\times q}^{\leq l}\subset
{\mathcal M}_{p\times q}^{\leq l+1}$.
\end{rem}

\begin{rem} The following analogue of the
tensor product yields a natural functor of
the category ${\mathbf K}^{\mathcal M}$: For $A\in{\mathbf
  K}^{{\mathcal M}_{p\times q}}$ and 
 $B\in{\mathbf
  K}^{{\mathcal M}_{p'\times q'}}$, define
$A\otimes B\in{\mathbf K}^{{\mathcal M}_{pp'\times qq'}}$
in the obvious way by considering 
tensor products $(A\otimes B)[{\mathcal M}_{pp'\times qq'}^l]=
A[{\mathcal M}_{p\times q}^l]\otimes 
B[{\mathcal M}_{p'\times q'}^l]$ of the graded parts.
This ``tensor product'' is bilinear and natural with respect to most 
constructions of this paper. The main difference with the 
usual tensor product is however the fact that
$A\otimes B$ can be zero even if $A$ and $B$ are both non-zero:
Consider $A$ and $B$ such that $A[U,W]=0$ if $(U,W)$ is of 
even length and $B[U',W']=0$ if $(U',W')$ is of odd length.
\end{rem}

\begin{rem} The category ${\mathbf K}^{\mathcal M}$,
realized as in Remark \ref{variationKM} can be embedded as
a full subcategory into a larger category with objects
given by graded vector spaces
$\mathcal{FS}_{(d_0,d_1,\dots)}=\bigoplus_{l=0}^\infty {\mathbf K}^{d_l}$ 
indexed by arbitrary sequences $d=(d_0,d_1,d_2,\dots)
\in{\mathbb N}^{\mathbb N}$. Morphisms
from $\mathcal{FS}_{(d_0,d_1,\dots)}$ to ${\mathcal
  FS}_{(e_0,e_1,\dots)}$ are linear applications preserving the grading
and correspond to elements in the direct product
of matrices $\prod_{l=0}^\infty {\mathbf K}^{{e}_l\times d_l}$. 
The category ${\mathbf K}^{\mathcal M}$ corresponds to the 
full subcategory with objects indexed by
geometric progressions.
\end{rem}
\section{The category $\hbox{Rec}({\mathbf K})$}
\label{sectcatrec}

A word $(S,T)\in{\mathcal M}_{p\times q}$ defines an endomorphism
$\rho(S,T)\in\hbox{End}(\mathbf K^{\mathcal M_{p\times q}})$ 
of the vector space
${\mathbf K}^{{\mathcal M}_{p\times q}}$ 
by setting
$$(\rho(S,T)A)[U,W]=A[US,WT]$$ 
for $A\in {\mathbf K}^{{\mathcal M}_{p\times q}}$ and
$(U,W)\in {\mathcal M}_{p\times q}$. The easy computation
$$\begin{array}{l}
\displaystyle \rho(S,T)\big(\rho(S',T')A)[U,W]=\rho(S',T')A[US,WT]\\
\displaystyle \qquad =A[USS',WTT']
=\rho(SS',TT')A[U,W]\end{array}$$
shows that $\rho:\mathcal M_{p\times
  q}\longrightarrow\hbox{End}(\mathbf K^{\mathcal M_{p\times q}})$
is a linear representation from the free monoid
${\mathcal M}_{p\times q}$ into the monoid $\hbox{End}
({\mathbf K}^{{\mathcal M}_{p\times q}})$ of all linear
endomorphisms of
${\mathbf K}^{{\mathcal M}_{p\times q}}$. 
\begin{defin} We call the
monoid $\rho({\mathcal M}_{p\times q})\subset\hbox{End}
({\mathbf K}^{{\mathcal M}_{p\times q}})$
the {\it shift-monoid}. An element $\rho(S,T)\in \rho({\mathcal
  M}_{p\times q})$ is a {\it shift-map}.
\end{defin}

\begin{rem} The terminology is motivated by the special 
case $p=q=1$: 
An element $A\in{\mathbf K}^{{\mathcal M}_{1\times 1}}$ is completely 
described by the sequence $\alpha_0,\alpha_1,\dots\in {
\mathbf K}^{\mathbb N}$ defined by the evaluation
$\alpha_l=A[0^l,0^l]$ 
on the unique word $(0^l,0^l)\in{\mathcal M}_{1\times 1}^l$ 
of length $l$. The generator $\rho(0,0)\in \rho({
\mathcal M}_{1\times 1})$ acts on $A$ by the usual shift
$(\alpha_0,\alpha_1,\alpha_2,\dots)\longmapsto
(\alpha_1,\alpha_2,\alpha_3,\dots)$ which erases the first element 
$\alpha_0$ of the sequence $\alpha_0,\alpha_1,\dots$ 
corresponding to $A$.
\end{rem}

\begin{prop} The linear representation $\rho:\mathcal M_{p\times q}
\longrightarrow\hbox{End}
({\mathbf K}^{{\mathcal M}_{p\times q}})$ is faithful. The
shift-monoid $\rho(\mathcal M_{p\times q})\subset\hbox{End}
({\mathbf K}^{{\mathcal M}_{p\times q}})$ is thus isomorphic to 
the free monoid ${\mathcal M}_{p\times q}$. 
\end{prop}

{\bf Proof} Otherwise there exists $(S,T)\not= (S',T')\in\mathcal
M_{p\times q}$ such that $\rho(S,T)=\rho(S',T')$. Consider
$A\in\mathbf K^{\mathcal M_{p\times q}}$ such that $A[S,T]=1$ and
$A[U,W]=0$ otherwise. The inequality
$$(\rho(S,T)A)[\emptyset,\emptyset]=1\not=0=(\rho(S',T')A)
[\emptyset,\emptyset]$$
yields then a contradiction. \hfill$\Box$

Given subsets ${\mathcal S}\subset
{\mathcal M}_{p\times q}$ and ${\mathcal A}\subset 
{\mathbf K}^{{\mathcal M}_{p\times q}}$, we write
$$\rho({\mathcal S}){\mathcal A}=\{\rho(S,T)A\ \vert\  
(S,T)\in {\mathcal S},A\in {\mathcal A}\}  \subset {\mathbf K}^{{\mathcal 
M}_{p\times q}}.$$
Since $\rho(\emptyset,\emptyset)$ acts as the identity,
we have ${\mathcal A}\subset
\rho({\mathcal S}){\mathcal A}$ for any subset
${\mathcal S}\subset {\mathcal M}_{p\times q}$ containing
the empty word $(\emptyset,\emptyset)$ of length $0$.
\begin{defin}
A subset ${\mathcal X}\subset {\mathbf K}^{{\mathcal
      M}_{p\times q}}$ is a {\it recursively closed set} if
$\rho({\mathcal M}_{p\times q}){\mathcal X}={\mathcal X}$.
The {\it recursive set-closure} of ${\mathcal X}\subset
{\mathbf K}^{{\mathcal M}_{p\times q}}$ is given by 
$\rho({\mathcal M}_{p\times q}){\mathcal X}$ and is the smallest
recursively closed subset of ${\mathbf K}^{{\mathcal M}_{p\times q}}$
containing ${\mathcal X}$. The {\it recursive closure}
$\overline{\mathcal X}^{rec}$ is the linear span of 
$\rho({\mathcal M}_{p\times q}){\mathcal X}$ and is the
smallest linear subspace of ${\mathbf K}^{{\mathcal M}_{p\times q}}$
which contains ${\mathcal X}$ and is recursively closed. 
\end{defin}
By definition, a subspace ${\mathcal A}$ is
recursively closed if and only if ${\mathcal A}=\overline{\mathcal
  A}^{rec}$. Intersections and sums (unions) of recursively closed
subspaces (subsets) in ${\mathbf K}^{{\mathcal M}_{p\times q}}$
are recursively closed (subsets).

Any recursively closed subspace ${\mathcal A}\subset {\mathbf
  K}^{{\mathcal M}_{p\times q}}$ is invariant under the shift-monoid
and we call the restriction 
$\rho_{\mathcal A}({\mathcal M}_{p\times q})\in\hbox{End}
({\mathcal A})$ of $\rho({\mathcal M}_{p\times q})$ 
to the invariant subspace ${\mathcal A}$
the {\it shift-monoid} of ${\mathcal A}$.

\begin{defin} Given an element $A\in
{\mathbf K}^{{\mathcal M}_{p\times q}}$ with
recursive closure $\overline{A}^{rec}$, we call the dimension
$\hbox{dim}(\overline{A}^{rec})\in{\mathbb N}\cup \{\infty\}$
the {\it (recursive) complexity} of $A$. A {\it recurrence matrix} 
is an element of the vector space 
$$\hbox{Rec}_{p\times q}({\mathbf K})=\{A\in {\mathbf K}^{{\mathcal M}_{p\times q}}
\ \vert\ \hbox{dim}(\overline{A}^{rec})<\infty\}$$
consisting of all elements having finite complexity.
\end{defin}
It is easy to check that $\hbox{Rec}_{p\times q}({\mathbf K})$ is
a recursively closed subspace of ${\mathbf K}^{{\mathcal M}_{p\times
    q}}$ containing the subspace
$\mathcal{FS}_{p\times q}(\mathbf K)\subset \mathbf K^{\mathcal
  M_{p\times q}}$ consisting of all elements with finite support. 
An element $A\in{\mathbf K}^{{\mathcal M}_{p\times q}}$ is 
a recurrence matrix, if and only if the shift-monoid 
$\rho_{\overline{A}^{rec}}
({\mathcal M}_{p\times q})$ of ${\overline A}^{rec}$ 
is a finite-dimensional linear representation of 
the free monoid ${\mathcal M}_{p\times q}$.

\begin{prop} \label{propproduct} We have
$$\hbox{dim}(\overline{AB}^{rec})\leq 
\hbox{dim}(\overline{A}^{rec})\ \hbox{dim}(\overline{B}^{rec})$$
for the matrix-product $AB\in{\mathbf K}^{{
\mathcal M}_{p\times q}}$ of $A\in{\mathbf K}^{{\mathcal M}_{p\times
r}}$ 
and $B\in{\mathbf K}^{{\mathcal M}_{r\times q}}$.
\end{prop}

\begin{cor} \label{corproduct}
The matrix-product $AB\in{\mathbf K}^{{
\mathcal M}_{p\times q}}$ of two recurrence matrices
$A\in\hbox{Rec}_{p\times r}({\mathbf K}),B\in 
\hbox{Rec}_{r\times q}({\mathbf K})$, is a recurrence matrix.
\end{cor}

Corollary \ref{corproduct} suggests the following definition.
\begin{defin} The {\it category $\hbox{Rec}({\mathbf K})$
of recurrence matrices} is the subcategory of ${\mathbf K}^{\mathcal M}$
containing only recurrence matrices as arrows. 
Its objects can be restricted to $\hbox{Rec}_{p\times 1}
({\mathbf K})$ or even to the recursively closed subspaces
$\mathcal{FS}_p=\bigoplus_{l=0}^\infty {\mathbf K}^{p^l}$ of
elements with finite support.
\end{defin}

{\bf Proof of Proposition \ref{propproduct}} Given bases
$A_1,A_2,\dots$ of
$\overline{A}^{rec}\subset \hbox{Rec}_{p\times r}(\mathbf K)$ and  
$B_1,B_2,\dots$ of
$\overline{B}^{rec}\subset \hbox{Rec}_{r\times q}(\mathbf K)$,
the computation
$$\begin{array}{l}
\displaystyle (\rho(s,t)(A_i B_j))[U,W]=(A_i B_j)[Us,Wt]\\
\displaystyle \quad =\sum_{v=0}^{r-1}\sum_{V\in \{0,\dots,r-1\}^l}
A_i[Us,Vv]B_j[Vv,Wt]\\
\displaystyle \quad =\sum_{v=0}^{r-1}\sum_{V\in \{0,\dots,r-1\}^l}
(\rho(s,v)A_i)[U,V](\rho(v,t)B_j)[V,W]\\ \displaystyle\quad
=\left(\sum_{v=0}^{r-1}(\rho(s,v)A_i)(\rho(v,t)B_j)
\right)[U,W]\end{array}$$
(with $(U,W)\in{\mathcal M}_{p\times q}^l$) shows that
${\mathcal C}=\sum_{1\leq i,j}{\mathbf K} A_iB_j$ is recursively closed in 
${\mathbf K}^{{\mathcal M}_{p\times p}}$. The obvious inclusion
$\overline{A B}^{rec}\subset{\mathcal C}$ finishes the proof.
\hfill $\Box$

\subsection{Other ring-structures on $\hbox{Rec}_{p\times
    q}({\mathbf K})$}

The vector-space ${\mathbf
  K}^{{\mathcal M}_{p\times q}}$ carries two 
natural ring-structures which are both inherited by
$\hbox{Rec}_{p\times q}({\mathbf K})$.

A first ring-structure on ${\mathbf
  K}^{{\mathcal M}_{p\times q}}$ comes from the usual commutative
product of functions $(A\circ B)[U,W]=(A[U,W])(B[U,W])$.

A second product (non-commutative if $pq>1$) is given by 
associating to 
$A\in{\mathbf K}^{{\mathcal M}_{p\times q} }$ the formal power series
$$\sum_{(U,W)=(u_1\dots u_n,w_1,w_n)\in{\mathcal M}_{p\times q}}
A[u_1\dots u_n,w_1\dots w_n]\ X_{u_1,w_1}\cdots X_{u_n,w_n}\in
\mathbf K[[\mathcal M_{p\times q}]]$$
in $pq$ non-commuting variables $X_{u,w},(u,w)\in\mathcal M_{p\times q}^1$.
Multiplication of non-commutative formal power-series
endows ${\mathbf K}^{{\mathcal M}_{p\times q}}$ with the {\it
  convolution-product}
$$(A*B)[U,W]=\sum_{(U,W)=(U_1,W_1)(U_2,W_2)}A[U_1,W_1]B[U_2,W_2]$$
where the sum is over all $l+1$ factorizations $(U,W)=
(U_1,W_1)(U_2,W_2)$ of a word 
$(U,W)\in{\mathcal M}_{p\times q}^l$.

The following result shows that both ring-structures restrict to
$\hbox{Rec}_{p\times q}({\mathbf K})$.

\begin{prop} \label{otherproducts}
(i) $\hbox{Rec}_{p\times q}({\mathbf K})$ is a
  commutative ring for the ordinary product of functions (given by) 
$(A\circ B)[U,W]=A[U,W]B[U,W]$. More precisely, we have
$\overline{A\circ B}^{rec}\subset \sum {\mathbf K}A_i\circ B_j$
where $\overline{A}^{rec}=\sum {\mathbf K}A_i$ and
$\overline{B}^{rec}=\sum {\mathbf K}B_j$.

\ \ (ii) $\hbox{Rec}_{p\times q}({\mathbf K})$ is a
ring for the convolution-product $(A*B)[U,W]=\sum_{(U,W)=(U_1,W_1)
(U_2,W_2)}A[U_1,W_1]B[U_2,W_2]$. More precisely, we have
$\overline{A*B}^{rec}\subset \sum {\mathbf K}A_i+\sum{\mathbf K}A*B_j$
where $\overline{A}^{rec}=\sum {\mathbf K}A_i$ and
$\overline{B}^{rec}=\sum {\mathbf K}B_j$.
\end{prop}

{\bf Proof} Assertion (i) follows from Corollary \ref{corproduct}
and the existence of a diagonal embedding of 
${\mathbf K}^{{\mathcal M}_{p\times q}}$ into 
${\mathbf K}^{{\mathcal M}_{(pq)\times (pq)}}$ preserving the
complexity.

The computation
$$\begin{array}{l}
\displaystyle \big(\rho(s,t)(A*B)\big)[U,W]=(A*B)[Us,Wt]\\
\displaystyle =\sum_{(U,W)=(U_1,W_1)(U_2,W_2)}A[U_1,W_1]
B[U_2s,W_2t]+A[Us,Wt]B[\emptyset,\emptyset]\\
\displaystyle 
=A*(\rho(s,t)B)[U,W]+\rho(s,t)A[U,W]B[\emptyset,\emptyset]
\end{array}$$
shows the identity
$$\rho(s,t)(A*B)=A*(\rho(s,t)B)+B[\emptyset,\emptyset]\rho(s,t)A
$$
and implies assertion (ii).\hfill$\Box$

\subsection{Convergent elements in ${\mathbf K}^{{\mathcal M}_{p\times
      q}}$ and $\hbox{Rec}_{p\times q}({\mathbf K})$}

Using the bijection 
$$s_1\dots s_n\longmapsto \sum_{j=1}^n s_jp^{j-1}$$
between $\{0,\dots, p-1\}^n$ and $\{0,\dots,p^n-1\}$, an infinite
matrix $\tilde A$ with coefficients $\tilde 
A_{s,t},0\leq s,t\in{\mathbb N}$
gives rise to an element $A\in{\mathbf K}^{{\mathcal M}_{
p\times q}}$ by setting
$$A[s_1\dots s_n,t_1\dots t_n]=\tilde A_{s,t}\hbox{ where }
s=\sum_{j=1}^n s_jp^{j-1},t=\sum_{j=1}^n t_jq^{j-1}.$$
We call such an element $A$ {\it convergent} with 
limit (the infinite matrix) $\tilde A$. If $p=1$ or $q=1$, the 
limit matrix $\tilde A$ degenerates to a row or column-vector and it 
degenerates to a single coefficient $\tilde A_{0,0}$ if $pq=0$.
Obviously, an element $A\in{\mathbf K}^{{\mathcal M}_{p\times
    q}}$ is convergent if and only if $\rho(0,0)A=A$.
The vector space spanned by all converging elements in 
${\mathbf K}^{{\mathcal M}_{p\times p}}$ (respectively in 
$\hbox{Rec}_{p\times p}({\mathbf K})$) is not preserved under
matrix-products (except if $p\leq 1$) but contains a few 
subalgebras given for instance
by converging matrices with limit an infinite lower (or upper)
triangular matrix.

\begin{rem} The vector-space of convergent elements contains (strictly)
a unique maximal subspace which is recursively closed.
This subspace is the set of all convergent elements $A\in 
\hbox{Rec}_{p\times q}({\mathbf K})$ such that $\overline{A}^{rec}$
consists only of convergent elements.
\end{rem}

\section{The quotient category $\hbox{Rec}({\mathbf K})/{\mathcal  FS}$ modulo elements of finite support and stable
complexity}\label{finsuppquotient}

We call a quotient ${\mathcal A}/
{\mathcal B}$ with vector-spaces ${\mathcal B}\subset {\mathcal A}\subset 
\hbox{Rec}_{p\times q}({\mathbf K})$ a 
{\it quotient-space of recurrence matrices} if
${\mathcal A}$ and $ {\mathcal B}$ are both recursively closed.

An important example is given by ${\mathcal A}=\hbox{Rec}_{
p\times q}({\mathbf K})$ and ${\mathcal B}=\mathcal{FS}_{p\times q}$ 
the vector space of all functions $B\in{\mathbf K}^{{\mathcal M}_{
p\times q}}$ with finite support. The vector space 
$\mathcal{FS}_{p\times q}$ is recursively closed since it consists
of all elements $B$ such that the shift map $\rho({\mathcal
  M}_{p\times q})$ acts in a nilpotent way on $\overline{B}^{rec}$
(more precisely, for every element $B\in \mathcal{FS}_{p\times q}$ 
there exists a natural integer $N$ such that 
$\rho(U,W)B=0$ if $(U,W)\in{\mathcal M}_{p\times q}$ is of length 
$\geq N$). Elements of $\mathcal{FS}_{p\times q}$ are in some sense 
trivial. Since the subset
$\mathcal{FS}=\cup_{p,q\in{\mathbb N}} \mathcal{FS}_{p\times q}
\subset \cup_{p,q\in{\mathbb N}}\hbox{Rec}_{p\times q}({\mathbf K})$ 
is in an obvious sense a ``two-sided ideal'' for the matrix-product
(whenever defined), 
we get a functor from $\hbox{Rec}({\mathbf K})$ onto
the ``quotient-category'' $\hbox{Rec}({\mathbf K})/
\mathcal{FS}$ by considering the projection $\pi_{\mathcal{FS}}:
\hbox{Rec}(\mathbf K)\longrightarrow \hbox{Rec}(\mathbf
K)/\mathcal{FS}$. 

Given a recursively closed finite-dimensional subspace
$\mathcal A\subset\hbox{Rec}_{p\times q}(\mathbf K)$,
we define its {\it stable complexity}  
$\hbox{dim}_{\mathcal{FS}}(\mathcal A)$ as the dimension
$\hbox{dim}(\pi_{\mathcal{FS}}(\mathcal A))\leq \hbox{dim}(\mathcal A)$ 
of its projection onto the quotient space 
$\hbox{Rec}_{p\times q}(\mathbf K)/
\mathcal{FS}_{p\times q}$. The {\it stable complexity} of an element
$A\in\hbox{Rec}_{p\times q}(\mathbf K)$ is the stable complexity
$\hbox{dim}(\pi_{\mathcal{FS}}(\overline{A}^{rec}))$
of its recursive closure. The {\it complexity} of an element
$\tilde A\in\hbox{Rec}_{p\times q}(\mathbf K)/\mathcal{FS}_{p\times
    q}$is defined as the stable complexity $\hbox{dim}_{\mathcal{FS}}
(A)$ of any lift $A\in\pi_{\mathcal{FS}}^{-1}(\tilde A)\subset
\hbox{Rec}_{p\times q}(\mathbf K)$.

\begin{prop} \label{propproductstable} (i) We have
$$\hbox{dim}_{\mathcal{FS}}(\overline{AB}^{rec})\leq 
\hbox{dim}_{\mathcal{FS}}(\overline{A}^{rec})\ 
\hbox{dim}_{\mathcal{FS}}(\overline{B}^{rec})$$
for the matrix-product $AB\in{\mathbf K}^{{
\mathcal M}_{p\times q}}$ of $A\in{\mathbf K}^{{\mathcal M}_{p\times
r}}$ 
and $B\in{\mathbf K}^{{\mathcal M}_{r\times q}}$.

\ \ (ii) We have 
$$\hbox{dim}_{\mathcal{FS}}(\overline{\pi_{\mathcal{FS}}^{-1}(\tilde{AB}
)}^{rec})\leq 
\hbox{dim}_{\mathcal{FS}}(\overline{\pi_{\mathcal{FS}}^{-1}(\tilde A)}^{rec})\ 
\hbox{dim}_{\mathcal{FS}}(\overline{\pi_{\mathcal{FS}}^{-1}(\tilde B)}^{rec})$$
for the product $\tilde{AB}\in{\mathbf K}^{{
\mathcal M}_{p\times q}}/\mathcal{FS}_{p\times q}$ of $\tilde A
\in{\mathbf K}^{{\mathcal M}_{p\times
r}}/\mathcal{FS}_{p\times r}$ 
and $\tilde B\in{\mathbf K}^{{\mathcal M}_{r\times q}}/\mathcal{FS}_{
r\times q}$.
\end{prop}

The proof is the same as for Proposition \ref{propproduct}.

\begin{rem} The subspace $\mathcal{FS}_{p\times q}=\hbox{Rec}_{p\times
    q}({\mathbf K})\cap \mathcal{FS}$ of
all elements with finite support in $\hbox{Rec}_{p\times q}({
\mathbf K})$ is also an ideal for the commutative product
of functions. Although $\mathcal{FS}_{p\times q}$
is also closed for the convolution-product (and
corresponds to polynomials in non-commutative variables)
it is not an ideal in the convolution-ring since
it
contains the convolutional identity
$Id_*$ given by $Id_*[\emptyset,\emptyset]=1$ and $
Id_*[U,W]=0$ for $(U,W)\in{\mathcal M}_{p\times q}\setminus (
\emptyset,\emptyset)$.
\end{rem}

\begin{rem} \label{remquotientKMFS}
One checks easily that the set $\mathcal{FS}$
is also an ideal (in the obvious sense) of the category 
$\mathbf K^{\mathcal M}$. One can thus consider the quotient
algebras $\mathbf K^{\mathcal M_{p\times p}}/\mathcal{FS}$ and the 
functor (still denoted) $\pi_{\mathcal{FS}}:\mathbf K^{\mathcal
  M}\longrightarrow \mathbf K^{\mathcal M}/\mathcal{FS}$
onto the quotient category $\mathbf K^{\mathcal M}/\mathcal{FS}$.
\end{rem}  


\section{Matrix algebras}\label{sectmatralg}

Given $A\in\mathbf K^{\mathcal M_{p\times p}}$, the elements
$\rho(S,T)A$ indexed by $(S,T)\in\mathcal M_{p\times p}^l$ 
can be considered as
coefficients of a square matrix of size $p^l\times p^l$ with values
in the ring $\mathbf K^{\mathcal M_{p\times p}}$. More precisely,
the application 
$$A\longmapsto \varphi^l(A)=\left(\rho(S,T)A\right)_{(S,T)\in\mathcal 
M_{p\times p}^l}$$
can be described as the restriction
$$A=\prod_{j=0}^\infty A[\mathcal M_{p\times p}^j]\longmapsto
\varphi^l(A)=\prod_{j=l}^\infty A[[\mathcal M_{p\times p}^j]$$
obtained by removing the initial terms $A[\mathcal M_{p\times p}^0],
A[\mathcal M_{p\times  p}^1],\dots,A[\mathcal M_{p\times p}^{l-1}]$ 
from the sequence of matrices
$A[\mathcal M_{p\times p}^0],A[\mathcal M_{p\times
  p}^1],\dots$ representing $A$.

We leave the proof of the following obvious assertions to the reader.

\begin{prop} (i) We have $\varphi^{l+k}=\varphi^{l}\circ \varphi^k$.

\ \ (ii) $\varphi^l$ defines a morphism of rings between the ring
$\mathbf K^{\mathcal M_{p\times p}}$ and the ring of $p^l\times p^l $
matrices with values in $\mathbf K^{\mathcal M_{p\times p}}$.

\ \ (iii) $\varphi^l$ restricts to a morphism of rings between the ring
$\hbox{Rec}_{p\times p}(\mathbf K)$ and the ring of $p^l\times p^l $
matrices with values in $\hbox{Rec}_{p\times p}(\mathbf K)$.

\ \ (iv) We have $\mathcal{FS}_{p\times p}=\cup_{l=0}^\infty
\ker(\varphi^l)\subset \hbox{Rec}_{p\times p}(\mathbf K)\subset
\mathbf K^{\mathcal M_{p\times p}}$ for the vector space 
$\mathcal{FS}_{p\times p}$ of elements with finite support in 
$\hbox{Rec}_{p\times p}(\mathbf K)$ (see section
\ref{finsuppquotient}).

Moreover, $\varphi^l$ induces an injective morphism
$\overline\varphi^l$ from the quotient ring $\mathbf K^{\mathcal
  M_{p\times p}}/\mathcal{FS}_{p\times p}$ (respectively $\hbox{Rec}_{p\times p}(\mathbf K)/\mathcal{FS}_{p\times p}$)
onto the ring of $p^l\times p^l$ matrices with values in 
$\mathbf K^{\mathcal  M_{p\times p}}/\mathcal{FS}_{p\times p}$ 
(respectively in
$\hbox{Rec}_{p\times p}(\mathbf K)/\mathcal{FS}_{p\times p}$)
(see Remark \ref{remquotientKMFS}
for the definition of the quotient ring $\mathbf K^{\mathcal
  M_{p\times p}}/\mathcal{FS}_{p\times p}$).
\end{prop}

\begin{rem} $\varphi^l(A)$ with $A\in \mathbf K^{
\mathcal M_{p\times q}}$ is defined for arbitrary $(p,q)\in \mathbb N^2$.
It can be considered as a functor from the category 
$\mathbf K^{\mathcal M}$ (or $\hbox{Rec}(\mathbf K)$)
into a category with objects indexed by 
$p\in\mathbb N$ and arrows given by matrices of size
$p^l\times q^l$ having coefficients in $\mathbf K^{\mathcal M_{p\times
    q}}$ (respectively in $\hbox{Rec}(\mathbf K)$).
\end{rem}


\section{Presentations}
We describe in this chapter two mehtods of defining 
elements of $\hbox{Rec}_{p\times q}(\mathbf K)$ by finite amounts
of data. The first method, {\it (monoidal) presentations}, puts
emphasis on the shift monoid. The second method, {\it 
recursive presentations}, is often more intuitive and sometimes
more concise.

\subsection{Monoidal presentations}

Let ${\mathcal A}=\bigoplus_{j=1}^a {\mathbf K}A_j\subset 
\hbox{Rec}_{p\times q}({\mathbf K})$ be a finite-dimensional 
recursively closed vector space with basis $A_1,\dots, A_a$. 
The isomorphism ${\mathbf K}^a\longrightarrow {\mathcal A}$
defined by $(\alpha_1,\dots,\alpha_a)\longmapsto
\sum_{j=1}^a \alpha_j A_j$ realizes the shift-monoid 
$\rho_{\mathcal A}({\mathcal M}_{p\times q})\subset {\mathbf K}^{
a\times a}$ of ${\mathcal A}$ as 
a matrix monoid in $\hbox{End}({\mathbf K}^a)$.
The generators $\rho_{\mathcal A}(s,t)\in {\mathbf K}^{a\times
  a},0\leq s<p,0\leq t<q$ 
(with respect to the fixed basis $A_1,\dots A_a$ of ${\mathcal A}$) are
called {\it shift-matrices}. Their coefficients $\rho_{\mathcal A}
(s,t)_{j,k}$ (for $1\leq j,k\leq a$) are given by
$$\rho(s,t)A_k=\sum_{j=1}^a\rho_{\mathcal A}(s,t)_{j,k} A_j\ .$$

More generally, we can define shift-matrices $\rho_{\mathcal A}(s,t)
\in {\mathbf K}^{d\times d}$ with respect to any
finite (not necessarily linearly independent) generating set 
$A_1,\dots ,A_d$ of a recursively closed finite-dimensional
vector-space ${\mathcal A}\subset \hbox{Rec}_{p\times q}({\mathbf K})$ 
by requiring the identity
$\rho(s,t)A_k=\sum_{j=1}^d\rho_{\mathcal A}(s,t)_{j,k} A_j$. These
equations define the matrices $\rho_{\mathcal A}(s,t)$ up to linear
applications ${\mathbf K}^d\longrightarrow K$
where $K=\{(\alpha_1,\dots, \alpha_d)\in{\mathbf K}^d\ \vert 
\sum_{j=1}^d \alpha_jA_j=0\}$ is the subspace of relations 
among the generators $A_1,\dots A_d\in{\mathcal A}$ (or, equivalently,
the kernel of the map $(\alpha_1,\dots,\alpha_d)\in{\mathbf
  K}^d\longmapsto \sum_{j=1}^d \alpha_jA_j$).

\begin{prop} \label{propcoeffs} Let ${\mathcal A}=
\sum_{j=1}^d {\mathbf K}A_j\subset \hbox{Rec}_{p\times q}({\mathbf
  K})$ be a recursively closed 
vector-space with shift-matrices $\rho_{\mathcal A}(s,t)\in{\mathbf
  K}^{d\times d}$ with respect to the (not necessarily free)
generating set $A_1,\dots A_d$. We have 
$$\left(\begin{array}{c}A_1[s_1\dots s_n,t_1\dots,t_n]\\
\vdots\\A_d[s_1\dots s_n,t_1\dots,t_n]\end{array}\right)=
\rho_{\mathcal A}(s_n,t_n)^t\cdots \rho_{\mathcal A}(s_1,t_1)^t
\left(\begin{array}{c}A_1[\emptyset,\emptyset]\\
\vdots\\ A_d[\emptyset,\emptyset]\end{array}\right).$$
\end{prop}

{\bf Proof} 
For $(U,W)=(s_1\dots s_n,t_1\dots t_n)\in{\mathcal M}_{p\times q}^n$,
the formula
$$A[U,W]=(\rho_{\mathcal A}(U,W)A)[\emptyset,\emptyset]=(\rho_{\mathcal A}(s_1,t_1)\cdots
\rho_{\mathcal A}(s_n,t_n)A)[\emptyset,\emptyset]$$
implies the result by duality.

A second proof is given by the computation
$$\left(\begin{array}{c}A_1[s_1\dots s_n,t_1\dots, t_n]\\
\vdots\\ A_d[s_1\dots s_n,t_1\dots, t_n]\end{array}\right)=
\left(\begin{array}{c}\rho_{\mathcal A}(s_n,t_n)A_1[s_1\dots s_{n-1},t_1\dots, 
t_{n-1}]\\
\vdots\\\rho_{\mathcal A}(s_n,t_n)A_d[s_1\dots s_{n-1},t_1\dots, 
t_{n-1}] \end{array}\right)$$
$$=\rho_{\mathcal A}(s_n,t_n)^t
\left(\begin{array}{c}A_1[s_1\dots s_{n-1},t_1\dots, 
t_{n-1}]\\
\vdots\\A_d[s_1\dots s_{n-1},t_1\dots, 
t_{n-1}] \end{array}\right)$$
and induction on $n$.\hfill $\Box$

\begin{defin} A {\it monoidal presentation} or {\it presentation} 
${\mathcal P}$ of complexity $d$ is given by the following data:

a vector $(\alpha_1,\dots,\alpha_d)\in
{\mathbf K}^d$ of $d$ {\it initial values},

$pq$ {\it shift-matrices} $\rho_{\mathcal P}(s,t)\in {\mathbf
  K}^{d\times d}$ of size $d\times d$ with coefficients
$\rho_{\mathcal P}(s,t)_{k,j},1\leq k,j\leq d$. 
\end{defin}

In the sequel, a presentation will always denote a monoidal presentation.

\begin{prop} A presentation ${\mathcal P}$
of complexity $d$ as above 
defines a unique set $A_1,\dots,A_d\in\hbox{Rec}_{p\times q}$
of $d$ recurrence matrices such that $A_k[\emptyset,\emptyset]=\alpha_k$
and $\rho(s,t)A_k=\sum_{j=1}^d \rho_{\mathcal P}(s,t)_{j,k}A_j$ for
$1\leq k\leq d$ and for all $(s,t)\in{\mathcal M}_{p\times q}^1$.
\end{prop}

{\bf Proof} For $(U,W)=(s_1\dots s_n,t_1\dots t_n)\in{\mathcal M}_{p\times
  q}^n$,
we define the evaluations $A_1[U,W],\dots A_d[U,W]\in{\mathbf K}$
by 
$$\left(\begin{array}{c}A_1[U,W]\\\vdots\\ A_d[U,W]\end{array}\right)=
\rho_{\mathcal P}(s_n,t_n)^t\cdots \rho_{\mathcal P}(s_1,t_1)^t  
\left(\begin{array}{c}A_1[\emptyset,\emptyset]\\\vdots\\ A_d[
\emptyset,\emptyset]\end{array}\right).$$
The result follows from Proposition \ref{propcoeffs}.\hfill $\Box$

In the sequel, we will often drop the subscript ${\mathcal P}$
for shift-matrices of a presentation. Thus, we identify 
(abusively) shift-matrices with the corresponding 
shift-maps, restricted to the subspace defined by the presentation.

A presentation ${\mathcal P}$ is {\it reduced} if the elements
$A_1,\dots,A_d\in\hbox{Rec}_{p\times q}({\mathbf K})$
defined by ${\mathcal P}$ are linearly independent. 
We say that a presentation ${\mathcal P}$ {\it presents}
(or is a presentation of) the recurrence matrix
$A=A_1\in\hbox{Rec}_{p\times q}({\mathbf K})$.
The empty presentation of complexity $0$ presents 
by convention the 
zero-element of $\hbox{Rec}_{p\times q}({\mathbf K})$. 
A presentation of $A\in \hbox{Rec}_{p\times q}({\mathbf K})$ 
with complexity $a=\hbox{dim}(\overline{A}^{rec})$
is {\it minimal}.
Every recurrence matrix $A\in\hbox{Rec}_{p\times q}({\mathbf K})$
has a minimal presentation ${\mathcal P}$:
Complete $0\not=A\in\hbox{Rec}_{p\times q}({\mathbf K})$ to a basis
$A_1=A,\dots,A_{a}$ of its recursive closure $\overline{A}^{rec}$.
For $1\leq k\leq a$, set $\alpha_k=A_k[\emptyset,\emptyset]
\in {\mathbf K}$ and define the shift matrices 
$\rho(s,t)_{\mathcal P}\in \hbox{End}({\mathbf K}^a)$ by
$$\rho(s,t)A_k=\sum_{j=1}^{a}\rho_{\mathcal P}(s,t)_{j,k}A_j\ .$$
Linear independency of $A_1,\dots,A_a$ implies that the shift
matrices $\rho(s,t)$ are well-defined.

In the sequel, a presentation denotes
often a finite
set of recurrence matrices $A_1,\dots,A_d\in\hbox{Rec}_{p\times q}
({\mathbf K})$ spanning a
recursively closed subspace $\sum_{k=1}^d {\mathbf K}A_k$
together with $pq$ suitable shift matrices in $\hbox{End}({\mathbf K}^d)$
(which are sometimes omitted if they are obvious).
We denote by $(A_1,\dots,A_d)[\emptyset,\emptyset]
\in {\mathbf K}^d$ the corresponding 
initial values.

\begin{prop} \label{proppressum}
Given presentations $\mathcal P_{\mathcal A},\mathcal P_{\mathcal B}$
with respect to generators
$A_1,\dots,A_d,B_1,\dots,B_e\in\hbox{Rec}_{p\times q}(\mathbf K)$, a
presentation of $\mathcal C=\mathcal A+\mathcal B$ with respect to 
the generators $A_1,\dots,A_d,B_,\dots,B_e$ is given by shift
matrices
$$\rho_{\mathcal C}(s,t)=\left(\begin{array}{cc}
\rho_{\mathcal A}(s,t)&0\\0&\rho_{\mathcal B}(s,t)\end{array}\right)
\in\mathbf K^{(d+e)\times d+e)}$$
consisting of diagonal blocs.
\end{prop}

The proof is obvious.

\begin{prop} \label{presentproduct}
Given presentations $\mathcal P_{\mathcal A},
\mathcal P_{\mathcal B}$ with generators $A_1,\dots,A_d\in
\hbox{Rec}_{p\times r}(\mathbf K),B_1,\dots,B_e\in
\hbox{Rec}_{r\times q}(\mathbf K)$, a presentation $\mathcal
P_{\mathcal C}$ of the recursively closed vector space 
$\mathcal C\subset \hbox{Rec}_{p\times
  q}(\mathbf K)$ with respect to the generators $C_{ij}=A_iB_j$
of all products among generators is given by the initial values 
$$C_{ij}[\emptyset,\emptyset]=A_i[\emptyset,\emptyset]B_j[
\emptyset,\emptyset]$$
and shift matrices $\rho_{\mathcal C}(s,t)\in\mathbf K^{de\times de}$ 
with coefficients $$\rho_{\mathcal C}(s,t)_{kl,ij}=
\sum_{u=1}^r \rho_{\mathcal A}(s,u)_{k,i}\rho_{\mathcal B}(u,t)_{l,j}
\ .$$
\end{prop}

{\bf Proof} There is nothing to prove for the initial values.

For the shift matrices, we have
$$\begin{array}{ll}
\displaystyle \rho_{\mathcal C}(s,t)C_{ij}&
\displaystyle =\rho(s,t)(A_iB_j)=
\sum_{u=1}^r \rho_{\mathcal A}(s,u)A_i\rho_{\mathcal B}(u,t)B_j\\
&\displaystyle =\sum_{u=1}^r\sum_{k=1}^d\sum_{l=1}^e
\rho_{\mathcal A}(s,u)_{k,i}\rho_{\mathcal B}(u,t)_{l,j}A_kB_l\\
&\displaystyle =\sum_{k=1}^d\sum_{l=1}^e\left(\sum_{u=1}^r
\rho_{\mathcal A}(s,u)_{k,i}\rho_{\mathcal B}(u,t)_{l,j}\right)
C_{kl}\end{array}$$
which ends the proof.\hfill $\Box$

\begin{rem} The presentation $\mathcal P_{\mathcal C}$ given by
proposition \ref{presentproduct} is in general not reduced,
even if $\mathcal P_{\mathcal A}$ and $\mathcal P_{\mathcal B}$
are reduced presentations. The reason for this are (possible) 
multiplicities (up to isomorphism) of submonoids in the shift monoid 
$\rho_{\mathcal C}(\mathcal M_{p\times q})$.
\end{rem}

\begin{prop} \label{propprestransposed}
Given a presentation  ${\mathcal P}$ 
of $d$ recurrence matrices
$A_1,\dots A_d\in\hbox{Rec}_{p\times q}({\mathbf K})$ spanning a
recursively closed subspace, a
presentation ${\mathcal P}^t$ of the transposed recurrence
matrices $A_1^t,\dots A_d^t\in\hbox{Rec}_{q\times p}(
{\mathbf K})$ is given by the same
initial data $(A_1^t,\dots A_d^t)[\emptyset,\emptyset]=
(A_1,\dots A_d)[\emptyset,\emptyset]$ and the same shift matrices 
$\tilde \rho(t,s)=\rho(s,t)$, up to ``transposition'' of their
labels.
\end{prop}

{\bf Proof} Use an easy induction on $l$ for
the restricted functions $A_1[{\mathcal M}_{p\times q}^l],
\dots,A_d[{\mathcal M}_{p\times q}^l]$.\hfill$\Box$

\begin{rem} Shift-matrices of 
transposed recurrence matrices are identical to the original
shift-matrices (and are not to be transposed), only their labels
are rearranged!
\end{rem}

\begin{rem} Endowing ${\mathcal M}_{p\times q}$ with a complete order,
the recursive closure $\overline{A}^{rec}$ of an element $A\in
{\mathbf K}^{{\mathcal M}_{p\times q}}$ has a unique basis of the form 
$\rho(U_1,W_1)A,\rho(U_2,W_2)A,\dots$ where the word $(U_j,W_j)\in 
{\mathcal M}_{p\times q}$ (if it exists) is inductively defined as
the smallest element such that  
$\rho(U_1,W_1)A,\dots,\rho(U_j,W_j)A $ are linearly independent.
\end{rem}

\begin{rem} The stable complexity
  $\hbox{dim}_{\mathcal{FS}}(\pi_{\mathcal{FS}}(\overline{A}^{rec}))$
introduced in chapter \ref{finsuppquotient}
of an element $A\in \hbox{Rec}_{p\times q}(\mathbf K)$ equals
$\hbox{dim}(\overline{A}^{rec})-\hbox{dim}(\overline{A}^{rec}
\cap\mathcal{FS})$
where $\overline{A}^{rec}\cap\mathcal{FS}$ is the 
maximal recursively closed subspace of
$\overline{A}^{rec}$ with nilpotent action of 
$\rho(\mathcal M_{p\times q})$. 
\end{rem}

\subsection{Recursive presentations}\label{subsrecpres}

We start by defining recursive symbols which are 
the key ingredients for recursive presentations.

Given a field (or ring) $\mathbf K$, fixed in the sequel, 
the set 
$$\mathcal{RS}_{p\times q}
(A_1,\dots,A_d)=\cup_{d=0}^\infty
\mathcal{RS}_{p\times q}^{\leq d}(A_1,\dots,A_d)$$ 
of recursive $p\times q-$symbols over $A_1,\dots,
A_d$ is recursively defined as follows:
$\mathcal{RS}_{p\times q}^{\leq -1}(A_1,\dots,A_d)=\emptyset$ and
$\mathcal{RS}_{p\times q}^{\leq d}(A_1,\dots,A_d)$ is
the set of symbols
$(\rho,R)$ where $\rho\in\mathbf K$ is a constant
and $R$ is a matrix of size $p\times q$ with coefficients
$R_{s,t},0\leq s<p,0\leq t<q$ in the free vector space spanned
by $A_1,\dots,A_d$ and elements of 
$\mathcal{RS}_{p\times q}^{\leq d-1}(A_1,\dots,A_d)$.

A symbol $(\rho,R)$ has {\it depth} $d(\rho,R)=d$ if 
$(\rho,R)\in \mathcal{RS}_{p\times q}^{d}(A_1,\dots,A_d)=
\mathcal{RS}_{p\times q}^{\leq d}(A_1,\dots,A_d)\setminus
\mathcal{RS}_{p\times q}^{\leq d-1}(A_1,\dots,A_d)$.

{\bf Examples}
An example of a recursive $2\times 2-$symbol of depth $2$
over $A,B$ (with groundfield $\mathbb Q$) is for instance given by
$$(2,\left(\begin{array}{cc}
3B-2A&-A+(-1,\left(\begin{array}{cc}B&A+B\\-A&2A-B\end{array}\right))\\
5A-B+(0,
\left(\begin{array}{cc}A&B\\-B&(\rho,R)\end{array}\right))&
5A-3B\end{array}\right))$$
where $(\rho,R)=(7,\left(\begin{array}{cc}
    2A-B&3B\\-A+B&5A+B\end{array}\right))\in \mathcal{RS}_{2
\times 2}^0(A,B)$.

The expresssion $R=(1,(\rho,-A+R))$ defines however no recursive 
$1\times 1-$symbol over $A$.

{\bf Definition} A {\it recursive presentation} for $A_1,\dots,
A_d$ is given by identities
$$A_1=(\rho_1,R_1),\ A_2=(\rho_2, R_2),\dots,A_d=(\rho_d, R_d)$$
where $(\rho_1,R(1)),\dots,(\rho_d,R(d))\in \mathcal{RS}_{p\times q}
(A_1,\dots,A_d)$
are recursive $p\times q-$symbols over $A_1,\dots,A_d$.

With a hopefully understandable abuse of notation we say that
$B,A_1,\dots,A_d\in\mathbf K^{\mathcal M_{p\times q}}$ is a solution
to the equation $B=(\rho,R)$ with 
$(\rho,R)$ a $p\times q-$symbol over $A_1,\dots,A_d$
(here is the abuse) if $B$ satisfies $B[\emptyset,\emptyset]=\rho$ and
we have recursively the identities $\rho(s,t)B=R_{s,t}$ for 
$0\leq s<p,0\leq t<q$.

\begin{prop} \label{proprecdef} The identities of a 
recursive presentation 
for $A_1,\dots,A_d$ have a unique common solution $A_1,\dots,A_d
\in\mathbf K^{\mathcal M_{p\times q}}$. 

Moreover, the vector space spanned by
$$\rho(\mathcal M_{p\times q}^{\leq d(\rho_1,R_1)})A_1,\dots,
 \rho(\mathcal M_{p\times q}^{\leq d(\rho_d,R_d)})A_d$$
is recursively closed and we have thus 
$A_1,\dots,A_d\in\hbox{Rec}_{p\times q}(\mathbf K)$.
\end{prop}

\begin{rem} \label{remrecpresnonrecclosed}
The vector space spanned by $A_1,\dots,A_d$ is generally
  not recursively closed except if all symbols
$(\rho_1,R_1),\dots,(\rho_d,R_d)$ are of depth $0$.
\end{rem}

{\bf Proof of Proposition \ref{proprecdef}}
We have obviously $A_k[\emptyset,\emptyset]=\rho_k$ for $k=1,\dots,d$.
An easy induction on $l$ shows now that the matrices 
$$A_1[\mathcal M_{p\times q}^{l+1}],\dots,A_d[\mathcal 
M_{p\times q}^{l+1}]$$
are uniquely defined by 
$$A_1[\mathcal M_{p\times q}^{\leq l}],\dots,A_d[\mathcal 
M_{p\times q}^{\leq l}]\ .$$

The second part of the Proposition is obvious.\hfill$\Box$

Recursive presentations of depth $0$ are particularly nice:
Given such a recursive presentation
$A_1=(\rho_1,R(1)),\dots,A_d=(\rho_d,
R(d))$, the subspace spanned 
by $A_1,\dots,A_d\subset\hbox{Rec}_{p\times q}
(\mathbf K)$ is already recursively closed and the identity
$$R(k)_{s,t}=\rho(s,t)A_k=\sum_{j=1}^d\rho(s,t)_{j,k}A_j$$
shows that the matrices 
$R(1),\dots,R(d)$ encode the same information as shift matrices 
with respect to the generating set $A_1,\dots,A_d$. In particular,
every element $A\in\hbox{Rec}_{p\times q}(\mathbf K)$ of complexity
$d$ admits a recursive presentation of depth $0$ for $A_1=A,
\dots,A_d$ a basis of $\overline{A}^{rec}$.

\begin{rem}
Since the importance of the role played by the shift monoid is 
not apparent in the definition of recursive presentations 
they are less natural than monoidal presentations 
from a theoretical point of view.
They are however often more ``intuitive''
and more compact (see Remark \ref{remrecpresnonrecclosed} 
for the reason)
than monoidal presentations. Moreover, 
several interesting examples are
defined in a natural way in terms of recursion matrices while
a definition using shift matrices looks more artificial.
\end{rem}

\section{Saturation}

The aim of this section is to present finiteness results 
for  computations
in the category $\hbox{Rec}({\mathbf K})$. More precisely, 
given presentations of two suitable recurrence matrices $A,B$,
we show how to compute a presentation of the sum $A+B$, or
of the matrix-product $A B$ (whenever defined)
of $A$ and $B$ in a finite number of steps.

For a vector space ${\mathcal A}\subset{\mathbf K}^{{\mathcal
    M}_{p\times q}}$
we denote by ${\mathcal A}[{\mathcal M}_{p\times q}^{\leq l}]\subset
{\mathbf K}^{{\mathcal M}_{p\times q}^{\leq l}}$ the image
of ${\mathcal A}$ under the projection $A\longmapsto
A[{\mathcal M}_{p\times q}^{\leq l}]$ associating to 
$A\in{\mathbf K}^{{\mathcal M}_{p\times q}}$ its restriction
$A[{\mathcal
    M}_{p\times q}^{\leq l}]\in 
{\mathbf K}^{{\mathcal M}_{p\times q}^{\leq l}}$ to the
subset of all $\frac{(pq)^{l+1}-1}{pq-1}$
words of length at most $l$ in ${\mathcal M}_{p\times q}$.

\begin{defin} The {\it saturation level} of a non-zero vector space
${\mathcal A}\subset{\mathbf K}^{{\mathcal M}_{p\times q}}$ is the
smallest natural integer $N\geq 0$, if it exists, 
for which the obvious projection
$${\mathcal A}[{\mathcal M}_{p\times q}^{\leq N+1}]\longrightarrow
{\mathcal A}[{\mathcal M}_{p\times q}^{\leq N}]$$
is an isomorphism.
A vector-space ${\mathcal A}$ has saturation level $\infty$ 
if such an integer $N$ does not exist and 
the trivial vector space $\{0\}$ has
saturation level $-1$.
The {\it saturation level} of $A\in{\mathbf K}^{{\mathcal
    M}_{p\times q}}$ 
is the saturation level of its recursive closure
$\overline{A}^{rec}$.
\end{defin}

\begin{prop} \label{propfond} Let ${\mathcal A}\subset {\mathbf K}^{
{\mathcal M}_{p\times q}}$
be a recursively closed
vector space of finite saturation level $N$. Then
${\mathcal A}$ and ${\mathcal A}[{\mathcal
  M}_{p\times q}^{\leq N}]$ are isomorphic.

In particular, ${\mathcal A}$ is of finite dimension and
contained in $\hbox{Rec}_{p\times q}({\mathbf K})$.
\end{prop}

\begin{cor} \label{stabilizeddim} A finite set 
$A_1,\dots,A_d\subset \hbox{Rec}_{p\times q}({\mathbf K})$ spanning a 
recursively closed subspace ${\mathcal A}=
\sum_{j=1}^{d}{\mathbf K}A_j\subset \hbox{Rec}_{p\times q}({\mathbf K})$ 
is linearly independent if and only if
$A_1[{\mathcal M}_{p\times q}^{\leq N}],\dots,A_d[{\mathcal M}_{
p\times q}^{\leq N}]\in {\mathbf K}^{{\mathcal M}_r^{\leq N}}$ 
are linearly independent where $N\leq\hbox{dim}
({\mathcal A})-1< d$ denotes 
the saturation level of 
$\mathcal A$.
\end{cor} 

{\bf Proof of Proposition \ref{propfond}} We denote by 
$K_l=\{A\in{\mathcal A}\ \vert \ A[U,W]=0\hbox{ for all }
(U,W)\in{\mathcal M}_{p\times q}^{\leq l}\}
\subset {\mathcal A}$ the kernel of the natural projection
${\mathcal A}\longrightarrow {\mathcal A}[{\mathcal M}_{p\times
  q}^{\leq l}]$. We have
$$K_{N+1}=\{A\in K_N\ \vert\ \rho({\mathcal M}_{p\times q}^1)A\subset
K_N\}=K_N$$
which shows the equality $\rho({\mathcal M}_{p\times q})K_N=K_N$.
For $A\in K_N\subset K_0$ we have thus
$A[U,W]=\rho(U,W)A[\emptyset,\emptyset]=0$
for all $(U,W)\in{\mathcal M}_{p\times q}$ which implies
$A=0$ and shows the result.\hfill $\Box$

Corollary \ref{stabilizeddim} follows immediately.

\begin{rem} The proof of Proposition \ref{propfond} shows the
inequality $N+1\leq \hbox{dim}(\overline{A}^{rec})$
for the saturation level $N$ of $A\in
  \hbox{Rec}_{p\times q}({\mathbf K})$. Equality is achieved eg. for 
the recurrence matrix defined by $A[U,W]=1$ if $(U,W)\in {\mathcal
  M}_{p\times q}^N$
and $A[U,W]=0$ otherwise.
\end{rem}

\section{Algorithms}

It is easy to extract a minimal
presentation from a presentation ${\mathcal P}$ of a recurrence
matrix $A\in
\hbox{Rec}_{p\times q}({\mathbf K})$: Compute the saturation level $N$ of 
the recursively closed space ${\mathcal A}=\sum_{j=1}^d {\mathbf
  K} A_j$ defined by the presentation ${\mathcal P}$. 
This allows to detect linear dependencies among $A_1,\dots A_d$. 

Using Proposition \ref{propcoeffs}, these computations can be done in
polynomial time with respect to the complexity $d$ of the
presentation ${\mathcal P}$ for $A$.

The algorithm described in subsection \ref{subssatlev} 
is useful for computing the
saturation level and a few other items attached to a recurrence matrix.

Adding two recurrence matrices $A,B$ in $\hbox{Rec}_{p\times q}
({\mathbf K})$ is now easy: given presentations
$A_1,\dots, A_a$ and $B_1,\dots, B_b$ of $A$ and $B$, 
one can write down a presention 
$C_1=A_1+B_1,C_2=A_2,\dots,C_a=A_a,C_{a+1}=B_1,C_{a+2}=B_2,\dots,
C_{a+b}=B_b$ of $A+B$, eg. by using Proposition \ref{proppressum}. 

Similarly, using Proposition \ref{presentproduct}
a presentation $A_1,\dots, A_a$ and
$B_1,\dots,B_b$ of $A\in\hbox{Rec}_{p\times r}
({\mathbf K})$ and $B\in\hbox{Rec}_{r\times q}
({\mathbf K})$ yields a presentation 
$C_{11}=A_1B_1,\dots, C_{ab}=A_aB_b$ of $C=C_{11}=A_1B_1$.

\begin{rem} Computing a presentation of $AB$ for $A\in
  \hbox{Rec}_{p\times r}({\mathbf K})$, $B\in \hbox{Rec}_{r\times
    q}({\mathbf K})$ is also possible using the matrices
$(A_iB_j)[{\mathcal M}_{p\times q}^l]$ for $l$ up to the 
saturation level $N$ of $\overline{A}^{rec}\overline{B}^{rec}$.
This method is however quickly infeasible for 
for $r\geq 2$ and $N$ not too small. Using 
Proposition \ref{presentproduct} as suggested above, one
gets around this difficulty and obtains essentially polynomial
algorithms for computations in the algebra $\hbox{Rec}_{p\times p}
(\mathbf K)$.
\end{rem}

\subsection{An algorithm for computing the saturation level}\label{subssatlev}

\begin{prop} \label{nicecoord}
Given a recursively closed vector space ${\mathcal
    A}\subset \hbox{Rec}_{p\times q}({\mathbf K})$ of dimension $a$,
there exists a subset ${\mathcal S}=\{S_1,\dots,S_a\}\subset {\mathcal
  M}_{p\times q}$ of words such that the restriction 
$$X\longmapsto X[{\mathcal S}]=(X[S_1],\dots,X[S_a])$$ 
of $X$ onto ${\mathcal S}$ induces an isomorphism between 
${\mathcal A}$ and ${\mathbf K}^a$. Moreover, one can choose 
${\mathcal S}$ in order to have 
$$\mathcal S=\subset (\emptyset,\emptyset)\cup{\mathcal S}{\mathcal
    M}_{p\times q}^1\ .$$
\end{prop}

{\bf Proof} The first part of the proof is obvious. In order to have
the inclusion $\mathcal S\subset (\emptyset,\emptyset)
\cup{\mathcal S}{\mathcal M}_{p\times q}^1$ suppose that the
evaluation of $X$ on $\{S_1,\dots,S_{k_j}\}\subset {\mathcal S}\cap 
{\mathcal M}_{p\times q}^{\leq j}$ induces a bijection between 
${\mathcal A}[{\mathcal M}_{p\times q}^{\leq j}]$ and ${\mathbf
  K}^{k_j}$
for $j\leq N$ with $N$ the saturation level of ${\mathcal A}$.
Choosing a presentation $A_1,\dots,A_d,\rho(s,t)\in{\mathbf
  K}^{d\times d}, 0\leq s<p,0\leq t<q$ and writing
$$\rho^t(\overline S)=\rho^t(s_{j_n},t_{j_n})\cdots
\rho^t(s_{j_1},t_{j_1})$$ for $S=(s_{j_1}\cdots s_{j_n},t_{j_1},\cdots
t_{j_n})\in{\mathcal M}_{p\times q}^n$, Proposition 
\ref{propcoeffs} shows that the $k_j$
vectors 
$$\rho^t(S_j)\left((A_1,\dots,A_d)[\emptyset,\emptyset]\right)^t,
1\leq j\leq k_j$$ form a basis of the vector space spanned by 
$$\rho^t({\mathcal M}_{p\times q}^{\leq j})
\left((A_1,\dots,A_d)[\emptyset,\emptyset]\right)^t\subset {\mathbf
  K}^d
\ .$$ This implies easily the result. \hfill$\Box$

Proposition \ref{nicecoord} implies that the following algorithme 
computes the saturation level of a presentation.

Input: a finite set $A_1,\dots,A_d\subset \hbox{Rec}_{p\times q}(
{\mathbf K})$ spanning a recursively closed subspace (eg. given by a
finite presentation).

Set $a:=0,\  N:=-1,\ {\mathcal S}:=\emptyset$,

If $(A_1,\dots,A_d)[\emptyset,\emptyset]=(0,\dots,0)$ then $
a:=0,\  N:=-1,\ {\mathcal S}:=\emptyset,\ $stop else
$a:=1,\  k:=1,\ u:=1,\ N:=0,\ S_1:=\{(\emptyset,\emptyset)\}\
  $ endif

For $j\in\{k,k+1,\dots,u\}$ do

For $(s,t)\in{\mathcal M}_{p\times q}^1$ do

If $\rho^t(s,t)\rho^t(\overline S_j)\left((A_1,\dots,A_d)[\emptyset,
\emptyset]\right)^t$ is not in the linear span of 
$\rho^t(\overline S)\left((A_1,\dots,A_d)[\emptyset,
\emptyset]\right)^t,\ S\in \{S_1,\dots,S_u\}$
then $a:=a+1,\ S_a:=S_j(s,t)$ endif (where 
$\overline{(s_1s_2\dots s_m,t_1t_2\dots t_m)}=(s_ms_{m-1}\dots s_2s_1,
t_mt_{m-1}\dots t_2t_1)$)

If $a=d$ then stop endif

endfor

If $a=u$ then stop else $k:=u+1, u:=a, N:=N+1$ endif

endfor

The final value of $a$ in this algorithm is the dimension
of the recursively closed vector space ${\mathcal A}=\sum_{j=1}^d
{\mathbf K}A_j$, the final value of $N$ is the saturation level
${\mathcal A}$ and the set ${\mathcal S}=\{S_1,\dots,S_a\}$ 
satisfies the conditions of Proposition \ref{nicecoord}.


\section{Topologies and a metric}\label{topologies}

For ${\mathbf K}$ a topological field, the vector space 
${\mathbf K}^{{\mathcal M}_{p\times q}}$ carries two natural
``obvious'' topologies.

The first one is the product topology on $\prod_{l=0}^\infty
{\mathbf K}^{{\mathcal M}_{p\times q}^l}$. It is defined as
the coarsest topology for which all projections
$A\longmapsto A[{\mathcal M}_{p\times q}^l]$ are continuous.
Its open subsets are generated by ${\mathcal O}_{\leq m}\times
\prod_{l>m}{\mathbf K}^{{\mathcal M}_{p\times q}^l}$ with
${\mathcal
  O}_{\leq m}\subset {\mathbf K}^{{\mathcal M}_{p\times q}^{\leq m}}$
open for the natural topology on the finite-dimensional
vector space ${\mathbf K}^{{\mathcal M}_{p\times
 q}^{\leq m}}$.

The second topology is the strong topology or box topology.
Its open subsets are generated by
$\prod_{l=0}^\infty {\mathcal O}_l$, with ${\mathcal O}_l\subset
{\mathbf K}^{{\mathcal M}_{p\times q}^l}$ open for all $l\in{\mathbb N}$.

The restrictions
to $\hbox{Rec}_{p\times q}({\mathbf K})$ of both topologies 
are not very interesting: The product topology is very coarse, 
the strong topology is too fine: it gives 
rise to the discrete topology on
the quotient $\hbox{Rec}_{p\times q}({\mathbb C})/{\mathcal
  FS}_{p\times q}$ considered in Chapter \ref{finsuppquotient}.

For a recurrence matrix $A\in\hbox{Rec}_{p\times q}({\mathbb C})$ 
defined over ${\mathbb C}$ we set
$$\parallel
A\parallel_\infty^\infty=\hbox{limsup}_{l\rightarrow\infty}
\hbox{sup}_{(U,W)\in{\mathcal M}_{p\times q}^l}\vert
A[U,W]\vert^{1/l}.$$

\begin{rem} $\parallel
A\parallel_\infty^\infty$ equals the supremum over ${\mathcal
  M}_{p\times q}$ of 
the numbers $\sigma(S,T)^{1/l}$ where $(S,T)\in{\mathcal M}_{
p\times q}^l$ has length $l$ and where $\sigma(S,T)$ 
is the spectral radius (largest modulus among eigenvalues)
of the shift matrix $\rho_{\overline{A}^{rec}}(S,T)$
with respect to a minimal presentation of $A$. We have thus
the inequality
$\sigma_{s,t}\leq \parallel A\parallel_\infty^\infty$
(which is in general strict) where $\sigma_{s,t}$ are the spectral
radii of shift-matrices $\rho(s,t)\in{\mathbb C}^{d\times d}$
associated to a minimal presentation of $A$.
\end{rem}

\begin{prop} \label{propnormineq}
The application $A\longmapsto \parallel A
\parallel_\infty^\infty$ defines a metric on ${\mathbb C}^*\backslash
\hbox{Rec}_{p\times q}({\mathbb C})/\mathcal{FS}_{p\times q}$
such that
$$\parallel A+\lambda B\parallel_\infty^\infty\leq \hbox{sup}
(\parallel A\parallel_\infty^\infty,
\parallel B\parallel_\infty^\infty)$$
for $\lambda\in{\mathbb C}^*,A,B
\in\hbox{Rec}_{p\times q}({\mathbb C}) $ (with equality holding 
generically) and 
$$\parallel AB\parallel_\infty^\infty\leq 
r \parallel A\parallel_\infty^\infty\ 
\parallel B\parallel_\infty^\infty$$
for $A\in\hbox{Rec}_{p\times r}({\mathbb C}), 
B\in\hbox{Rec}_{r\times q}({\mathbb C})$.
\end{prop}

{\bf Proof} The proof of the first inequality is easy and left to the
reader. The second inequality follows from
$$\parallel AB\parallel_\infty\leq \parallel \tilde A\tilde B
\parallel_\infty=r\parallel
A\parallel_\infty\parallel B\parallel_\infty$$
where $\tilde A,\tilde B$ are of complexity $1$ and have coefficients
$\tilde A[S,T]=\parallel A\parallel_\infty^l$ for
$(S,T)\in{\mathcal M}_{p\times r}^l$ and
$\tilde B[S,T]=\parallel B\parallel_\infty^l$ for
$(S,T)\in{\mathcal M}_{r\times q}^l$ which depend only on the length of
$(S,T)$.
The inequality $\parallel AB\parallel_\infty\leq 
\parallel \tilde A\tilde B\parallel_\infty$
follows now easily from the equalities
$\parallel A\parallel_\infty=\parallel \tilde A\parallel_\infty$, 
$\parallel B\parallel_\infty=\parallel \tilde B\parallel_\infty$
and from the observation that all coefficients of $\tilde A$
(respectively $\tilde B$) are upper bounds for the corresponding
coefficients of $A$ (respectively $B$).\hfill $\Box$

\begin{rem} The set of all elements $A\in {\mathbb C}^{{\mathcal
      M}_{p\times q}}, p,q\in{\mathbb N}$, such that 
$\parallel A\parallel_\infty^\infty<\infty$ form a subcategory
containing $\hbox{Rec}({\mathbb C})$ of ${\mathbb C}^{\mathcal M}$. 

The vector space $\hbox{Rec}_{p\times q}({\mathbb C})$
can be normed by
$$\parallel A\parallel =\sum_{l=0}^\infty \frac{\parallel A[{\mathcal
    M}_{p\times q}^l]\parallel_\infty}{l!}$$
where $\parallel A[{\mathcal
    M}_{p\times q}^l]\parallel_\infty$ denotes the largest absolute 
value of all coefficients $A[U,W],(U,W)\in{\mathcal M}_{p\times q}^l$.
However, matrix-multiplication is unfortunately
not continuous for this norm.

This norm has different obvious variations: 

The factorials
$l!$ of the denominators can be replaced by any other sequence $s_0,
s_1,\dots $ of strictly positive numbers such that 
$\hbox{lim}_{n\rightarrow \infty}\lambda^n/s_n=0$ for all $\lambda 
>0$.

The sup-norm $\parallel A[{\mathcal
    M}_{p\times q}^l]\parallel_\infty$ can be replaced by many other
``reasonable'' norms (like $l_1$ or $l_2$ norms) on 
${\mathbb C}^{{\mathcal M}_{p\times q}^l}$.

It would be interesting to find a norm on the algebra
$\hbox{Rec}_{p\times p}({\mathbb C})$ for which matrix 
products are continuous.
\end{rem}


\section{Criteria for non-recurrence matrices}

This short section lists a few easy properties (which have
sometimes obvious generalizations, eg by replacing
${\mathbb Q}$ with a number field) of recurrence matrices and
non-recurrence matrices in ${\mathbf K}^{{\mathcal M}_{p\times q}}$.
Proofs are straightforward and left to the reader or only sketched. 

\begin{prop} An element $A\in {\mathbf K}^{{\mathcal M}_{p\times q}}$
is not in $\hbox{Rec}_{p\times q}({\mathbf K})$ if and only if there
exist two sequences $(S_i,T_i),(U_j,W_j)\in ({\mathcal M}_{p\times
  q})^{\mathbb N}$ such that for all $n\in{\mathbb N}$, the $n\times n$
matrix with coefficients
$c_{i,j}=\rho(S_i,T_i)A[U_j,W_j]=A[U_jS_i,W_jT_i], 1\leq i,j\leq n,$
has non-zero determinant.
\end{prop}

\begin{prop} \label{majorationcriterion}
Consider $A\in \hbox{Rec}_{p\times q}({\mathbb C})$. Then
there exists a constant $C\geq 0$ such that
$\vert A[U,W]\vert \leq C^{n+1}$ for all $(U,W)\in{\mathcal
  M}_{p\times q}^n$.
\end{prop}

{\bf Proof} Given a finite presentation $A_1,\dots A_d$ with shift
matrices $\rho(s,t)\in  {\mathbb C}^{d\times d}$,  
choose $C\geq 0$ such that $C\geq\vert
A_i[\emptyset,
\emptyset]\vert$ for $i=1,\dots,d$ and $C\geq d\parallel \rho(s,t)
\parallel_\infty$ for $0\leq s<p,0\leq t<q$.\hfill$\Box$

\begin{rem} \label{majorationcriterionrem}
Associate to $A\in{\mathbb C}^{{\mathcal M}_p}$ the
formal power series
$$f_A=A[\emptyset]+\sum_{n=1}^\infty \sum_{0\leq u_1,\dots ,u_n<p}
A[u_1\dots u_n]Z_{u_1}\cdots Z_{u_n}$$
in $p$ commuting variables $Z_0,\dots, Z_{p-1}$. 
Proposition \ref{majorationcriterion} shows
that $f_A$ defines a holomorphic function in a 
neighbourhood of $(0,\dots,0)\in{\mathbb C}^p$ if $A\in
\hbox{Rec}_p({\mathbb C})$.
\end{rem}

\begin{prop} \label{critfingenring}
The set of values $\{A[U,W]\ \vert\ 
(U,W)\in{\mathcal M}_{p\times q}\}$ of a recurrence 
matrix $A\in \hbox{Rec}_{p\times q}({\mathbf K})$
is contained in a subring $\tilde{\mathbf K}\subset {\mathbf
  K}$ which is finitely generated (as a ring).
\end{prop}

{\bf Proof} The subring $\tilde {\mathbf K}$ generated by all values
and coefficients involved in a finite presentation
$A=A_1,\dots,A_d$ of $A$ works.\hfill$\Box$

\begin{prop} \label{critfinprimes}
For $A\in\hbox{Rec}_{p\times q}({\mathbb Q})$, there exists
a natural integer $N$ such that 
$N^{n+1}A[U,W]\in{\mathbb Z}$ for all
$(U,W)\in{\mathcal M}_{p\times q}^n$. 
\end{prop}

{\bf Proof} Given a presentation $A_1,\dots,A_d\in \hbox{Rec}_{p\times
  q}({\mathbb Q})$, a non-zero integer $N\in {\mathbb N}$ such that
$N(A_1,\dots,A_d)[\emptyset,\emptyset]\in{\mathbb Z}^d$ and 
$N\rho(s,t)\in{\mathbb Z}^{d\times d}$ for all $0\leq s<p,0\leq t<q$
works.\hfill$\Box$


\section{Diagonal and lower triangular
  subalgebras in ${\mathbf K}^{{\mathcal M}_{p\times p}}$
and $\hbox{Rec}_{p\times p}({\mathbf K})$}\label{diagsubalg}

This section describes
a maximal commutative subalgebra formed by diagonal elements,
the center (which is contained as a subalgebra in the diagonal algebra) 
and the
lower triangular subalgebra of the algebras 
${\mathbf K}^{{\mathcal M}_{p\times p}}$
and $\hbox{Rec}_{p\times p}({\mathbf K})$.

\subsection{The diagonal subalgebra and the center}

An element $A\in {\mathbf K}^{{\mathcal M}_{p\times p}}$
is {\it diagonal} if $A[U,W]=0$ for all $(U,W)\in{\mathcal M}_{p\times
  p}$ such that $U\not= W$. 
We denote by ${\mathcal D}_p({\mathbf K})\subset {\mathbf K}^{{\mathcal M}_{p\times 
p}}$ the vector space of all diagonal elements. It is easy to show
that ${\mathcal D}_p({\mathbf K})$ is a commutative algebra which is 
isomorphic to the function ring underlying ${\mathbf K}^{{\mathcal
    M}_p}$.

The center of the algebra ${\mathbf K}^{{\mathcal M}_{p\times p}}$
is the subalgebra ${\mathcal C}_p({\mathbf K})\subset {\mathcal D}_p({\mathbf K})$
formed by all diagonal matrices $A$ with diagonal coefficients 
$A[U,U]$ depending
only on the length $l$ of $(U,U)\in{\mathcal M}_{p\times p}^l$.
For $p\geq 1$, the map
$${\mathcal C}_p({\mathbf K})\ni C\longmapsto
(C[\emptyset,\emptyset],C[0,0],C[0^2,0^2],\dots)$$
defines an isomorphism between ${\mathcal C}_p({\mathbf K})$ 
and the algebra
${\mathbf K}^{{\mathcal M}_{1\times 1}}$ corresponding to all
sequences ${\mathbb N}\longrightarrow {\mathbf K}$,
endowed with the coefficient-wise (or Hadamard) product.

We denote by ${\mathcal D}_{p-rec}({\mathbf K})=
{\mathcal D}_p({\mathbf K})\cap \hbox{Rec}_{p\times p}({\mathbf K})$ 
and by ${\mathcal C}_{p-rec}({\mathbf K})=
{\mathcal C}_p({\mathbf K})\cap \hbox{Rec}_{p\times p}({\mathbf K})$
the subalgebras of ${\mathcal D}_p({\mathbf K})$ and ${\mathcal C}_p({\mathbf K})$
formed by all recurrence matrices. Associating to $C\in{\mathcal C}_
{p-rec}({\mathbf K})$ the generating series $\sum_{l=0}^\infty
C[0^l,0^l]z^l$ (where $0^0=\emptyset$) yields an isomorphism between
the algebra ${\mathcal C}_{p-rec}(
{\mathbf K})$ and the algebra $\hbox{Rec}_{1\times
  1}({\mathbf K})$ corresponding to the vector-space
of rational functions in one variable
without singularity at the origin. The product is the coefficientwise
product $\left(\sum_{l=0}^\infty \alpha_lz^l\right)\cdot
\left(\sum_{l=0}^\infty \beta_lz^l\right)=\left(\sum_{l=0}^\infty
  \alpha_l\beta_lz^l\right)$ of the corresponding series-expansions.

Diagonal and central recurrence matrices in $\hbox{Rec}_{p\times
  p}({\mathbf K})$ can be characterized by the following result.

\begin{prop} (i) A recurrence matrix $A\in\hbox{Rec}_{p\times
    p}({\mathbf K})$ is diagonal if and only if it can be given by 
a presentation with shift-matrices $\rho(s,t)=0$ whenever $s\not= t,
0\leq s,t<p$.

(ii) A diagonal recurrence matrix $A\in\hbox{Rec}_{p\times
    p}({\mathbf K})$ is central if and only if it can be given by 
a presentation with shift-matrices $\rho(s,t)=0$ whenever $s\not= t,
0\leq s,t<p$ and $\rho(s,s)=\rho(0,0)$ for all $s, 0\leq s<p$.
\end{prop}

We leave the easy proof to the reader.

\begin{rem} The vector spaces ${\mathcal D}_p,{\mathcal D}_{p-rec},
{\mathcal C}_p,{\mathcal C}_{p-rec}$ are also (non-commutative for
${\mathcal D}_p$ and ${\mathcal D}_{p-rec}$, if $p>1$) algebras for the
convolution-product.
\end{rem}

\subsection{Lower triangular subalgebras}\label{Toeplitz}

Lower (or upper) triangular elements in ${\mathbf K}^{{\mathcal
    M}_{p\times p}}$ can be defined using the bijection
$$(u_1\dots u_l,w_1\dots w_l)\longmapsto(\sum_{j=1}^l
u_jp^{j-1},\sum_{j=1}^lw_jp^{j-1})$$
between ${\mathcal M}_{p\times p}^l$ and $\{0,\dots,p^l-1\}\times
\{0,\dots,p^l-1\}$. More precisely, $A\in {\mathbf K}^{{\mathcal
    M}_{p\times p}}$ is {\it lower triangular} if 
for all $l\in{\mathbb N}$ the equality 
$A[U,W]=0$ holds for
$(U,W)=(u_1\dots u_l,w_1\dots w_l)\in {\mathcal M}_{p\times p}$
such that $\sum_{j=1}^l u_jp^{j-1}<\sum_{j=1}^l w_jp^{j-1}$.

Similarly, an element $A\in{\mathbf K}^{{\mathcal M}_{p\times p}}$ is
{\it upper triangular} if the transposed element $A^t$ (defined by
$A^t[U,W]=A[W,U]$) is lower triangular. 

We denote by ${\mathcal L}_{p}({\mathbf K})\subset {\mathbf K}^{{
\mathcal M}_{p\times p}}$ the vector-space of all lower triangular
elements in ${\mathbf K}^{{\mathcal M}_{p\times p}}$. It is easy to
check that ${\mathcal L}_{p}({\mathbf K})$ is closed under the
matrix-product and the vector space ${\mathcal L}_p({\mathbf K})$ is
thus an algebra. The subspace 
of all convergent elements in ${\mathcal L}_{p}({\mathbf K})$
is also closed under the matrix-product and forms a subalgebra.

A lower triangular matrix $A\in{\mathcal L}_{p}({\mathbf K})$ is
{\it unipotent} if $A[U,U]=1$ for all ``diagonal
words'' $(U,U)\in{\mathcal M}_{p\times p}$ and {\it strictly
lower triangular} if $A[U,U]=0$ for all diagonal
words $(U,U)\in{\mathcal M}_{p\times p}$. The subset 
${\mathcal N}_p({\mathbf K})$ of all
lower strictly-triangular matrices is a two-sided ideal in
${\mathcal L}_{p}({\mathbf K})$. The associated quotient
${\mathcal L}_{p}({\mathbf K})/{\mathcal N}_p({\mathbf K})$
is isomorphic to the commutative algebra ${\mathcal D}_p({\mathbf K})$
of diagonal elements. Unipotent lower triangular matrices 
form a multiplicative subgroup and correspond to the
(multiplicative) identity of the quotient algebra
${\mathcal L}_{p}({\mathbf K})/{\mathcal N}_p({\mathbf K})$.

The following proposition
is useful for recognizing triangular recurrence matrices: 

\begin{prop} \label{triangrecmat}
A recurrence matrix $A\in \hbox{End}_{p-rec}({\mathbf K})$
is triangular if and only if it admits a presentation of the form $A=A_1,
\dots,A_{k},A_{k+1},\dots,A_a$ such that
$\rho(s,s)A_1,\dots,\rho(s,s)A_{k}\in \sum_{j=1}^{k} A_j$
for $0\leq s<p$ and
$\rho(s,t)A_1,\dots,\rho(s,t)A_{k}=0$ for $0\leq s<t<p$.

Such a recurrence matrix $A$ is unipotent if and only if 
it admits a presentation as 
above which satisfies moreover $A_1[\emptyset,\emptyset]=1,A_2[
\emptyset,\emptyset]=\dots=A_{k}[\emptyset,\emptyset]=0$,
$\rho(s,s)A_1$ is in the affine 
space $A_1+\sum_{j=2}^{k} {\mathbf K}A_j$ for $0\leq s<p$, and
$\rho(s,s)A_h\in\sum_{j=2}^{k} {\mathbf K}A_j$ for $2\leq h\leq k$
and $0\leq s<p$.
\end{prop}

{\bf Proof} These conditions are clearly sufficient
since they imply by an easy induction on $l$ that
all matrices $A_1[{\mathcal M}_{p\times p}^l],\dots,A_{k}
[{\mathcal M}_{p\times p}^l]$ are lower triangular.
In the unipotent case they imply
that $A_1[{\mathcal M}_{p\times p}^l]$ is unipotent and that 
$A_2[{\mathcal M}_{p\times p}^l],\dots,A_{k}[{\mathcal M}_{p\times p}^l]$
are strictly lower triangular for all $l\in{\mathbb N}$.

In order to see that they are also necessary, 
consider a basis $A_1=A,A_{k}$ of the vector space ${\mathcal T}
\subset \overline{A}^{rec}$ spanned by all lower
triangular recurrence matrices of the form
$\rho(S,S)A$, with $(S,S)\in{\mathcal M}_{p\times p}$.
In the unipotent case, all 
diagonal coefficients of an element $X\in{\mathcal T}$
are equal and we can consider a basis $A_2,\dots,A_{k}$ 
spanning the subspace of all strictly inferior triangular matrices 
in ${\mathcal T}$.\hfill $\Box$

\begin{rem} The proof of Proposition \ref{triangrecmat} constructs the
  subspace 
${\mathcal T}\subset \overline A^{rec}$ spanned by all elements of the form
$\rho(S,S)A\in\overline A^{rec}$. This subspace is contained in the
maximal vector space spanned by all lower
triangular matrices in $\overline A^{rec}$.
\end{rem}


\section{Elements of Toeplitz type}
\label{toeplitzsection}

In this section we identify the set ${\mathcal M}_p^l$ 
with $\{0,\dots, p^l-1\}$ using the bijection 
$$u_1\dots u_l\longmapsto \sum_{j=1}^l u_jp^{j-1}.$$
Analogously, we identify ${\mathcal M}_{p\times p}$ 
in the obvious way with $\{0,\dots, p^l-1\}\times\{0,\dots, p^l-1\}$.
This identification yields a natural isomorphism between 
${\mathbf K}^{{\mathcal M}_p^l}$ and vectors
$(\alpha_0^l,\dots, \alpha_{p^l-1}^l)\in {\mathbf K}^{p^l}$,
respectively between ${\mathbf K}^{{\mathcal M}_{p\times p}^l}$
and matrices with coefficients $(\alpha_{i,j}^l)_{0\leq i,j<p^l}$.

\subsection{Toeplitz matrices}\label{Toeplitzmatrices}

A (finite or infinite) matrix $T$ of square-order $n\in {\mathbb
  N}\cup\{\infty\}$ with coefficients $t_{u,w},
0\leq u,w<n$ is a 
{\it Toeplitz matrix} if $t_{u,w}$ depends only on the difference 
$u-w$ of its indices. A Toeplitz matrix is thus described by a
(finite or infinite) 
sequence $\dots,\alpha_{-2},\alpha_{-1},\alpha_0,\alpha_1,\dots$
defined by $\alpha_{u-w}=t_{u,w}$.

We call a matrix algebra ${\mathcal A}\subset M_{n\times
  n}({\mathbf K})$ of $n\times n$ square matrices (with $n\in{\mathbb
  N}\cup\{\infty\}$ finite or infinite) an algebra of  {\it Toeplitz
  type} if all elements of ${\mathcal A}$ are Toeplitz matrices.
For finite $n$, the typical example is the
$n-$dimensional commutative algebra $\hbox{Toep}_\rho(n)=
\sum_{j=0}^{n-1}{\mathbf K}T_\rho^j$
generated by the $n\times n$ matrix
$$T_\rho=\left(\begin{array}{cccccc}
0& &\dots&0&\rho\\
1&0& &     & 0\\
0&1&0&     &\vdots  \\
\vdots& &\ddots & & \\
0& &\dots&1&0\end{array}\right)$$
with minimal polynomial $T_\rho^n-\rho$ obtained by
considering the product $T_1D_{n,\rho}$ where
$T_1$ is the obvious cyclic permutation matrix of order $n$ and where 
$D_{n,\rho}$ is the diagonal matrix with diagonal entries
$1,1,\dots,1,\rho$.

An example of an element in $\hbox{Toep}_1(2^n)$ is the matrix
with coefficients $t_{u,w}=\gamma_{v_2(u-w)},\ 0\leq u,w<2^n$ 
depending only on the highest power 
$2^{v_2(u-w)}\in \{1,2,4,8,16,\dots,\}\cup\infty$ dividing $u-w$. 
If such a matrix 
$$\left(\begin{array}{ccccccc}
\gamma_\infty&\gamma_0&\gamma_1&\gamma_0&\gamma_2&\gamma_0&\dots\\
\gamma_0&\gamma_\infty&\gamma_0&\gamma_1&\gamma_0&\gamma_2\\
&&\ddots\end{array}\right)$$
is invertible, then its inverse is of the same type.

In the sequel, we will mainly consider the algebra $\hbox{Toep}(n)=
\hbox{Toep}_0(n)$ of all lower triangular Toeplitz matrices
having finite or infinite square order $n\in {\mathbb N}\cup\{\infty\}$.
We call $\hbox{Toep}(n)$ the {\it (lower triangular)
Toeplitz algebra} of order $n$.

The proof of the following well-known result is easy and 
left to the reader.

\begin{prop} \label{isotoeppoly} 
For finite $n$, the algebra
$\hbox{Toep}(n)$ is isomorphic to the ring
${\mathbf K}[x]\pmod {x^n}$ of polynomials modulo $x^n$ and
$\hbox{Toep}(\infty)$ is isomorphic to 
the ring ${\mathbf K}[[x]]$ of formal power series.

In both cases, the isomorphism is given by considering the 
generating series $\sum_{j=0}^{n-1} t_{j,0} x^j$ 
associated to the first column of a lower triangular
Toeplitz matrix $t_{i,j},0\leq i,j$.
\end{prop} 

An element $A\in{\mathbf K}^{{\mathcal M}_{p\times p}}$ is
of {\it Toeplitz type} if all matrices
$A[{\mathcal M}_{p\times p}^l]$ are
Toeplitz matrices.
{\it Algebras of Toeplitz type} in ${\mathbf K}^{{\mathcal M}_{p\times p}}$
or $\hbox{Rec}_{p\times p}({\mathbf K})$ are defined in the obvious
way as containing only elements of Toeplitz type.

The algebra ${\mathcal T}_p({\mathbf K})$ of lower triangular 
elements of Toeplitz type in ${\mathbf K}^{{\mathcal M}_{p\times p}}$
will be studied below.

Shift maps $\rho(s,t)$
preserve the vector space of (recurrence) matrices of Toeplitz type
in ${\mathbf K}^{{\mathcal M}_{p\times p}}$. 
The recursive closure $\overline{T}^{rec}$ of an element of 
Toeplitz type in ${\mathbf K}^{{\mathcal M}_{p\times p}}$ contains
thus only elements of Toeplitz type. We have the following result:

\begin{prop} \label{Toeplitzprop} All recurrence matrices
  $T_1,\dots,T_{d}\in\hbox{Rec}_{p\times p}
({\mathbf K})$ of  a reduced presentation are of Toeplitz type if
and only if the shift matrices $\rho(s,t)$ 
satisfy the following two conditions:

\ \ (1) The shift matrices $\rho(s,t),\ 0\leq s,t<p$ depend only 
on $s-t$.

\ \ (2) Assuming condition (1) and writing (somewhat abusively)
$\rho(s-t)$ for a shift matrix $\rho(s,t)$ we have
the identities
$$\rho(s+p)\rho(t)=\rho(s)\rho(t+1)$$
for $1-p\leq s<0$ and $1-p\leq t<p-1$.
\end{prop}

For $p=2$, the conditions of Proposition \ref{Toeplitzprop} boil down to 
the identities $\rho(0,0)=\rho(1,1)$ and 
$$\begin{array}{l}
\displaystyle \rho(1,0)\rho(0,0)=\rho(0,1)\rho(1,0),\\
\displaystyle \rho(1,0)\rho(0,1)=\rho(0,1)\rho(0,0).\end{array}$$

{\bf Proof of Proposition \ref{Toeplitzprop}} We show by induction on
$l$ that all matrices $T_i[{\mathcal M}_{p\times p}^l]$ are
Toeplitz matrices.

For $l=0$, there is nothing to do. Condition (1) implies that 
$T_1[{\mathcal M}_{p\times p}^1],\dots,T_d[{\mathcal
  M}_{p\times p}^1]$ are Toeplitz matrices. 
Setting $T=T_i$ and 
denoting $(\rho(s,t)T)[{\mathcal M}_{p\times p}^l]$ by $\rho(s-t)T$ we get
$$T[{\mathcal M}_{p\times p}^{l+1}]=\left(\begin{array}{ccccc}
\rho(0)T&\rho(-1)T&\rho(-2)T&\cdots&\rho(1-p)T\\
\rho(1)T&\rho(0)T&\rho(-1)T&\cdots&\rho(2-p)T\\
\vdots&          & \ddots       &      &   \vdots\\
\rho(p-2)T&\rho(p-3)T&\rho(p-4)T&\cdots&\rho(1)\\
\rho(p-1)T&\rho(p-2)T&\rho(p-3)T&\cdots&\rho(0)T\end{array}\right)$$
where all matrices 
$\rho(1-p)T,\dots,\rho(p-1)T$ are Toeplitz 
matrices by induction. 
We have to show that two horizontally adjacent
blocks $(\rho(t+1)T\ \rho(t)T),\ 1-p\leq t<p-1$ 
define a $p^l\times (2\cdot p^l)$
matrix of ``Toeplitz-type'' (the case of vertically adjacent blocks
gives rise to the same conditions and is left to the reader).
We have
$$\big(\rho(t+1)T\vert \rho(t)T\big)=\left(\begin{array}{cr|lc}
\cdots&\rho(1-p)\rho(t+1)T&\rho(0)\rho(t)T&\cdots\\
\cdots&\rho(2-p)\rho(t+1)T&\rho(1)\rho(t)T&\cdots\\
      &\vdots\qquad             &\rho(2)\rho(t)T&\cdots\\
&\vdots\qquad&\qquad\vdots&\\
\cdots&\rho(-1)\rho(t+1)T&\qquad\vdots&\\
\cdots&\rho(0)\rho(t+1)T&\rho(p-1)\rho(t)T&\cdots\end{array}\right)$$
(the vertical line separates the matrix into two square blocks)
where $\rho(\alpha)\rho(\beta)T=(\rho(\alpha)\rho(\beta)T)[
{\mathcal M}_{p\times p}^{l-1}]$. If
$$\rho(s-p)\rho(t+1)T=\rho(s)\rho(t)T,\ 1\leq s<p-1,$$
the matrix $T=T_i[{\mathcal M}_{p\times p}^{l+1}]$ is of Toeplitz
type.

The opposite direction is easy and left to the reader.\hfill $\Box$

\subsection{The algebra ${\mathcal T}_p({\mathbf K})\subset 
{\mathbf K}^{{\mathcal
      M}_{p\times p}}$ formed by lower triangular elements of Toeplitz
  type}
\label{substoeplalg}

We denote by ${\mathcal T}_p({\mathbf K})$ the algebra given by all
lower triangular 
elements of Toeplitz type in ${\mathbf K}^{{\mathcal M}_{p\times p}}$.
This algebra is isomorphic
to $\prod_{l=0}^\infty \left({\mathbf K}[x]
\pmod{x^{p^l}}\right)$ by Proposition \ref{isotoeppoly}. 
The subalgebra of converging elements
in ${\mathcal T}_p({\mathbf K})$ can thus be identified with
the algebra ${\mathbf K}[[x]]$ of formal power series.
We denote by ${\mathcal T}_{p-rec}({\mathbf K})=
{\mathcal T}_p({\mathbf K})\cap \hbox{Rec}_{p\times p}({\mathbf K})$ 
the subalgebra of recurrence matrices in ${\mathcal T}_p({\mathbf K})$.

Using the sequence of vectors 
$$A[{\mathcal M}_p^l]=(\alpha_0^l,\dots,\alpha_{p^l-1}^l)\in{\mathbf
  K}^{p^l}$$
defined by $A\in{\mathbf K}^{{\mathcal M}_p}$, we associate to $A\in
{\mathbf K}^{{\mathcal M}_p}$, the lower triangular element
$L_A\in{\mathbf K}^{{\mathcal M}_{p\times p}}$ of Toeplitz type
given by
$$L_A[{\mathcal M}_{p\times p}^l]=\left(\begin{array}{ccccc}
\alpha_0^l\\
\alpha_1^l&\alpha_0^l\\
\alpha_2^l&\alpha_1^l&\alpha_0^l\\
\vdots&&&\ddots\\
\alpha_{p^l-1}^l&\dots&&\alpha_1^l&\alpha_0^l\end{array}\right).$$

The map $A\longmapsto
L_A$ is clearly an isomorphism of vector spaces between ${\mathbf
  K}^{{\mathcal M}_p}$ and the algebra
${\mathcal T}_p({\mathbf K})$.

\begin{prop} \label{propisoToep}
The application $A\longmapsto L_A$ satisfies
$$\hbox{dim}(\overline{A}^{rec})\leq \hbox{dim}(\overline{L_A}^{rec})
\leq 2\ \hbox{dim}(\overline{A}^{rec})$$
and induces thus an isomorphism of vector spaces
between $\hbox{Rec}_p({\mathbf K})$ and the algebra
${\mathcal T}_{p-rec}({\mathbf K})\subset
\hbox{Rec}_{p\times p}({\mathbf K})$.
\end{prop}

\begin{cor} \label{corToepconstr} The map $\hbox{Rec}_p(\mathbf K)\times
  \hbox{Rec}_p(\mathbf K)\longrightarrow \hbox{Rec}_{p\times p}
(\mathbf K)$ defined by $(A,B)\longmapsto L_A+L_B^t$ with $L_A,L_B$ 
as above for $A,B\in\hbox{Rec}_p(\mathbf K)$ 
is a surjection onto the vector space of all recurrence
matrices of Toeplitz type with kernel the subspace of all pairs
$(A,B)\in \left(\hbox{Rec}_p(\mathbf K\right))^2$ such that
$A[U]=B[U]=0$ for $U\in\mathcal M_p\setminus \{0^*\}$
and $A[0^l]=-B[0^l]$ for all $l\in \mathbb N$.

In particular, an element of Toeplitz type in $\mathbf K^{\mathcal
  M_{p\times p}}$ is a recurrence matrix if and only if its
first row and column vectors are in $\hbox{Rec}_p(\mathbf K)$.
\end{cor}  

\begin{rem} The inequality 
$\hbox{dim}(\overline{A}^{rec})\leq \hbox{dim}(\overline{L_A}^{rec})$
in Proposition \ref{propisoToep}
is in general strict as illustrated by the example
$A\in{\mathbb Q}^{{\mathcal M}_{2\times 2}}$ with coefficients
$A[U]=1$ for all $U\in{\mathcal M}_2$. The recurrence vector $A\in
\hbox{Rec}_2({\mathbb Q})$ is thus of complexity $1$ while $L_A$ 
has complexity $2$ since $\rho(0,0)L_A=\rho(1,1)L_A
=L_A$ and $\rho(1,0)L_A$ 
(with coefficients $(\rho(1,0)L_A)[U,W]=1$ for all $(U,W)\in{\mathcal
  M}_{2\times 2}$) are linearly independent.
\end{rem}

{\bf Proof of Proposition \ref{propisoToep}} The inequality 
$\hbox{dim}(\overline{A}^{rec})\leq \hbox{dim}(\overline{L_A}^{rec})$
follows easily from the observation that the first column
of $L_A[{\mathcal M}_{p\times p}^l]$ is given by $A[{\mathcal M}_p^l]$.

For $A\in {\mathbf K}^{{\mathcal M}_p}$, we denote by $U_A\in
{\mathbf K}^{{\mathcal M}_{p\times p}}$ the strictly upper triangular 
matrix  of Toeplitz type defined by
$$U_A[{\mathcal M}_{p\times p}^l]=\left(\begin{array}{cccccccccc}
0&\alpha_{p^l-1}^l&\dots& \alpha_3^l&\alpha_2^l&\alpha_1^l\\
 &0&\alpha_{p^l-1}^l&\dots&\alpha_3^l&\alpha_2^l\\
 & &\ddots          &     &          &\alpha_3^l\\
 & &                &\ddots     &         & \\
 & &                &     &0        &\alpha_{p^l-1}^l\\
 & &                &      &         &0\end{array}\right)$$
where $A[{\mathcal
  M}_p^l]=(\alpha_0^l,\alpha_1^l,\dots,\alpha_{p^l-1}^l)$.

We consider now the vector spaces 
${\mathcal L}_{\overline{A}^{rec}}=\{L_A\ \vert\
A\in\overline{A}^{rec}\}$ and 
${\mathcal U}_{\overline{A}^{rec}}=\{U_A\ \vert\
A\in\overline{A}^{rec}\}$ spanned by all lower triangular,
respectively
strictly upper triangular elements of Toeplitz type in ${\mathbf
  K}^{{\mathcal M}_{p\times p}}$.

It is now straightforward to check
that ${\mathcal L}_{\overline{A}^{rec}}\oplus
{\mathcal U}_{\overline{A}^{rec}}
\subset{\mathbf K}^{{\mathcal
    M}_{p\times p}}$
is recursively closed. This shows the inequalities
$$\hbox{dim}(\overline{L_{\mathcal A}}^{rec})\leq
\hbox{dim}({\mathcal L}_{\overline{A}^{rec}}\oplus
{\mathcal U}_{\overline{A}^{rec}})\leq 
2\ \hbox{dim}(\overline{A}^{rec}).$$

The last part of Proposition \ref{propisoToep} follows from
the observation that the map $A\longmapsto L_A$ defines an isomorphism
of vector-spaces between $K^{{\mathcal M}_p}$ and 
${\mathcal T}_p({\mathbf K})$.\hfill $\Box$

The proof of Corollary \ref{corToepconstr} is immediate using 
the observation that the sequence $X[0^l,0^l]$ defines an
element of $\hbox{Rec}_1(\mathbf K)$ for $X\in\hbox{Rec}_{p\times
  p}(\mathbf K)$.

\subsection{The polynomial ring-structure on ${\mathbf K}^{{\mathcal
      M}_p}$ and $\hbox{Rec}_p({\mathbf K})$}\label{sectionpolynomialring}

The isomorphisms of vector spaces
${\mathbf K}^{{\mathcal M}_p}\sim {\mathcal T}_p\sim \prod_{l=0}^\infty 
\left({\mathbf K}[x]\mod {x^{p^l}}\right)$
given by the map $A\longmapsto L_A$ considered above and 
Proposition \ref{isotoeppoly} endow ${\mathbf K}^{{\mathcal  M}_p}$ 
with the {\it polynomial product}. More precisely, we consider the
map given by
$$A\longmapsto \psi(A)=\prod_{l=0}^\infty \psi_l(A)\in
\prod_{l=0}^\infty \left({\mathbf K}[x]\pmod{x^{p^l}}\right)$$
where $\psi_l(A)=\sum_{k=0}^{p^l-1}
\alpha_k^lx^k$ if
$A[{\mathcal M}_p^l]=(\alpha_0^l,\dots,\alpha_{p^l-1}^l)$
for $A\in {\mathbf K}^{{\mathcal M}_p}$.
We have then $\psi(C)=\psi(A)\psi(B)\in \prod_{l=0}^\infty 
\left({\mathbf K}[x]\mod {x^{p^l}}\right)$ if and only if 
$\psi_l(C)\equiv \psi_l(A)\psi_l(B)\pmod{x^{p^l}}$
for all $l\in{\mathbb N}$ or equivalently if and only if 
$L_C=L_AL_B\in{\mathbf K}^{{\mathcal M}_{p\times p}}$.

In particular, if $A$ and $B$ are convergent and correspond to 
formal power series $g_A,g_B$, then the polynomial product $\psi(C)=
\psi(A)\psi(B)$ corresponds to a convergent element $C\in{\mathbf
  K}^{{\mathcal M}_p}$ with associated 
formal power series $g_C=g_Ag_B$ defined as the usual 
product of the formal power series $g_A$ and $g_B$.

Using the isomorphism $A\longmapsto L_A$ 
of the previous section, we identify
the polynomial product algebra ${\mathbf K}^{{\mathcal M}_p}$
with the subalgebra ${\mathcal T}_p({\mathbf K})\subset {\mathbf
  K}^{{\mathcal M}_{p\times p}}$ of lower triangular 
matrices of Toeplitz type. Similarly, we identify
$\hbox{Rec}_p({\mathbf K})$ with the commutative subalgebra 
${\mathcal T}_{p-rec}({\mathbf K})\subset \hbox{Rec}_{p\times
  p}({\mathbf K})$.

We denote by $\mathcal{IT}_p\subset {\mathcal T}_p({\mathbf K})$ 
the subspace corresponding to
all elements $A\in {\mathbf K}^{{\mathcal M}_p}$
such that $\rho(s)A$ is of finite support for $s=0,1,\dots p-2$
and $\rho(p-1)A\in \mathcal{IT}_p$ (where we identify
${\mathbf K}^{{\mathcal M}_p}$ and ${\mathcal T}_p({\mathbf K})$
using the map $A\longmapsto L_A$). Although this definition 
is recursive, it makes sense: An element $A\in\mathcal{IT}_p$ corresponds
to a sequences $\psi_0(A),\psi_1(A),\dots $ of polynomials where $
\psi_k(A)\in{\mathbf K}[x]$ is
of degree $< p^l$ and almost all polynomials
$\psi_l(A)$ satisfy $\psi_l(A)\equiv 0\pmod{x^{p^l-p^k}}$ 
for every fixed natural integer $k\in {\mathbb N}$.

It is easy to check that $\mathcal{IT}_p$ is an ideal of the 
algebra ${\mathcal T}_p({\mathbf K})$. We denote by 
$\tilde{\mathcal T}_p({\mathbf K})={\mathcal T}_p({\mathbf K})/
\mathcal {IT}_p$ and
$\tilde{\mathcal T}_{p-rec}({\mathbf K})={\mathcal T}_{p-rec}
({\mathbf K})/\mathcal
{IT}_{p-rec}$ (where $\mathcal{IT}_{p-rec}=\mathcal{IT}_p\cap
{\mathcal T}_{p-rec}({\mathbf K})$) the obvious quotient algebras.
It is also easy to check that the differential operator
$\frac{d}{dx}$ (acting in the obvious way on the
polynomial sequence $\psi(A)=(\psi_0(A),\psi_1(A),\dots)$
associated to $A\in{\mathbf K}^{{\mathcal M}_p}$)
is well-defined on the quotient algebra 
$\tilde{\mathcal T}_p({\mathbf K})$.

\begin{prop} \label{propdiffring}
The algebra $\tilde{\mathcal T}_{p-rec}({\mathbf K})$ 
is a differential
subalgebra of the differential algebra $\tilde{\mathcal T}_{p}
({\mathbf K})$.

In particular, converging recurrence vectors correspond to 
a differential subring of the differential ring 
$({\mathbf K}[[x]],\frac{d}{dx})$ of formal power-series.
\end{prop}

{\bf Proof} Consider the factorization $\frac{d}{dx}=
\frac{1}{x}\left(x\frac{d}{dx}\right)$
 of the differential operator $\frac{d}{dx}$ into the differential 
operator $(x\frac{d}{dx})$, followed by multiplication by $x^{-1}$.
Given an element $A\in{\mathbf K}^{{\mathcal M}_p}$ (identified
in the obvious way with the corresponding sequence 
$\psi(A)=\prod_{l=0}^\infty \psi_l(A)$ of polynomials where
$\psi_l(A)\in{\mathbf K}[x]$ is of degree $<p^l$),
the differential operator $\psi(A)\longmapsto (x\frac{d}{dx})\psi(A)$ 
corresponds to the map $A\longmapsto B=NA$ multiplying a row-vector
$A\in {\mathbf K}^{{\mathcal M}_{p\times 1}}$ on the left with the
converging recurrence matrix $N$ with limit the diagonal matrix
having diagonal entries $0,1,2,3,\dots$. Multiplication of $\psi(B)$
by $^{-1}$ corresponds now to the map 
$B\longmapsto C=PB$ where 
$$P[{\mathcal M}_{p\times p}^l]=\left(\begin{array}{ccccc}
0&1\\
 &0&1\\
\vdots&&\ddots&\ddots\\
1&&&&0\end{array}\right)$$
is the obvious cyclic permutation of order $p^l$. Since
$N$ and $P$ are both of complexity $2$ in $\hbox{Rec}_{p\times p}(
{\mathbf K})$ this shows the inequality
$$\hbox{dim}(\overline{\frac{d}{dx}A}^{rec})\leq 4\ 
\hbox{dim}(\overline{A}^{rec})$$
for $A\in \hbox{Rec}_p({\mathbf K})$ and implies
the first part. 

The second part follows from the obvious observation that the
quotient map ${\mathcal T}_p({\mathbf K})\longmapsto
\tilde {\mathcal T}_p({\mathbf K})$ restricts to an injection
on the subalgebra of converging elements in 
$\mathcal{T}_{p}({\mathbf K})$.\hfill$\Box$

\begin{rem} \label{remholToepl}
Over the field ${\mathbf K}={\mathbb C}$ of complex numbers
and for $p\geq 2$, the formal power series $g_A$ associated to 
a converging recurrence vector $A\in\hbox{Rec}_p({\mathbb C})$
defines by Proposition \ref{majorationcriterion}
a holomorphic function in the open unit disc.
\end{rem}

\section{Elements of Hankel type}

A (finite or infinite) matrix $H$ with coefficients $h_{s,t},\ 0\leq s,t$
is a {\it Hankel matrix} if $h_{s,t}=\alpha_{s+t}$ for some 
sequence $\alpha_0,\alpha_1,\alpha_2,\dots$. An element
$H\in{\mathbf K}^{{\mathcal M}_{p\times p}}$ is {\it of 
Hankel type} if for all $l\in{\mathbb N}$, the matrix 
$H[{\mathcal M}_{p\times p}^l]$ is a 
Hankel matrix (for the usual total order
on rows and columns induced by $(u_1\dots u_l)\longmapsto 
\sum_{j=1}^l u_j p^{j-1}$). 

Consider the involutive element 
$S\in\hbox{Rec}_{p\times p}({\mathbf K})$ where the coefficients of
$S[{\mathcal M}_{p\times p}^l]$ are $1$ on the antidiagonal 
and zero elsewhere. More precisely, 
$S[u_1\dots u_n,w_1\dots w_n]=1$ if $w_1=p-1-u_1,\dots,w_n=p-1-u_n$
and $S[U,W]=0$ otherwise. (In particular,
$S$ is a recurrence matrix of Hankel type with complexity $1$.) 
The following result reduces the study of Hankel matrices 
to the study of Toeplitz matrices.

\begin{prop} (Left- or right-)multiplication by the involution 
$S\in\hbox{Rec}_{p\times p} ({\mathbf K})$
defines a bijection preserving complexities between elements of
Hankel type and elements of Toeplitz type
in ${\mathbf K}^{{\mathcal M}_{p\times p}}$. In particular, the maps
$A\longmapsto SA$ and $A\longmapsto AS$ yield bijections between
recurrence matrices of Hankel type and recurrence matrices of 
Toeplitz type in $\hbox{Rec}_{p\times p}({\mathbf K})$.
\end{prop}

The easy proof is left to the reader.

Elements of Hankel type satisfy $H[U,W]=H[W,U]$ for all $(U,W)\in
{\mathcal M}_{p\times p}$ and are thus examples of symmetric
elements in ${\mathbf K}^{{\mathcal M}_{p\times p}}$.

Let $H\in{\mathbf K}^{{\mathcal M}_{p\times p}}$ be an element
of Hankel type. Since $\rho(s,t)H,\ 0\leq s,t<r$
is of Hankel type, 
the recursive closure $\overline{H}^{rec}$ contains only elements
of Hankel type. The following result
is the exact analogue of Proposition \ref{Toeplitzprop}:

\begin{prop} \label{Hankelprop} All recurrence matrices
  $H_1,\dots,H_d\in\hbox{Rec}_{p\times p}
({\mathbf K})$ of a reduced presentation are of Hankel 
type if and only if the shift matrices $\rho(s,t)$ 
satisfy the following two conditions:

\ \ (1) The shift matrices $\rho(s,t),\ 0\leq s,t<p$ depend only 
on $s+t$.

\ \ (2) Assuming condition (1) and writing 
$\rho(s+t)$ for a shift matrix $\rho(s,t)$, we have
the identities
$$\rho(s+p)\rho(t)=\rho(s)\rho(t+1)$$
for $0\leq s\leq p-2$ and $0\leq t\leq 2p-3$.
\end{prop}

For $p=2$, Proposition \ref{Hankelprop} boils down to 
the identities $\rho(1,0)=\rho(0,1)$ and 
$$\begin{array}{l}
\displaystyle \rho(0,0)\rho(1,1)=\rho(1,1)\rho(0,1),\\
\displaystyle \rho(1,1)\rho(0,0)=\rho(0,0)\rho(0,1)\end{array}$$
for shift-matrices of a presentation containing only recurrence
matrices of Hankel type.

{\bf Proof of Proposition \ref{Hankelprop}} We show by induction on
$l$ that all matrices $H_i[{\mathcal M}_{p\times p}^l]$ are
Hankel matrices.

For $l=0$, there is nothing to do. Condition (1) implies that
$H_1[{\mathcal M}_{p\times p}^1],\dots,H_d[{\mathcal
  M}_{p\times p}^1]$ are Hankel matrices. 
Setting $H=H_i$ and 
denoting $(\rho(s,t)H)[{\mathcal M}_{p\times p}^l]$ by $\rho(s+t)H$ we get
$$H[{\mathcal M}_{p\times p}^{l+1}]=\left(\begin{array}{ccccc}
\rho(0)H&\rho(1)H&\rho(2)H&\cdots&\rho(p-1)H\\
\rho(1)H&\rho(2)H&\rho(3)H&\cdots&\rho(p)H\\
\vdots&          &        &      &   \vdots\\
\rho(p-1)H&\rho(p)H&\rho(p+1)H&\cdots&\rho(2p-2)H\end{array}\right)$$
where all matrices 
$\rho(0)H,\dots,\rho(2p-2)H$ are Hankel 
matrices by induction. 
We have to show that two horizontally adjacent
blocks $(\rho(t)H\ \rho(t+1)H),\ 0\leq t<2p-2$ 
define a $p^l\times (2\cdot p^l)$
matrix of ``Hankel-type'' (the argument for vertically adjacent blocks
gives rise to the same conditions and is left to the reader).
We have
$$\big(\rho(t)H\vert \rho(t+1)H\big)=\left(\begin{array}{cr|lc}
\cdots&\rho(p-1)\rho(t)H&\rho(0)\rho(t+1)H&\cdots\\
\cdots&\rho(p)\rho(t)H&\rho(1)\rho(t+1)H&\cdots\\
\cdots&\rho(p+1)\rho(t)H&\qquad\vdots&\\
&\vdots\qquad&\qquad\vdots&\\
&\vdots\qquad&\rho(p-2)\rho(t+1)H&\cdots\\
\cdots&\rho(2p-2)\rho(t)H&\rho(p-1)\rho(t+1)H&\cdots\end{array}\right)$$
(the vertical line separates the matrix into two square blocks)
using the shorthand notation $\rho(\alpha)\rho(\beta)H=
(\rho(\alpha)\rho(\beta)H)[{\mathcal M}_{p\times p}^{l-1}]$.
Such a matrix is of Hankel-type if
$$\rho(s+p)\rho(t)H=\rho(s)\rho(t+1)H\hbox{ for } 0\leq s\leq p-2\ .$$

The opposite direction is left to the reader.\hfill $\Box$

Given two vectors $A,B\in\mathbf K^{\mathcal M_{p\times p}}$ 
with coefficients
$$\begin{array}{l}
\displaystyle A[\mathcal
M_p^l]=\left(\alpha_0^l,\dots,\alpha_{p^l-1}^l\right)\in\mathbf
K^{p^l}\\
\displaystyle B[\mathcal
M_p^l]=\left(\beta_0^l,\dots,\beta_{p^l-1}^l\right)\in\mathbf
K^{p^l}\end{array}$$
we consider the Hankel matrices $H_A,\tilde H_B\in\mathbf K^{\mathcal
  M_{p\times p}}$
with coefficients 
$$H_A[\mathcal M_{p\times p}^l]=\left(\begin{array}{ccccc}
\alpha_0^l&\alpha_1^l&\dots&\alpha_{p^l-1}^l\\
\alpha_1^l&\alpha_2^l&&0\\
\vdots&&&\vdots\\
\alpha_{p^l-1}^l&0\dots&0\end{array}\right)\ ,$$
$$\tilde H_B[\mathcal M_{p\times p}^l]=\left(\begin{array}{cccccccc}
0&0&0&\dots&0&0\\
0& & &     &0&\beta_0^l\\
\vdots&&&&&\vdots\\
0&0&\beta_0^l&&\beta_{p^l-4}^l&\beta_{p^l-3}^l\\
0&\beta_0^l&\beta_1^l&\dots&\beta_{p^l-3}^l&\beta_{p^l-2}^l\end{array}
\right)\ .$$

The following results are close analogues of Proposition 
\ref{propisoToep} and Corollary \ref{corToepconstr}.

\begin{prop} \label{propHankelvector}
The applications $A\longmapsto H_A,\tilde H_A$ defined above satisfy
the inequalities
$$\hbox{dim}(\overline{A}^{rec})\leq
\hbox{dim}(\overline{H_A}^{rec})$$
and 
$$\hbox{dim}(\overline{H_A}^{rec}),\hbox{dim}(\overline{
\tilde H_A}^{rec})\leq 2\ \hbox{dim}(\overline{A}^{rec})\ .$$
\end{prop}

\begin{cor} \label{corHankelvect}
The application $\hbox{Rec}_p(\mathbf K)\times
  \hbox{Rec}_p(
\mathbf K)\longrightarrow \hbox{Rec}_{p\times p}(\mathbf K)$
defined by $(A,B)\longmapsto H_A+\tilde H_B$ with $H_A,\tilde H_B
\in \hbox{Rec}_{p\times p}(\mathbf K)$
as above for $A,B\in\hbox{Rec}_p(\mathbf K)$ is a surjection onto the 
vector space of all Hankel matrices in $\hbox{Rec}_{p\times p}(\mathbf
K)$ with kernel spanned by $(A,B)\in\left(\hbox{Rec}_p(\mathbf
  K)\right)^2$ such that $A=0$ and $B[U]=0$ for 
$U\in \mathcal M_p\setminus \{(p-1)^*\}$.

In particular, a Hankel matrix in $\mathbf K^{\mathcal M_{p\times p}}$
is a recurrence matrix if and only if its first column vector and
its last row vector are elements of $\hbox{Rec}_p(\mathbf K)$.
\end{cor}

\begin{rem} The ``missing'' inequality 
$$\hbox{dim}(\overline{A}^{rec})\leq
\hbox{dim}(\overline{\tilde H_A}^{rec})$$
in Proposition \ref{propHankelvector} does not necessarily hold. 
Its possible failure is due to the fact that the last coefficient of
$A[\mathcal M_p^l]$ is not involved in $\tilde H_A[\mathcal M_{p\times
  p}^l]$.
\end{rem}

{\bf Proof of Proposition \ref{propHankelvector}} Check that 
the vector space spanned by all elements $H_X,\tilde H_X$ 
for $X\in\overline{A}^{rec}$ is recursively closed.
\hfill$\Box$

The proof of Corollary \ref{corHankelvect} is obvious.

\subsection{Hankel matrices and continued fractions of Jacobi type}

An infinite Hankel matrix $H$ associated to the generating function
$\gamma=\sum_{n=0}^\infty \gamma_n x^n\in\mathbf K[[x]]$
(where we suppose $\gamma_0=1$ in order to simplify subsequent
statements)
is {\it non-degenerate} if all finite
Hankel matrices $H(n)$ with coefficients $H(n)_{i,j}=\gamma_{i+j},
0\leq i,j<n$ are invertible. It is then well-known (see for instance
\cite{Flajolet}) that we have a continued fraction of Jacobi type
$$\gamma=\frac{1}{1-\alpha_0x-\beta_0\frac{x^2}{1-\alpha_1x-\beta_1
\frac{x^2}{1-\dots}}}$$
for $\gamma=1+\gamma_1x+\dots$ 
where the sequences $\alpha_0,\alpha_1,\dots,
\beta_0,\beta_1,\dots$ appear also in the recursive definition
$$p_{n+1}=(x-\alpha_n) p_n-\beta_np_{n-1}$$ 
of the monic orthogonal polynomials $p_0=1,\dots$ for $H$.
We have then the following result.

\begin{prop} \label{propJacobifraction} Let $H\in\hbox{Rec}_{p\times
    p}(\mathbf K)$ be a converging non-degenerate Hadamard
matrix such that $H=LDL^t$ with $L\in \hbox{Rec}_{p\times p}(\mathbf K)$
lower triangular unipotent and invertible in 
$\hbox{Rec}_{p\times p}(\mathbf K)$. Then the sequences
$\alpha_0,\alpha_1,\dots$ and $\beta_0,\beta_1,\dots$ involved in
the continued fraction of Jacobi type associated to $H$
define converging elements of $\hbox{Rec}_{p}(\mathbf K)$.
\end{prop}

{\bf Proof} Since $L^{-1}H\left(L^{-1}\right)^t=D$
is a diagonal matrix with invertible diagonal coefficients,
the $n-$th row $(m_{n,0},\dots,m_{n,n},0,\dots$ of 
$M=L^{-1}$ (identified with its limit) consists of the coefficients
of the $n-$th orthogonal polynomial $p_n=\sum_{k=0}^n m_{n,k}x^k$
for $H$. Since $\langle xp_n,p_k\rangle=\langle p_n,xp_k\rangle$,
there exists constants $\alpha_n,\beta_n$ such that
we have the identities $$xp_n=p_{n+1}+\alpha_n p_n+\beta_n p_n$$
which can be written in matrix-form as
$$\tilde M=\left(\begin{array}{ccccc}
\alpha_0&1\\\beta_1&\alpha_1&1\\
&\beta_2&\alpha_2&1\\
&&&\ddots\end{array}\right)M$$
where $\tilde M$ is defined by adding a first zero column
to the infinite lower triangular matrix $M=L^{-1}$.
The tridiagonal matrix 
$$S=\left(\begin{array}{ccccc}
\alpha_0&1\\\beta_1&\alpha_1&1\\
&\beta_2&\alpha_2&1\\
&&&\ddots\end{array}\right)$$
is called the {\it Stiltjes matrix} of $H$ and a small computation
shows that it is also given by $S=L^{-1}\tilde L$ where $\tilde L$
is obtained by removing the first row of the infinite
lower triangular matrix $L$. For $L$ invertible 
in $\hbox{Rec}_{p\times p}(\mathbf K)$
both matrices $\tilde L$ and $L^{-1}$ are in $\hbox{Rec}_{p\times
  p}(\mathbf K)$. This shows $S\in\hbox{Rec}_{p\times p}(\mathbf K)$
and implies the result.\hfill$\Box$

\section{Groups of recurrence matrices}

The (multiplicative) identity $Id$ of the algebra ${\mathbf K}^{
{\mathcal M}_{p\times p}}$ satisfies
$\rho(s,s)Id=Id,0\leq s <p$ and $\rho(s,t)Id=0, 0\leq s\not= t<p$,
and is thus an element of complexity $1$ in $\hbox{Rec}_{p\times
  p}({\mathbf K})$.  

\begin{defin} We denote by
$\hbox{GL}_{p-rec}({\mathbf K})$ the {\it general linear group of all 
recurrence matrices} in $\hbox{Rec}_{p\times p}({\mathbf K})$ 
which are invertible (for the matrix-product)
in $\hbox{Rec}_{p\times p}({\mathbf K})$.

A {\it group of recurrence matrices} is a subgroup of 
$\hbox{GL}_{p-rec}({\mathbf K})$ for some $p\in {\mathbb N}$.
\end{defin}

\begin{prop} \label{groupconstr} The set
  $\left(\hbox{GL}_{p-rec}({\mathbf K})\right)^p\times
\left(\hbox{Rec}_{p\times p}({\mathbf K})\right)^{p\choose 2}$
injects into $\hbox{GL}_{p-rec}({\mathbf K})$.
\end{prop}

{\bf Proof} For $A_0,\dots,A_{p-1}\in \hbox{GL}_{p-rec}({\mathbf K})$
and $B_{s,t}\in\hbox{Rec}_{p\times p}({\mathbf K}), 0\leq t<s<p$, 
consider the matrix $G\in \hbox{Rec}_{p\times p}({\mathbf K})$ 
defined by $\rho(s,s)G=A_s$, $\rho(s,t)G=B_{s,t}$ for $0\leq t<s<p$
and $\rho(s,t)G=0$ for $s<t$ and make an arbitrary choice in 
${\mathbf K}^*$ for $G[\emptyset,\emptyset]$. This defines an element
$G\in\hbox{Rec}_{p\times p}({\mathbf K})$ which is invertible.
\hfill$\Box$

For $p\geq 2$, the group $\hbox{GL}_{p-rec}(\mathbf K)$ 
is complicated: Proposition
\ref{groupconstr} shows already that it is very huge. One can also 
show that it contains every finite group and every finite-dimensional
matrix group over $\mathbf K$. Moreover, the following result shows
that there is (in general) no relation between the complexity 
of an element $G\in\hbox{GL}_{p-rec}(\mathbf K)$ and the complexity
of its inverse element $G^{-1}\in   \hbox{GL}_{p-rec}(\mathbf K)$.

\begin{prop} \label{propnobound}
For $p>1$ and $N\in\mathbb N$ arbitrary, the group
$\hbox{GL}_{p-rec}(\mathbb C)$ contains an element $G$ of
complexity $2$ with inverse of complexity $\geq N$.
\end{prop}

{\bf Proof} For a natural integer $p>1$, the determinant of the
symmetric $N\times N-$Hankel matrix 
$$\left(\begin{array}{ccccc}x^{p}&x^{p^2}&x^{p^3}&\dots&x^{p^N}\\
x^{p^2}& &&&\vdots\\
\vdots&&&&\vdots\\
x^{p^N}&&&&x^{p^{2N-1}}\end{array}\right)$$
is a monic polynomial $P$ of degree
$p+p^2+p^{2N-1}=\frac{p^{2N}-p}{p-1}>0$.
Choose $n\in\mathbb N$ such that $\omega=e^{2i\pi/n}\in \mathbb C$
is not a root of $P$ and consider the converging element 
$G\in\hbox{Rec}_{p\times p}(\mathbb C)$ defined by the lower
triangular 
Toeplitz matrix 
$$\left(\begin{array}{cccc}
1\\-\omega&1\\
0&-\omega&1\\
&&&\ddots\end{array}\right)$$
associated to $1-\omega x\in\mathbb C[x]\subset \mathbb C[[x]]$
as in section \ref{sectionpolynomialring}. 
It is easy to check that $G$ is of
complexity $2$. Since $\omega$ is of finite order,
$G$ admits an inverse element $G^{-1}\in\hbox{Rec}_{p\times p}(\mathbb
C)$, given by the 
converging lower triangular Toeplitz matrix with limit
$$\left(\begin{array}{ccccc}
1\\\omega&1\\
\omega^2&\omega&1\\
\omega^3&\omega^2&\omega&1\\
&&&&\ddots\end{array}\right)$$
associated to $\sum_{j=0}^\infty (\omega x)^j\in\mathbb C[[x]]$.
For $l\geq 1$, we consider the sequence 
$$\begin{array}{rcl}
S(l)&=&((\rho(1^l,0^l)G^{-1})[\emptyset,\emptyset],
 (\rho(1^{l},0^{l})G^{-1})[1,0], (\rho(1^{l},0^{l})
G^{-1})[1^2,0^2],\dots)\\
&=&(\omega^{p^l},\omega^{p^{l+1}},\omega^{p^{l+2}},\dots)\end{array}\ .$$
The submatrix 
$$\left(\begin{array}{cccccccccc}
\omega^p&\omega^{p^2}&\omega^{p^3}&\dots\\
\omega^{p^2}&\omega^{p^4}&\omega^{p^3}&\dots\\
\vdots\end{array}\right)$$
formed by the $N$ initial coefficients of the
sequences $S(1),\dots,S(N)$ has non-zero determinant by
construction. This shows that $G^{-1}$ has complexity $\geq N$.
\hfill$\Box$

\begin{rem} Proposition \ref{propnobound} holds also for instance
for the algebraic closure of finite fields. It is however 
trivially wrong (for $p$ fixed) over
  finite fields: Since there are only finitely many elements
of complexity $\leq b$ in $\hbox{Rec}_{p\times p}(\mathbb F)$
for $\mathbb F$ a finite field, the number of invertible 
elements in $\hbox{Rec}_{p\times p}(\mathbb F)$
having complexity $\leq b$ is also finite and there is
thus an upper bound for the complexities 
of their inverses.
\end{rem}

\begin{rem} \label{remnoinverse} An element $X\in \hbox{Rec}_{p\times p}({\mathbf K})$ 
having a multiplicative inverse $X^{-1}\in {\mathbf K}^{{\mathcal
      M}_{p\times p}}$ is not necessarily in
$\hbox{GL}_{p-rec}({\mathbf K})$ since $X^{-1}$ has in general
infinite complexity. A simple example is the ``diagonal'' matrix 
$X\in {\mathbb Q}^{{\mathcal M}_{p\times p}}$ which is in the
center of the algebra ${\mathbb Q}^{{\mathcal M}_{p\times p}}$
and has diagonal coefficients $X[U,U]=n+1
\in\{1,2,\dots\}$ depending only on the length $n$
of $(U,U)\in{\mathcal M}_{p\times q}^n$ (and off-diagonal coefficients
$X[U,W]=0$ for $U\not=W$ with $(U,W)\in{\mathcal M}_{p\times p}$). 
We have $\rho(s,t)X=0$ for $0\leq s\not= t<p$ and 
$\rho(s,s)X=X+Id$ where $Id\in\hbox{Rec}_{p\times p}({\mathbb Q})$ 
denotes the multiplicative identity.
The element $X$ has thus complexity $2$ and is in $\hbox{Rec}_{
p\times p}({\mathbf K})$. Since the denominators of
$X^{-1}$ involve all primes of
${\mathbb Z}$, it cannot be a recurrence matrix by 
Proposition \ref{critfingenring} or \ref{critfinprimes}.
A second proof of this fact 
is given by the observation that the infinite matrix $H$
with coefficients
$c_{i,j}=(\rho(0^i,0^i)X^{-1})[0^j,0^j]=\frac{1}{i+j+1},0\leq i,j$
(using the convention $0^0=\emptyset$) is the Hilbert matrix
$$
\left(\begin{array}{ccccc}
\frac{1}{1}&\frac{1}{2}&\frac{1}{3}&\dots\\
\frac{1}{2}&\frac{1}{3}&\frac{1}{4}&\dots\\
\frac{1}{3}&\frac{1}{4}&\frac{1}{5}&\dots\\
\vdots&\vdots&\vdots&\ddots
\end{array}
\right).$$
The finite submatrix formed by the first $n$ rows and columns of 
$H$ is non-singular for all $n$. This implies that 
$X^{-1}$ is of infinite complexity.

Another simple example for $p\geq 2$ is the converging lower 
triangular unipotent recurrence matrix $L\in\hbox{Rec}_{2\times 2}({
\mathbb Q})$ of Toeplitz type with associated generating series
$1-2x$. Its coefficients are $1$ on the diagonal, 
$-2$ on the subdiagonal
and $0$ everywhere else. Its inverse $L^{-1}\in{\mathbb Q}^{{\mathcal
    M}_{2\times 2}}$ is lower triangular of Toeplitz type 
with first row given by successive powers of 2
and corresponds to the generating series $\frac{1}{1-2x}$. 
Proposition \ref{majorationcriterion} shows that $L^{-1}$ cannot
be a recurrence matrix (see also Remark \ref{remholToepl}).
However, the element $L^{-1}$ 
is a recurrence matrix over 
a field of positive characteristic. (In the case of odd
characteristic $\wp$, its complexity depends on the 
multiplicative order of $2\in\left({\mathbb Z}/\wp{\mathbb Z}
\right)^*$.)
\end{rem}

\begin{rem}
For ${\mathbf K}$ a topological field, the group 
of invertible elements in ${\mathbf K}^{{\mathcal M}_{p\times p}}$
is a topological group for both topologies considered in \S 
\ref{topologies}. The subgroup 
 $\hbox{GL}_{p-rec}({\mathbf K})$ is thus also a topological group (with
 discrete topology for the second topology described in
\S \ref{topologies}).

It would be interesting to have answers to the following questions: 
Does $\hbox{GL}_{p-rec}({\mathbf K})$ admit other,
more interesting topologies? Is the map
$X\longmapsto X^{-1}$ of invertible elements in the quotient
algebra $\hbox{Rec}_{p\times p}({\mathbf
  K})/\mathcal{FS}_{p\times p}$ continuous for the metric 
$\parallel X\parallel_\infty^\infty$ described in 
\S \ref{topologies}?
\end{rem}

\subsection{A few homomorphisms and characters}

Given an element $A\in \hbox{GL}_{p-rec}({\mathbf K})$,
the projection
$$A\longmapsto A[{\mathcal M}_{p\times p}^l]\in \hbox{GL}(
{\mathbf K}^{p^l})$$
yields a homomorphism of groups. All these homomorphisms
are surjective and have sections.
In particular, $\hbox{GL}_{p-rec}({\mathbf K})$
contains (an isomorphic image of) 
every finite-dimensional matrix-group over ${\mathbf K}$ for $p\geq 2$.

The {\it determinant} is the homomorphism
$\hbox{GL}_{p-rec}({\mathbf K})\longrightarrow 
\prod_{l=0}^\infty {\mathbf K}^*$ defined by
$$A\longmapsto\det(A)=(\det(A[{\mathcal M}_{p\times p}^0]),\det(A[{\mathcal M}_{p\times p}^1]),\det(A[{\mathcal M}_{p\times p}^2]),\dots)\in\prod_{l=0}^\infty 
{\mathbf K}^*\ .$$
Its image is in general not surjective: For a finite (countable)
field ${\mathbf K}$ there exist only finitely (countably) many elements 
of complexity $\leq d$ in $\hbox{End}_{p-rec}({\mathbf K})$ and the group
$\hbox{GL}_{p-rec}({\mathbf K})$ is thus countable. If ${\mathbf K}^*$
contains at least two elements, the subgroup
$\det(\hbox{GL}_{p-rec}({\mathbf K}))$ is thus a proper subgroup of the 
uncountable group $\prod_{l=0}^\infty {\mathbf K}^*$.

\begin{rem}
It would probably be interesting to have a description of
the countable subgroup $$\det(\hbox{GL}_{p-rec}({\mathbf K}))\subset
\prod_{l=0}^\infty {\mathbf K}^*\ .$$
An obvious restriction is for instance given by the trivial remark that 
all numerators and denominators of $\hbox{det}(A)$ involve only 
a finite number of distinct prime-factors for 
$A\in\hbox{GL}_{p-rec}({\mathbb Q})$  
(see Proposition \ref{critfinprimes}). The largest subgroup of 
$\prod_{l=0}^\infty {\mathbb Q}^*$ with this property is however 
still uncountable.
\end{rem} 

Setting
$\tau_n(A)=\sum_{U\in{\mathcal M}_p^n}A[U,U]$ 
we get the {\it trace-map} $A\longmapsto
\hbox{tr}(A)=(\tau_0(A),\tau_1(A),\dots,\dots)$ defining a linear 
application 
$$\hbox{tr}:\hbox{Rec}_{p\times p}({\mathbf K})\longrightarrow
\hbox{Rec}_{1\times 1}({\mathbf K}).$$
The trace map satisfies the identity
$\hbox{tr}(AB)=\hbox{tr}(BA)$ and is
thus constant on conjugacy classes under the action of 
$\hbox{GL}_{p-rec}({\mathbf K})$. The trace
$\hbox{tr}(A)\in\hbox{Rec}_{1\times 1}({\mathbf K})$
of a recurrence matrix
$A\in \hbox{Rec}_{p\times p}({\mathbf K})$ is easy to compute 
from a presentation
$(A_1,\dots,A_{a})[\emptyset,\emptyset],\rho_{\overline{A}^{rec}}(s,t),
\ 0\leq s,t<p$: 
It admits the 
presentation $(\hbox{tr}(A_1),\dots,\hbox{tr}(A_{a}))[\emptyset,
\emptyset]=(A_1,\dots,A_{a})[\emptyset,
\emptyset]$ with shift-matrix 
$\sum_{s=0}^{p-1} \rho_{\overline{A}^{rec}}(s,s)$. 
The associated generating series
$$\sum_{n=0}^\infty \tau_n(A) t^n\in {\mathbf K}[[t]]$$
is always a rational function.

\begin{rem}
The properties of the traces $\hbox{tr}(A^0),\hbox{tr}(A^1),
\hbox{tr}(A^2),\hbox{tr}(A^3),\dots$ are essentially shared by
the coefficients of $x^{p^l-k}$ of the characteristic polynomial
of $A[{\mathcal M}_{p\times p}^l]$.

This suggests the following question: What can be said of the
spectra of the matrices $A[{\mathcal M}_{p\times p}^l],l=0,1,2,\dots$?
Is there sometimes (perhaps after a suitable normalization)
a ``spectral limit'', eg. a nice spectral measure, etc?
(The answer to the last question is trivially yes for
recurrence matrices of complexity $1$ since they are tensor-powers,
up to a scalar factor.)
\end{rem}

\subsection{A few other properties}

Recall that a group $\Gamma$ is {\it residually finite} if for every
element $\gamma\in \Gamma$ different from the identity
there exists a homomorphism
$\pi:\Gamma\longrightarrow F$ into a finite group $F$ such that 
$\pi(\gamma)\not=1\in F$.

\begin{prop} \label{residfinite}
(i) For ${\mathbf K}$ a finite field, the group
$\hbox{GL}_{p-rec}({\mathbf K})$ is residually finite.

\ \ (ii) For ${\mathbf K}$ an arbitrary field, a finitely
generated group $\Gamma\subset
GL_{p-rec}({\mathbf K})$ is residually finite.
\end{prop}

{\bf Proof.} Given a finite field ${\mathbf K}$ and an element 
$1\not=A\in  \hbox{GL}_{p-rec}({\mathbf K})$, there exists $l$ such that 
$A[{\mathcal M}_{p\times p}^l]\not= 1\in 
\hbox{Aut}({\mathbf K}^{p^l})$ thus proving assertion (i).

For proving assertion (ii),
we remark that presentations for a finite set of generators
$\gamma_1^{\pm 1},\dots,\gamma_m^{\pm 1}$ generating
$\Gamma\subset \hbox{GL}_{p-rec}({\mathbf K})$ 
involve only finitely many elements
$k_1,\dots,k_N\in {\mathbf K}$ (appearing in the initial data
and as coefficients of the shift-matrices). The group $\Gamma$ 
is thus defined over the finitely generated field extension
${\mathbf F}$ containing $k_1,\dots,k_N$ over the primary
field ${\mathbf k}$ (which is either ${\mathbb Q}$ or the finite 
primary field ${\mathbb F}_{\wp}$ with $\wp$ elements for $\wp$ a
prime number) 
of ${\mathbf K}$. Given an element $\gamma\not=\hbox{1d}$ in $\Gamma$,
choose $l\in{\mathbb N}$ such that $\gamma[{\mathcal M}_{p\times p}^l]
\not=\hbox{1d}$.
Choose a transcendental basis
$t_1,\dots,t_f$ of ${\mathbf F}$ such that ${\mathbf F}=\tilde 
{\mathbf F}(t_1,\dots,t_f)$ where $\tilde {\mathbf F}$ is the maximal
subfield of ${\mathbf F}$ which is algebraic over ${\mathbf k}$.
Choose now a maximal ideal 
${\mathcal J}\subset {\mathbf A}=\tilde
{\mathbf F}[t_1,\dots,t_f]$ such that
$\Gamma$ is defined over ${\mathbf A}/{\mathcal J}$ and 
$\gamma[{\mathcal M}_{p\times p}^l]\not\equiv 
\hbox{1d}\pmod {\mathcal J}$. The quotient
field ${\mathbf A}/{\mathcal J}$ is either finite or a number field.
Reducing modulo a suitable prime $\wp\in {\mathbf A}/{\mathcal J}$ 
in the case of a number field we get a quotient group
(with $\overline \gamma\not=1$ in the quotient)
of $p^l\times p^l-$matrices defined over a finite field.\hfill $\Box$

\begin{rem} \label{reductionmodp}
It follows from Proposition \ref{critfinprimes} or from the proof of 
Proposition \ref{residfinite}, that any finitely generated group
$\Gamma\subset \hbox{GL}_{p-rec}({\mathbb Q})$ (or more generally
$\Gamma\subset \hbox{GL}_{p-rec}({\mathbf K})$ where ${\mathbf K}$
is a number field) can be reduced modulo $\wp$ to a quotient-group
$\Gamma_\wp\subset \hbox{GL}_{p-rec}({\mathbb Z}/\wp{\mathbb Z})$
for almost all primes $\wp$ (or prime-ideals $\wp$
of the number field ${\mathbf K}$).
\end{rem}

\begin{rem}
Considering the group $\left(\hbox{Rec}_{p\times p}({\mathbf K})\Big/
\mathcal{FS}_{p\times p}\right)^*$ of invertible elements in the 
quotient algebra
$\hbox{Rec}_{p\times p}({\mathbf K})\Big/
\mathcal{FS}_{p\times p}$ (where $\mathcal{FS}_{p\times p}
\subset \hbox{Rec}_{p\times p}({\mathbf K})$ denotes the 
two-sided ideal of all
finitely supported elements $X\in\hbox{Rec}_{p\times p}({\mathbf K})$ 
with $X[U,W]=0$ except for 
finitely many words $(U,W)\in{\mathcal M}_{p\times p}$) 
we get a group homomorphism
$$\hbox{GL}_{p-rec}({\mathbf K})\longrightarrow
\left(\hbox{Rec}_{p\times p}({\mathbf K})\Big/
\mathcal{FS}_{p\times p}\right)^*
$$
with kernel $\oplus_{n=0}^\infty \hbox{GL}({\mathbf K}^{p^n})$.\end{rem}

\begin{rem} There are three projective versions of 
the group $\hbox{GL}_{p-rec}({\mathbf K})$. One can either consider the
quotient-group $\hbox{GL}_{p-rec}({\mathbf K})/({\mathbf K}^* \hbox{Id})$
or the quotient-group 
$\hbox{GL}_{p-rec}({\mathbf K})/(\hbox{GL}_{p-rec}({\mathbf K})\cap 
{\mathcal C}_{p-rec}({\mathbf K}))$ where ${\mathcal C}_{p-rec}
({\mathbf K})\subset \hbox{Rec}_{p\times
  p}({\mathbf K})$ denotes the center of $\hbox{Rec}_{p\times p}
({\mathbf K})$.
Finally, one can also consider equivalence classes
by invertible central elements in ${\mathbf K}^{{\mathcal M}_{p\times
    p}}$ of all recurrence matrices
$X\in \hbox{Rec}_{p\times p}({\mathbf K})$ for which there
exists a recurrence matrix 
$Y\in \hbox{Rec}_{p\times p}({\mathbf K})$ such that
$XY$ is central and invertible
in ${\mathbf K}^{{\mathcal M}_{p\times p}}$.

The obvious homomorphism
$\hbox{GL}_{p-rec}({\mathbf K})$ into this last group
is in general neither injective 
nor surjective as can be seen as follows:
For $p\geq 2$, we denote by $\hbox{Id}\in \hbox{Rec}_{p\times p}
({\mathbf K})$ the identity recurrence matrix
and by $J\in\hbox{Rec}_{p\times p}
({\mathbf K})$ the recurrence matrix of complexity $1$ defined by
$J[U,W]=1$ for all $(U,W)\in {\mathcal M}_{p\times p}$. 
Choose $\alpha,\beta\in{\mathbb C}$ such that 
$\alpha(\alpha+p^l\beta)\not=0$ for all $l\in{\mathbb N}$
and define $X\in \hbox{Rec}_{p\times p}(\mathbb C)$ by
$$X[{\mathcal M}_{p\times p}^l]=\alpha\ \hbox{Id}[{\mathcal
  M}_{p\times p}^l]+\beta\ J[{\mathcal M}_{p\times p}^l].$$ 
The computation
$$\begin{array}{c}
\displaystyle \left(\alpha\ \hbox{Id}[{\mathcal M}_{p\times p}^l]+\beta\ J[
{\mathcal M}_{p\times p}^l]\right)\left((\alpha+p^l\beta)
\hbox{Id}[{\mathcal M}_{p\times p}^l]-\beta\ J[{\mathcal M}_{p\times
  p}^l]\right)\\ \displaystyle 
=\alpha(\alpha+p^l\beta)\ \hbox{Id}[{\mathcal M}_{p\times p}^l]
\end{array}$$
shows that $X$ is invertible in the quotient of $\hbox{Rec}_{p\times
  p}(\mathbb C)$ by central elements of $\mathbf
K^{\mathcal M_{p\times p}}$.

Choosing $\alpha=\beta=1$, the eigenvalues of $\left(\hbox{Id}+J
\right)[{\mathcal M}_{p\times p}^l]$ are $1+p^l$and $1$ with
multiplicity $p^l-1$. We have thus $\det(\left(\hbox{Id}+J
\right)[{\mathcal M}_{p\times p}^l])=1+p^l$. Choosing an
odd prime $\wp$ such that $p\pmod \wp$ is a non-square in 
${\mathbb Z}/\wp{\mathbb Z}$ and setting $l=(\wp-1)/2$, we have
 $\det(\left(\hbox{Id}+J
\right)[{\mathcal M}_{p\times p}^l])\equiv 0\pmod \wp$.
Since the number of such primes $\wp$ is infinite, Proposition
\ref{critfinprimes} or Remark \ref{reductionmodp}
imply that $\left(\hbox{Id}+J
\right)\not\in \hbox{GL}_{p-rec}({\mathbb Q})$. A similar
argument implies finally that the class of $X$ in the projective 
quotient corresponds to no element in $\hbox{GL}_{p-rec}({\mathbb Q})$.
\end{rem}

\section{Examples of groups of recurrence matrices}
\subsection{$\hbox{GL}_p({\mathbf K})\subset \hbox{GL}_{p-rec}({\mathbf K})$}
Any matrix $g$ of size $p\times p$ with 
coefficients $g_{u,w},0\leq u,w<p$ in ${\mathbf K}$ gives rise to a 
recurrence matrix $\mu(g)=G\in \hbox{Rec}_{p\times p}({\mathbf K})$
of complexity $1$ by considering the $n-$th tensor-power
$g\otimes g\otimes\cdots \otimes g$ and setting
$$G[u_1\dots u_n,w_1\dots w_n]=g_{u_1,w_1}g_{u_2,w_2}\cdots
g_{u_n,w_n}.$$
Since $\mu(g)\mu(g')=\mu(gg')$ (and $\mu(\hbox{Id}(\hbox{GL}_p))=
\hbox{Id}(\hbox{GL}_{p-rec})$), the map $g\longmapsto \mu(g)$ induces
an injective homomorphisme 
$\hbox{GL}_{p}({\mathbf K})\longmapsto \hbox{GL}_{p-rec}({\mathbf K})$.

\subsection{An infinite cyclic group related to the shift}
\label{shiftexample}

The recurrence matrix $A=A_1\in \hbox{Rec}_{2\times 2}({\mathbf K})$ 
presented by 
$(A_1,A_2)[\emptyset,\emptyset]=(1,1)$ and shift-matrices
$$\rho(0,0)=\left(\begin{array}{cc}1&0\\-1&0\end{array}\right),\qquad
\rho(0,1)=\left(\begin{array}{cc}0&0\\1&1\end{array}\right),$$
$$\rho(1,1)=\left(\begin{array}{cc}0&0\\1&0\end{array}\right),\qquad
\rho(1,1)=\left(\begin{array}{cc}1&0\\-1&0\end{array}\right)_ .$$
yields permutation matrices $A_1[{\mathcal M}_{2\times 2}^l]$ associated 
to the cyclic permutation $(0\ 1\ 2\ \dots 2^l-1)$
defined by addition of $1$ modulo $2^l$. The first few matrices
$A_1[{\mathcal M}_{2\times 2}^l], l=0,1,2$ (using the bijection
$s_1\dots s_n\longmapsto \sum_{j=1}^n s_j2^{j-1}$ for rows and
columns) are
$$\left(\begin{array}{c} 1\end{array}\right),
\left(\begin{array}{cc} 0&1\\1&0\end{array}\right),
\left(\begin{array}{cccc} 0&0&0&1\\1&0&0&0\\0&1&0&0\\0&0&1&0
\end{array}\right)\ .$$
All matrices $A_2[{\mathcal M}_{2\times 2}^l]$ have a coefficient 
$1$ in the upper-right corner and zero coefficients everywhere else.

The inverse of $A_1$ is given by $A_1^t$. The two converging 
elements $\rho(0,0)A_1^t,\rho(0,0)A_1\in\hbox{End}_{2-rec}({\mathbb
  Q})$
correspond to the shift $(x_0,x_1,\dots)\longmapsto(x_1,x_2,\dots)$
and its section $(x_0,x_1,\dots)\longmapsto (0,x_0,x_1,\dots)$ 
on converging elements in ${\mathbf K}^{{\mathcal M}_2}$. 
This example has an obvious generalization to 
$\hbox{GL}_{p-rec}({\mathbf K})$ for all $p\geq 2$.

\subsection{Diagonal groups}

The set ${\mathcal D}_{p-rec}({\mathbf K})\cap \hbox{GL}_{p-rec}
({\mathbf K})$ of all diagonal recurrence matrices in 
$\hbox{GL}_{p-rec}({\mathbf K})$ forms a commutative subgroup
containing the maximal central subgroup
${\mathcal C}^*=
{\mathcal C}_{p-rec}({\mathbf K})\cap \hbox{GL}_{p-rec}({\mathbf
  K})$ of $\hbox{GL}_{p-rec}({\mathbf K})$.
For ${\mathbf K}$ of characteristic $\not= p$
the trace map establishes a bijection between ${\mathcal C}^*$
and 
$\hbox{GL}_{1-rec}({\mathbf K})$. An element 
of ${\mathcal C}^*$ (in characteristic
$\not= p$) corresponds to an element $A\in\hbox{Rec}_{1\times 1}({
\mathbf K})$ which is invertible in the algebra 
(or, equivalently, in the function-ring) $\hbox{Rec}_{1\times 1}({
\mathbf K})$. Such an element is thus encoded by a rational function
$\sum_{n=0}^\infty \alpha_n x^n\in {\mathbf K}(x)$ having only
non-zero coefficients such that 
$\sum_{n=0}^\infty \frac{x^n}{\alpha_n}$ is
also rational. Examples of such functions are eg. $\frac{1}{1-\lambda
  z}$ with $\lambda\not= 0$, functions having only
non-zero, ultimately periodic coefficients $\alpha_0,\alpha_1,\dots$,
and, more generally, generating functions of 
recurrence matrices in $\hbox{Rec}_{1\times 1}({\mathbf K})$ with 
values in a finite subset of ${\mathbf K}^*$.

\subsection{Lower (or upper) triangular groups}

Lower (or upper) triangular elements in $\hbox{GL}_{p-rec}(
{\mathbf K})$ form a group ${\mathcal L}^*={\mathcal L}_{p\times
  p}({\mathbf K})\cap \hbox{GL}_{p-rec}({\mathbf K})$ containing 
the subgroup consisting of all converging lower triangular 
elements in $\hbox{GL}_{p-rec}(
{\mathbf K})$. A still smaller subgroup is given by considering
all lower triangular elements of Toeplitz type
in $\hbox{GL}_{p-rec}(
{\mathbf K})$. A slight modification of the
proof of Proposition \ref{groupconstr} shows that the set 
$\left({\mathcal L}^*\right)^p\times \left(
\hbox{Rec}_{p\times p}({\mathbf K})\right)^{p\choose 2}$ injects into 
${\mathcal L}^*$. It would be interesting to know if every conjugacy
class in $\hbox{GL}_{p-rec}({\mathbf K})$ intersects 
${\mathcal L}^*\left({\mathcal L}^*\right)^t$ where 
$\left({\mathcal L}^*\right)^t$ denotes the group of upper triangular 
recurrence matrices obtained by transposing $\mathcal L^*$.

Consider the commutative subgroup ${\mathcal T}^*$ of all
converging lower triangular elements of Toeplitz type 
in $\hbox{GL}_{p-rec}(
{\mathbf K})$. Converging elements of ${\mathcal T}_{p-rec}({\mathbf
    K})$ form a differential subring of ${\mathbf K}[[x]]$, 
(see \S \ref{Toeplitz}), and the logarithmic derivative
$G\longmapsto G'/G$ defines a group homomorphism
from ${\mathcal T}^*$ into an additive subgroup
of ${\mathcal T}_{p\times p}({\mathbf K})$.
In the case ${\mathbf K}\subset {\mathbb
  C}$, elements of ${\mathcal T}^*$ correspond
to some holomorphic functions on the open unit disc of ${\mathbb C}$
which have no zeroes or poles in the open unit disc, cf. Remark 
\ref{remholToepl}. 

Examples of such elements in ${\mathcal T}^*$ are given
by rational functions involving only cyclotomic polynomials. The 
rational function $\frac{1}{1-t}$ for instance corresponds
to the invertible recurrence matrix of Toeplitz-type with limit
the unipotent Toeplitz matrix consisting only of $1'$s below the 
diagonal. 
A more exotic example is given
by the element of ${\mathcal T}^*$ associated to the
power series $\prod_{n=0}^\infty (1+it^{2^n})\in{\mathbb C}[[t]]$
(where $i^2=1$ is a square root of $1$) with inverse series
$(1-t^2)\prod_{n=0}^\infty (1-it^{2^n})$.

If ${\mathbf K}$ is contained in the algebraic closure of 
a finite field, the group ${\mathcal T}^*$ contains the 
subring of all rational functions without zero or pole at the origin. 


\subsection{Orthogonal groups}\label{subsorthog}
A recurrence matrix $A\in \hbox{Rec}_{p\times p}({\mathbf K})$ 
is {\it symmetric} if 
$A[U,W]=A[W,U]$ for all $(U,W)\in{\mathcal M}_{p\times p}$.
A complex recurrence matrix $A\in\hbox{Rec}_{p\times p}({\mathbf K})$ 
is {\it hermitian} if $A[U,W]=\overline{A[W,U]}$, for all
$(U,W) \in {\mathcal M}_{p\times p}$, with $\overline x$ denoting the 
complex conjugate of $x\in {\mathbb C}$.
(More generally, one can define ``hermitian'' matrices over a field
admitting an involutive automorphism.)
A real symmetric recurrence matrix is 
{\it positive definite} if all matrices
$A[{\mathcal M}_{p\times p}^l]$ are positive definite.

A symmetric or hermitian matrix $A\in\hbox{Rec}_{p\times p}({\mathbf K})$
is {\it non-singular} if $\det(A[{\mathcal M}_{p\times p}^l])
\in{\mathbf K}^*$ for all $l\in 
{\mathbb N}$ and {\it non-degenerate} if
$A\in\hbox{GL}_{p-rec}({\mathbf K})$. 
An example of a non-singular positive definite symmetric
recurrence matrix 
is the diagonal recurrence matrix $A\in \hbox{Rec}_{p\times p}
({\mathbb Q})$ with diagonal coefficients $A[U,U]=l+1$
for $(U,U)\in{\mathcal M}_{p\times p}^l$ of length $l$.
An example of a positive definite non-degenerate symmetric
recurrence matrix is the identity matrix $\hbox{Id}\in
\hbox{Rec}_{p\times p}({\mathbb Q})$.

Such a non-singular matrix $A \in\hbox{Rec}_{p\times p}({\mathbf K})$
defines the {\it orthogonal group}
$$\hbox{O}(A)=\{B\in\hbox{GL}_{p-rec}({\mathbf K})\ \vert \ B^tAB=A\}
\subset \hbox{GL}_{p-rec}({\mathbf K})$$
of $A$ in the symmetric case and the {\it unitary group}
$$\hbox{U}(A)=\{B\in\hbox{GL}_{p-rec}({\mathbf K})\ \vert \ 
\overline{B}^tAB=A\}$$
of $A$ in the case of a hermitian recurrence matrix $A$. 
For ${\mathbf K}$ a real field,
one speaks also of {\it Lorenzian groups} if
the symmetric matrix $A$ is not positive definite.

The following obvious proposition relating presentations of $A,A^t\in
\hbox{Rec}_{p\times p}({\mathbf K})$
and $\overline{A}$ (over ${\mathbf K}={\mathbb C}$)
is mainly a restatment of Proposition \ref{propprestransposed}, 
recalled for the convenience of the reader. Its easy proof is omitted.

\begin{prop} \label{Atransposed}
The following assertions are equivalent:

\ \ (i) $(A_1,\dots,A_{a})[\emptyset,\emptyset]=(\alpha_1,\dots,
\alpha_a)\in{\mathbf K}^a,\rho(s,t)\in{\mathbf K}^{a\times a},0\leq s,t<r$ 
is a presentation of $A=A_1$.

\ \ (ii) $(A_1,\dots,A_{a})[\emptyset,\emptyset]=(\alpha_1,\dots,
\alpha_a),\tilde \rho(s,t)\in{\mathbf K}^{
a\times a},0\leq s,t<p$ with
$\tilde \rho(s,t)=\rho(t,s)$ is a presentation of $A^t=A_1$.

\ \ (iii) (Over ${\mathbf K}={\mathbb C})$) 
${(A_1,\dots,A_{a})}[\emptyset,\emptyset]=(\overline{\alpha_1},
\dots,\overline{\alpha_a}),\overline{\rho}(s,t)\in{\mathbf K}^{
a\times a},0\leq s,t<r$
is a presentation of $\overline{A}=A_1$ with $\overline{X}$ denoting 
complex conjugation applied to all coefficients of $X$.
\end{prop}

\subsection{Symplectic groups}
A recurrence matrix $A\in \hbox{Rec}_{p\times p}({\mathbf K})$ 
is {\it antisymmetric} if $A[U,W]=-A[W,U]$ 
for all $(U,W)\in{\mathcal M}_{p\times p}$.
An antisymmetric recurrence matrix $A\in\hbox{Rec}_{p\times p}({\mathbf
  K})$ is {\it symplectic} if $\det(A[{\mathcal M}_{p\times p}^l])\in
{\mathbf K}^*$ for all $l\geq 1$. 
If ${\mathbf K}$ is of characteristic $2$, 
we require moreover $A[U,U]=0$ for all $(U,U)\in
{\mathcal M}_{p\times p}$. Symplectic recurrence matrices exist 
only for $p$ even.

 A symplectic recurrence matrix $A$ is 
{\it non-degenerate} if $\tilde A\in\hbox{GL}_{p-rec}({\mathbf K})$
where $\tilde A[\emptyset,\emptyset]=1$ and
$\tilde A[U,W]=A[U,W]$ if $(U,W)\not=[\emptyset,\emptyset]$.

For $A\in\hbox{Rec}_{p\times p}({\mathbf K})$ a symplectic recurrence
matrix, the associated {\it symplectic group} $Sp(A)\subset 
GL_{p-rec}({\mathbf K})$ of recurrence matrices is defined as
$$\hbox{Sp}(A)=\{B\in\hbox{GL}_{p-rec}({\mathbf K})\ \vert \
B^tAB=A\}$$
where the value $B[\emptyset,\emptyset]$ can be 
neglected.

\subsection{Groups generated by elements of bounded complexity}
\label{boundedsection}

Denote by $\Gamma_{a,b}\subset \hbox{GL}_{p-rec}({\mathbf K})$
the group generated by all elements $A\in \hbox{GL}_{p-rec}({\mathbf K})$
of complexity $\leq a$ with inverse $B=A^{-1}$ of complexity $\leq b$.
We have $\Gamma_{a,b}=\Gamma_{b,a}$, $\Gamma_{a,b}\subset
\Gamma_{a',b'}$ if $a\leq a',b\leq b'$ and Proposition 
\ref{propnobound} implies that many of these inclusions are strict.

Moreover, the set of generators of $\Gamma_{a,b}$ (elements
in  $\hbox{GL}_{p-rec}({\mathbf K})$ of complexity $\leq a$ with inverse
of complexity $\leq b$) is a union of algebraic sets 
since they can be described by 
(a finite union of) polynomial equations. 
 
\begin{rem}
The group $\Gamma_{1,1}=\Gamma_{1,2}=\Gamma_{1,3}=\dots$
is isomorphic to ${\mathbf K}^*\times \hbox{GL}({\mathbf K}^p)$.

For ${\mathbf K}$ a finite field,
the group $\Gamma_{a,b}$ is finitely generated and the sequence
$$\Gamma_{a,1}\subset\Gamma_{a,2}\subset\Gamma_{a,3}\subset \dots$$
stabilizes. It would be interesting to determine the smallest
integer
$A=A(a,\mathbf K)$ such that $\Gamma_{a,A}=\Gamma_{a,b}$ for all
$b\geq A$. The first non-trivial case is the determination of
$A(2,\mathbb F_2)$.
\end{rem}

\begin{rem} One can similarly consider the subalgebra
${\mathcal R}_a\subset\hbox{Rec}_{p\times p} ({\mathbf K})$
generated as an algebra by all recurrence matrices 
of complexity $\leq a$
in $\hbox{Rec}_{p\times p} ({\mathbf K})$. The subalgebra
${\mathcal R}_a$ of $\hbox{Rec}_{p\times p} ({\mathbf K})$
is always recursively closed. 

The following example shows that many inclusions 
${\mathcal R}_a\subset{\mathcal R}_{a+1}$
are strict for $\hbox{Rec}_{p\times p}({\mathbb Q})$: 
Consider a central diagonal recurrence matrix 
with diagonal coefficients $A[U,U]=\alpha_l,(U,U)\in{\mathcal
  M}_{p\times p}^l$ for $\alpha_0,\alpha_1,\alpha_2,\dots\subset\mathbb
Q$ a periodic sequence of minimal period-length a prime number $\wp$. 
Such an element has complexity $(\wp -1)$ and cannot be contained 
in the sub-algebra generated by recurrence matrices of lower complexities
in  $\hbox{Rec}_{p\times p}({\mathbb Q})$.
\end{rem}


\section{Computing $G^{-1}$ for $G\in\hbox{GL}_{p-rec}({\mathbf K})$}

The aim of this section is to discuss a few difficulties when 
computing inverses of recurrence matrices in
$\hbox{GL}_{p-rec}({\mathbf K})$.

\begin{defin} The depth of $A\in{\mathbf K}^{{\mathcal M}_{p\times
        q}}$
is the smallest element $D\in{\mathbb N}\cup\{\infty\}$ such that 
$\rho({\mathcal M}_{p\times q}^{\leq D})A$ spans $\overline{A}^{rec}$
(where ${\mathcal M}_{p\times q}^\infty={\mathcal M}$).
\end{defin}

It is easy to check that $A$ is a recurrence matrix if and only if
$A$ has finite depth $D<\infty$.

The following result is closely related to Proposition \ref{nicecoord}.

\begin{prop} \label{depthdefin}
We have 
$$\hbox{dim}\left(\sum_{(U,W)\in{\mathcal M}_{p\times q}^{\leq k}}
{\mathbf K}\rho(U,W)A\right)<
\hbox{dim}\left(\sum_{(U,W)\in{\mathcal M}_{p\times q}^{\leq k+1}}
{\mathbf K}\rho(U,W)A\right)$$
if $k$ is smaller than the depth $D$ of $A$.
\end{prop}

\begin{cor} The depth of any non-zero recurrence matrix
$A$ is smaller than its complexity $\hbox{dim}(\overline{A}^{rec})$.
\end{cor}

{\bf Proof of Proposition \ref{depthdefin}} (See also the proof
of Proposition \ref{nicecoord}.) The equality
$$\sum_{(U,W)\in{\mathcal M}_{p\times q}^{\leq k}}
{\mathbf K}\rho(U,W)A=
\sum_{(U,W)\in{\mathcal M}_{p\times q}^{\leq k+1}}
{\mathbf K}\rho(U,W)A$$
implies that these vector-spaces are recursively closed and
coincide thus with $\overline{A}^{rec}$. 
\hfill $\Box$

Given a presentation ${\mathcal P}$
of $G\in\hbox{GL}_{p-rec}({\mathbf K})$,
there are two obvious methods for computing a presentation 
of its inverse $G^{-1}$. The first method analyzes the matrices
$$\left(G[{\mathcal M}_{p\times p}^0]\right)^{-1},
 \left(G[{\mathcal M}_{p\times p}^1]\right)^{-1},\dots,
\left(G[{\mathcal M}_{p\times p}^a]\right)^{-1}$$
with $a$ huge enough in order to guess a presentation 
$\tilde {\mathcal P}$ 
of $G^{-1}$. This can be done if $a\geq
1+D+N$ where $D$ is the depth and $N$ the saturation level of 
$G^{-1}$. It is then straightforward to check if the presentation
$\tilde P$ is correct by computing the matrix-product of $G$ 
and the recurrence matrix presented by $\tilde {\mathcal P}$.
The limitation of this method is the need of 
inverting the square-matrix $G[{\mathcal M}_{p\times p}^a]$
of large order $p^a\times p^a$. Below, we will describe an algorithm
based on this method.

The second method is to guess an upper bound $b$ for the complexity
of $G^{-1}$, to write down a generic presentation $\tilde P$
of complexity
$b$ where the initial data and shift matrices involve a set
of $d+p^2d^2$ unknowns. Equating the matrix-product of $G$ with the
recurrence matrix presented by $\tilde {\mathcal P}$ to the 
identity yields polynomial equations. 
We will omit a detailed discussion of this method since it
seems to be even worse than the first one.

\begin{rem} Another important issue which we will not adress
here is the existence of a finite algorithm which is able to tell if an
element $A\in\hbox{Rec}_{p\times p}({\mathbf K})$ (given,
say, by a minimal presentation) is or is not invertible in 
$\hbox{GL}_{p\times p}({\mathbf K})$? The algorithm presented below
will always succeed (with finite time and memory requirements)
in computing an inverse of $A$ for $A\in\hbox{GL}_{p-rec}({\mathbf
  K})$. It will however fail to stop (or more likely, use up all
memory on your computing device) if $A\not\in 
\hbox{GL}_{p-rec}({\mathbf K})$ 
is invertible in ${\mathbf K}^{{\mathcal M}_{p\times p}}$.

Proposition \ref{propnobound} is perhaps an obstruction to the
existence of such an algorithm.
\end{rem}

\subsection{An algorithm for computing
  $G^{-1}\in\hbox{GL}_{p-rec}({\mathbf K})$}

Given a presentation ${\mathcal P}$ of $G\in\hbox{GL}_{p-rec}({\mathbf
  K})$, the following algorithm computes a presentation 
$\tilde {\mathcal P}$ of $G^{-1}$ in a finite number of steps.

{\bf Step 1} Set $D=0$.

{\bf Step 2} Compute the saturation level $N$ of the 
(not necessarily recursively closed) vector space
spanned by $\rho({\mathcal M}_{p\times p}^{\leq D})G^{-1}$. 
(This needs inversion of the $p^{D+N+1}\times p^{D+N+1}$ matrix
$G[{\mathcal M}_{p\times p}^{\leq D+N+1}]$ where
$N<\frac{p^{D+1}-1}{p-1}$ 
is the saturation level of the vector space spanned by
$\rho({\mathcal M}_{p\times p}^{\leq D})G^{-1}$.)

{\bf Step 3} Supposing the depth $D$ correct, use the 
saturation level $N$ of step 2 for computing a presentation 
$\tilde P_D$ using the finite-dimensional vector space 
spanned by $\left(\rho({\mathcal M}_{p\times p}^{\leq D})G^{-1}\right)
[{\mathcal M}_{p\times p}^{\leq N+1}]$. 

{\bf Step 4} If the recurrence matrix $\tilde G$ defined by the
presentation $\tilde {\mathcal P}_D$ satisfies $G\tilde G=\hbox{Id}$, 
stop and print the presentation $\tilde {\mathcal
  P}$. Otherwise, increment $D$ by $1$ and return to step 2.

\begin{rem} The expensive part (from a computational view) of
the algorithm are steps 2 and 3 and involve computations with
large matrices (for $p\geq 2$). A slight improvement is to merge
step 2 and 3 and to do the computations for guessing the presentation
$\tilde {\mathcal P}_D$ of step 3 at once during step 2.

One could avoid a lot of iterations by running step 2
simultaneously for $D$ and $D+1$. The cost of this ``improvement'' 
is however an extra factor of $p$ in the size of the 
involved matrices and should thus be avoided since
step 4 is faster than step 2.
\end{rem} 


\section{Lie algebras}

The Lie bracket $[A,B]=AB-BA$ turns the algebra $\hbox{Rec}_{p\times
  p}
({\mathbf K})$ (or ${\mathbf K}^{{\mathcal M}_{p\times p}}$) into a
Lie algebra. For ${\mathbf K}$ a suitable
complete topological field (say ${\mathbf K}={\mathbb R}$
or ${\mathbf K}={\mathbb C}$), the exponential function 
$X\longmapsto \hbox{exp}(X)=\sum_{k=0}^\infty \frac{X^k}{k!}$ is well-defined and 
continuous
for both topologies defined in \S \ref{topologies}.
The exponential function does however
not preserve the subspace $\hbox{Rec}_{p\times p}({\mathbf
  K})$ of recurrence matrices and the associated Lie group 
$$\{\hbox{exp}(A)\in {\mathbf K}^{{\mathcal M}_{p\times p}}\ \vert
\ A\in \hbox{Rec}_{p\times
  p}({\mathbf K})\}$$
is thus only a group in ${\mathbf K}^{{\mathcal M}_{p\times p}}$.

The Lie algebra $\hbox{Rec}_{p\times p}({\mathbb C})$ contains
analogues of the classical Lie algebras of type $A,B,C$ and $D$.

The analogue of type $A$ is given by the vector space
of recurrence matrices $X\in\hbox{Rec}_{p\times p}({\mathbb C})$ such that 
$\hbox{tr}(X)=0\in\hbox{Rec}_{1\times 1}$.

The $B$ and $D$ series are defined as the set of all recursive 
antisymmetric matrices
of $\hbox{Rec}_{p\times p}({\mathbb C}),p\geq 3$ odd, for the $B$
series
and of $\hbox{Rec}_{p\times p}({\mathbb C}),p\geq 4$ even, 
for the $D$ series.

The $C$ series is defined only for even $p$ and consists of
all recursive matrices $X\in \hbox{Rec}_{p\times p}({\mathbb C})$
such that $\Omega X=(\Omega X)^t$ where
$\Omega\in \hbox{Rec}_{p\times p}({\mathbb C})$ has complexity $2$
and $\Omega[{\mathcal M}_{p\times p}^l]$ is of the form
$\left(\begin{array}{cc} 0&\hbox{Id}\\-\hbox{Id}&0\end{array}\right)$
for $l\geq 1$ with $\hbox{Id}$ and $0$ denoting the identity matrix and
the zero matrix of size $p^l/2\times p^l/2$ (the value
of $\Omega[\emptyset,\emptyset]$ is irrelevant).

A different and perhaps more natural way to define
the $C$ series is to consider triplets $A,B,C\in\hbox{Rec}_{p\times
  p}({\mathbf K})$ (for all $p\in{\mathbb N}$) 
of recurrence matrices with $B=B^t$ and $C=C^t$ symmetric.
The $C$ series is then the Lie subalgebra in 
$\prod_{l=0}^\infty {\mathbf K}^{2p^l\times 2p^l}$ of elements
given by
$$\left(\begin{array}{cc} A[{\mathcal M}_{p\times p}^l]&
B[{\mathcal M}_{p\times p}^l]\\
C[{\mathcal M}_{p\times p}^l]&
-A^t[{\mathcal M}_{p\times p}^l]\end{array}\right).$$

I ignore if there are natural ``recursive'' analogues of (some of) 
the exceptional simple Lie algebras.

\begin{rem} Analogues of type $B,C,D$ Lie algebras can also
be defined using arbitrary non-singular
symmetric or symplectic recurrence matrices in 
$\hbox{Rec}_{p\times p}({\mathbf K})$. 
\end{rem}

\begin{rem} A Lie-algebra ${\mathcal L}\subset \hbox{Rec}_{p\times p}
({\mathbf K})$ is in general not a recursively closed subspace
of $\hbox{Rec}_{p\times p}
({\mathbf K})$.
\end{rem}

\begin{rem}
It would be interesting to understand the 
algebraic structure of Lie algebras in
$\hbox{Rec}_{p\times p}({\mathbf K})$.
Given such a Lie algebra ${\mathcal L}
\subset \hbox{Rec}_{p\times p}({\mathbf
  K})$, the intersection ${\mathcal L}\cap \mathcal{FS}_{p\times p}$ 
with the vector-space $\mathcal{FS}_{p\times p}$ of all
finitely supported elements defines an ideal in ${\mathcal L}$.
The interesting object is
thus probably the quotient Lie algebra ${\mathcal L}/({\mathcal L}\cap 
\mathcal{FS}_{p\times p})$.
What is the structure of this quotient algebra for the
analogues in $\hbox{Rec}_{p\times p}({\mathbb C})$ of the $A-D$ series?
\end{rem} 

\subsection{Lie algebras in the convolution ring ${\mathbf
    K}^{{\mathcal M}_p}$}

Using the convolution-ring structure on 
${\mathbf K}^{{\mathcal M}_{p\times q}}$ we get other, different
Lie algebras. Since the convolution-structure depends only
on the product $pq$, we restrict ourself to 
${\mathbf K}^{{\mathcal M}_p}$ for simplicity. 

The
three obvious subalgebras in ${\mathbf K}^{{\mathcal M}_p}$
(corresponding to formal power series in $p$ non-commuting 
variables)
given by ${\mathbf K}^{{\mathcal M}_p}$, $\hbox{Rec}_p({\mathbf K})$
and the vector space $\mathcal{FS}_p\subset
\hbox{Rec}_{p-rec}({\mathbf K})$ of finitely supported elements
give rise to three Lie-algebras with Lie bracket
$[A,B]=A*B-B*A$ (where $*$ stands for the convolution product in
${\mathbf K}^{{\mathcal M}_p}$).

The Lie algebra resulting from $\mathcal{FS}_p$ 
(where $\mathcal{FS}_p$ as a convolution algebra is isomorphic to 
the polynomial algebra in $p$ non-commuting variables) contains 
the free Lie algebra on $p$ generators.

All these Lie algebras are filtrated:
We have $[A,B][\emptyset,\emptyset]=0$ for all
$A,B\in{\mathbf K}^{{\mathcal M}_{p}}$ and
$[A,B][{\mathcal M}_p^{< \alpha+\beta}]=0$ 
if $A[{\mathcal M}_p^{< \alpha}]=0$ and 
$B[{\mathcal M}_p^{< \beta}]=0$.


\section{Integrality and lattices of recurrence vectors}

A recurrence matrix $A\in \hbox{Rec}_{p\times q}({\mathbb Q})$
is {\it integral} if all its coefficients $A[U,W]$ are integers.
This notion can easily be generalized by considering coefficients
in the integral ring ${\mathcal O}_{\mathbf K}$ of algebraic 
integers over a number field ${\mathbf K}$.

We denote by $\hbox{Rec}_{p\times q}({\mathbb Z})$ the ${\mathbb
  Z}-$module of all integral recurrence matrices in 
$\hbox{Rec}_{p\times q}({\mathbb Q})$. Since products
of integral matrices are integral, one can consider the subcategory
${\mathbb Z}^{\mathcal M}\subset{\mathbb Q}^{\mathcal M} $ 
having only integral recurrence matrices as arrows and the subalgebra 
$\hbox{Rec}_{p\times p}({\mathbb Z})\subset
\hbox{Rec}_{p\times q}({\mathbb Q})$ consisting of all integral
recurrence matrices.

Call a presentation $(A_1,\dots,A_d)[\emptyset,\emptyset]$ with 
shift matrices $\rho(s,t)\in\mathbb Q^{d\times d}$ of a 
finite-dimensional recursively closed subspace in 
$\hbox{Rec}_{p\times q}(\mathbb
Q)$ {\it integral}
if $(A_1,\dots,A_d)[\emptyset,\emptyset]\in\mathbb Z^d$ and 
$\rho(s,t)\in\mathbb Z^{d\times d}$ for $0\leq s<p,0\leq t<q$.

\begin{prop} Every integral recurrence matrix $A\in 
\hbox{Rec}_{p\times q}({\mathbb Z})$ admits an 
an integral minimal presentation of $\overline{A}^{rec}$
such that $A_1,\dots,A_a$ is a $\mathbb Z-$basis of 
$\hbox{Rec}_{p\times q}({\mathbb Z})\cap \overline{A}^{rec}$.
In particular, $A$ can be written as
$A=\sum_{j=1}^a \lambda_j A_j$ with
$\lambda_1,\dots,\lambda_a\in\mathbb Z$.
\end{prop}

{\bf Proof} Proposition \ref{propfond} shows that the set 
$(\hbox{Rec}_{p\times q}({\mathbb Z})\cap\overline{A}^{rec})
\subset \overline{A}^{rec}\subset \hbox{Rec}_{p\times q}({\mathbb Q})$
is a free ${\mathbb Z}-$mo\-dule. Since it contains
the recursive set-closure $\rho({\mathcal M}_{p\times q})A$
spanning $\overline{A}^{rec}$, it is of maximal rank $a=\hbox{dim}
(\overline{A}^{rec})$. A ${\mathbb Z}-$basis $A_1,\dots,A_a$
of $(\hbox{Rec}_{p\times q}({\mathbb Z})\cap
\overline{A}^{rec})$ has the required properties.
\hfill$\Box$

We call an integral recurrence matrix 
$A\in\hbox{Rec}_{p\times p}({\mathbb Z})$ {\it unimodular}
if $\det(A[{\mathcal M}_{p\times p}^l])\in\{\pm 1\}$ for all
$l\in{\mathbb Z}$. The set of of all integral unimodular
recurrence matrices in $\hbox{GL}_{p-rec}({\mathbb Q})$
is the {\it unimodular subgroup of integral recurrence matrices}
in $\hbox{Rec}_{p\times p}({\mathbb Q})$.

\begin{rem} The ${\mathbb Q}-$vector space $\hbox{Rec}_{p\times
    q}({\mathbb Z})\otimes_{\mathbb Z}{\mathbb Q}$ is in general a
strict subspace of $\hbox{Rec}_{p\times
    q}({\mathbb Q})$: An element $A\in \hbox{Rec}_{p\times
    q}({\mathbb Q})\setminus\hbox{Rec}_{p\times
    q}({\mathbb Z})\otimes_{\mathbb Z}{\mathbb Q}$ is for instance
  given by $A[{\mathcal M}_{p\times q}^l]=\frac{1}{2^l},\ l
\in{\mathbb N}$. For $q=p$, the vector space $\hbox{Rec}_{p\times
    p}({\mathbb Z})\otimes_{\mathbb Z}{\mathbb Q}$ is
thus a proper subalgebra of the algebra $\hbox{Rec}_{p\times
    p}({\mathbb Q})$. However, given an arbitrary 
recurrence matrix $A\in\hbox{Rec}_{p\times
  q}({\mathbb Q})$, there exists by Proposition \ref{critfinprimes}
an integer $\alpha\geq 1$ and an integer $\lambda\geq 1$
with $\alpha H_\lambda A\in\hbox{Rec}_{p\times
  q}({\mathbb Z})$ where $H_\lambda\in{\mathcal C}_{p-rec}
({\mathbb Q})$ is the integral diagonal matrix in the center 
${\mathcal C}_{p-rec}({\mathbb Q})$ of 
$\hbox{Rec}_{p\times p}({\mathbb Q})$ with diagonal coefficients
$H_\lambda[U,U]=\lambda^l$ for $(U,U)\in{\mathcal M}_{p\times p}^l$.
\end{rem} 

\begin{rem}\label{remintnoninvert} 
An integral unimodular recurrence matrix 
in $\hbox{Rec}_{p\times p}({\mathbb Z})$ is
generally not invertible in $\hbox{Rec}_{p\times p}({\mathbb Z})$. 
An example is given by the converging
lower triangular recurrence matrix $A\in{\mathcal T}_2({\mathbb Z})
\subset \hbox{Rec}_{2\times 2}({\mathbb Z})$ of Toeplitz type
(cf. Section \ref{substoeplalg})
with generating series $1-2z\in {\mathbb Z}[[z]]$
already mentioned in the second part of Remark \ref{remnoinverse}.
Its inverse $A^{-1}\in{\mathcal T}_2({\mathbb Z})\subset
{\mathbb Q}^{{\mathcal M}_{2\times 2}}$
corresponds to the generating series $\frac{1}{1-2z}=1+2z+4z^2+8z^3+\dots
\in{\mathbb Z}[[z]]$ and cannot be a recurrence matrix 
by Proposition \ref{majorationcriterion}.
\end{rem}

\subsection{Finite-dimensional lattices}

\begin{defin} A {\it finite-dimensional lattice} of
  $\hbox{Rec}_{p\times q}({\mathbb C})$ is a free ${\mathbb Z}-$module
of finite rank in $\hbox{Rec}_{p\times q}({\mathbb C})$.
\end{defin}

A finite-dimensional lattice of
  $\hbox{Rec}_{p\times q}({\mathbb C})$ is discrete for both 
topologies introduced in \S \ref{topologies}.

Since the multiplicative structure is irrelevant for lattices,
it is enough to consider lattices in $\hbox{Rec}_{p}
({\mathbb C})$. A particularly beautiful set
of lattices is given by lattices which are recursively
closed sets (ie. they satisfy the inclusion $\rho({\mathcal M}_p)
\Lambda\subset \Lambda$). An example of such a lattice
is the subset $\hbox{Rec}_p({\mathbb Z})\cap \overline{A}^{rec}\subset
\overline{A}^{rec}$ for $A\in\hbox{Rec}_p({\mathbb Q})$.

Given two recurrence vectors $A,B\in\hbox{Rec}_p({\mathbb R})$, 
we define their 
{\it scalar-product} as the element 
$$\langle A,B\rangle=A^tB\in\hbox{Rec}_{1}({\mathbb R}).$$
More precisely, such a scalar-product can be represented by the
generating series
$$\langle A,B\rangle_z=
\sum_{l=0}^\infty z^l\sum_{U\in{\mathcal M}_p^l}A[U]B[U]
\in{\mathbb R}[[z]]$$
and defines a rational function. For $\Lambda\subset
\hbox{Rec}({\mathbb R})$ a lattice, we call the rational function
$z\longmapsto \det(\langle A_i,A_j\rangle_z)$ (with $A_1,A_2,\dots$ a
${\mathbb Z}-$basis of $\Lambda$) the {\it determinant} of $\Lambda$.
Given a finite dimensional subspace ${\mathcal A}\subset 
\hbox{Rec}_{p}({\mathbb R})$ (spanned eg. by a lattice) 
there exists an open interval $(0,\alpha({\mathcal A}))\subset{\mathbb R}$
such that bilinear map ${\mathcal A}^2\ni(A,B)\longmapsto
\langle A,B\rangle_{z_0}\in{\mathbb R}$ defines
an Euclidean scalar-product on ${\mathcal A}$ for all
$z_0\in(0,\alpha({\mathcal A}))$. 

In particular, such an evaluation yields an isometry between
a lattice $\Lambda\subset \hbox{Rec}_p({\mathbb R})$ 
and a lattice of the ordinary Euclidean vector space.

\begin{rem} It is of course also possible to define 
scalar-products using a symmetric positive definite recurrence matrix
of $\hbox{Rec}_{p\times p}({\mathbb R})$ 
(the case considered above corresponds to the identity). 

One defines similarly Hermitian products and Hermitian lattices
over ${\mathbb C}$.
\end{rem}


\section{Monoids and their linear representations}

Any quotient monoid 
${\mathcal Q}=\langle {\mathcal M}_{p\times q}^1:
{\mathcal R}\subset {\mathcal M}_{p\times q}
\times {\mathcal M}_{p\times q}\rangle$
of the monoid ${\mathcal M}_{p\times q}$ (see chapter 
\ref{xectmonoids} for definitions) defines a subspace
$V_{\mathcal Q}\subset{\mathbf K}^{{\mathcal M}_{p\times q}}$
by setting
$$V_{\mathcal Q}=\{A\in {\mathbf K}^{{\mathcal M}_{p\times q}}\ \vert\ 
\rho(UL_iV)A=\rho(UR_iV)A,U,V\in {\mathcal M}_{p\times
  q},(L_i,R_i)\in{\mathcal R}\}.$$
The space $V_{\mathcal Q}$ is by construction recursively closed, and  
the shift monoid $\rho_{{\mathcal Q}}({\mathcal
  M}_{p\times q})$ acts on $V_{\mathcal Q}$
with a ``kernel'' generated by the relations ${\mathcal R}$.

In the case $p=q$, the subspace $V_{\mathcal Q}$ is in general not
multiplicatively closed. In the next section we will 
however describe a particular case where this happens.

Let us recall here a few elementary facts already discussed in
chapter \ref{sectcatrec}:

A linear representation of a monoid is a morphisme
$\pi:{\mathcal Q}\longrightarrow\hbox{End}(V)$
of some abstract monoid ${\mathcal Q}$ into a submonoid of 
$\hbox{End}(V)$ where $\hbox{End}(V)$ denotes the monoid of all linear
endomorphisms of a vector space $V$. As in the case
of linear representations of groups, one can define 
indecomposable representations 
(without proper non-trivial invariant
subspace), direct sums of representations, irreducible representations
(not a non-trivial direct sum of representations) etc and a
linear representation of a monoid ${\mathcal Q}$
gives rise to a linear representation of its monoid-algebra
${\mathbf K}[{\mathcal Q}]$. 

A linear representation of
the free monoid on $r$ generators is simply a set
$\pi(g_1),\dots,\pi(g_r)\subset \hbox{End}(V)$ of $r$ endomorphisms
corresponding to the free generators $g_1,\dots,g_r$.
If a quotient monoid ${\mathcal Q}$ of a free monoid
has relations ${\mathcal R}$, then
$\pi(L_i)=\pi(R_i)$ for every relation $(L_i,R_i)\in{\mathcal R}$.
Conjugate linear representations of a monoid are in generally 
considered as equivalent. Obviously, every recursively closed 
subspace of ${\mathbf K}^{{\mathcal M}_{p\times q}}$ gives rise
to a representation of the free monoid ${\mathcal M}_{p\times q}$ on
$pq$ generators. Reciprocally, every finite-dimensional 
linear representation $\rho:{\mathcal M}_{p\times q}\longrightarrow
\hbox{End}({\mathbf K}^d)$ defines a
recursively closed finite-dimensional subspace of 
$\hbox{Rec}_{p\times q}({\mathbf K})$ by considering the direct sum (of
dimension $\leq d^2$) of all recursively closed subspaces
with shift-matrices $\rho({\mathcal M}_{p\times q}^1)$ and arbitrary
initial values. The precise description of such subspaces 
will be given in section \ref{birecursivity}.

\begin{rem} A minimal
presentation $A_1,A_2,\dots,A_a$
of an element $A=A_1\in\hbox{Rec}_{p\times q}$ having complexity
$\hbox{dim}({\overline A}^{\hbox{rec}})=a$ yields a linear representation
$\rho({\mathcal M}_{p\times q})$, up to conjugation by elements of
$\hbox{GL}_a({\mathbf K})$ fixing the first basis vector $A_1$.
\end{rem}


\section{Abelian monoids}

For $p,q\in{\mathbb N}$ we consider the {\it free abelian} quotient
monoid ${\mathcal Ab}_{p\times q}\sim {\mathbb N}^{pq}$ 
generated by all $pq$ elements
of ${\mathcal M}_{p\times q}^1$ with relations ${\mathcal R}=\{(XY,YX)
\ \vert X,Y\in {\mathcal M}_{p\times q}^1\}$. The corresponding 
subspace $V_{{\mathcal Ab}_{p\times q}}$ consists of all functions 
$A\in {\mathbf K}^{{\mathcal M}_{p\times q}}$ such that
$$A[u_1\dots u_n,w_1\dots w_n]=A[u_{\pi(1)}\dots u_{\pi(n)},
w_{\pi(1)}\dots w_{\pi(n)}]$$
for all $(u_1\dots u_n,w_1\dots w_n)\in{\mathcal M}_{p\times q}$
and all permutations $\pi$ of $\{1,\dots ,n\}$. The vector space
$V_{{\mathcal Ab}_{p\times q}}$ can thus be identified with the vector 
space ${\mathbf K}[[Z_{0,0},\dots,Z_{p-1,q-1}]]$   
of formal power series in $pq$ commuting variables
$Z_{u,w},0\leq u<p,0\leq w<q$. We denote the vector space 
$V_{{\mathcal Ab}_{p\times q}}$
by ${\mathbf K}^{{\mathcal Ab}_{p\times q}}$.

For $A\in
{\mathbf K}^{{\mathcal Ab}_{p\times q}}$ we have
$$\rho(s,t)\rho(s',t')A=\rho(s',t')\rho(s,t)A$$
for all $0\leq s,s'< p,0\leq t,t'<q$.
The easy computation 
$$\begin{array}{l}
\displaystyle \rho(s,t)\rho(s',t')(AB)=\rho(s,t)\left(\sum_{v'}
(\rho(s',v')A)(\rho(v',t')B)\right)\\
\displaystyle \quad
=\sum_{v,v'}(\rho(s,v)\rho(s',v')A)(\rho(v,t)\rho(v',t')B)\\
\displaystyle \quad
=\sum_{v,v'}(\rho(s',v')\rho(s,v)A)(\rho(v',t')\rho(v,t)B)\\
\displaystyle \quad =\dots=\rho(s',t')\rho(s,t)(AB)\end{array}$$
shows that $AB\in {\mathbf K}^{{\mathcal Ab}_{p\times q}}$
if $A\in {\mathbf K}^{{\mathcal Ab}_{p\times r}}$, 
$B\in {\mathbf K}^{{\mathcal Ab}_{r\times q}}$. We get thus subcategories
${\mathbf K}^{{\mathcal Ab}}$ of ${\mathbf K}^{{\mathcal M}}$ 
and
$\hbox{Rec}({\mathbf K})\cap{\mathbf K}^{{\mathcal Ab}}$ of 
$\hbox{Rec}({\mathbf K})$.

These subcategories contain all elements of complexity $1$. 
The algebra formed by all recurrence matrices in 
${\mathbf K}^{{\mathcal Ab}_{p\times p}}$
is thus not commutative if $p>1$. 

It would be interesting to have other examples of ``natural''
quotient monoids of ${\mathcal M}_{p\times q}, p,q\in{\mathbf N}$, 
giving rise to subcategories in ${\mathbf K}^{\mathcal M}$ and
$\hbox{Rec}({\mathbf K})$.

\begin{rem} The association 
$$A\longmapsto f_A=A[\emptyset]+\sum_{n=1}^\infty \sum_{0\leq
  u_1,\dots,u_n<p} A[u_1\dots u_n]Z_{u_1}\cdots Z_{u_n}$$
of a recurrence vector $A\in{\mathbb C}^{{\mathcal Ab}_p}$
to the formal power series $f_A$ in $p$ commuting variables
is here completely natural (cf. Remark
\ref{majorationcriterionrem}). Does this have interesting
analytic consequences for
$f_A$ (which is holomophic in a neighbourhood of $(0,\dots ,0)$
by Proposition \ref{majorationcriterion})?
\end{rem} 


\section{Birecursivity} \label{birecursivity}

The {\it palindromic involution}
$$\iota(s_1s_2\dots s_{n-1}s_n,
t_1\dots t_n)=(s_ns_{n-1}\dots s_2s_1,t_n\dots t_1)$$
defines an involutive antiautomorphism of the monoid
${\mathcal M}_{p\times q}$ and we get the {\it palindromic automorphism}
$\iota:{\mathbf K}^{{\mathcal M}_{p\times q}}\longrightarrow 
{\mathbf K}^{{\mathcal M}_{p\times q}}$ by considering 
the involutive automorphism defined by
$(\iota X)[U,W]=X[\iota(U,W)]$ for $X\in{\mathbf K}^{
{\mathcal M}_{p\times q}}$. Since we have the identity
$(\iota X)(\iota(Y)=\iota(XY),X\in {\mathbf K}^{{\mathcal M}_{p\times
    r}},
Y\in {\mathbf K}^{{\mathcal M}_{r\times
    q}}$, the palindromic automorphism defines an involutive
automorphic functor
of the category ${\mathbf K}^{\mathcal M}$.

In characteristic $\not=2$, the palindromic automorphism $\iota$ 
endows the algebra ${\mathbf K}^{{\mathcal M}_{p\times q}}$
with a ${\mathbb Z}/2{\mathbb Z}-$ grading in the usual way
by considering the decomposition $X=X_++X_-$ into its {\it even} and
{\it odd} (palindromic) parts defined by 
$$X_+=\frac{X+\iota X}{2}\hbox{ and }X_-=\frac{X-\iota X}{2}$$
for $X\in {\mathbf K}^{{\mathcal M}_{p\times q}}$
and we have the sign-rules $$(XY)_+=X_+Y_++X_-Y_-,
(XY)_-=X_+Y_-+X_-Y_+$$ whenever the matrix product $XY$ is defined. 
In particular, we get a subcategory 
$\big({\mathbf K}^{{\mathcal M}}\big)_+$ consisting of
all even parts in ${\mathbf K}^{{\mathcal M}}$.

Conjugating the shift-monoid $\rho({\mathcal M}_{p\times q})$ by 
$\iota$, we get
a second morphism $\lambda:{\mathcal M}_{p\times q}
\longrightarrow \hbox{End}({\mathbf K}^{{\mathcal M}_{p\times q}})$,
called the {\it left-shift-monoid}. It is defined by
$$(\lambda(s_1\dots s_n,t_1\dots t_n)X)[U,W]=X[s_n\dots s_1 U,
t_n\dots t_1 W].$$
The obvious commutation rule 
$\rho(S,T)\lambda(S',T')=\lambda(S',T')\rho(S,T)$
yields an action of the (direct) product-monoid 
${\mathcal M}_{p\times q}\times {\mathcal M}_{p\times q}$
on ${\mathbf K}^{{\mathcal M}_{p\times q}}$. 
An element $X\in {\mathbf K}^{{\mathcal M}_{p\times q}}$ is a
{\it birecurrence matrix} if the linear span $\overline{X}^{birec}$
of the orbit $\lambda({\mathcal M}_{p\times q})
\rho({\mathcal M}_{p\times q})X$ has finite dimension.
The {\it birecursive complexity} of $X$ is defined as 
$\hbox{dim}(\overline{X}^{birec})\in {\mathbb N}\cup \{\infty\}$.

\begin{prop} \label{propbirec} 
We have
$$\hbox{dim}(\overline A^{birec})=\hbox{dim}({\mathbf
  K}[\rho_{\mathcal A}({\mathcal M}_{p\times q})])$$
with ${\mathbf
  K}[\rho_{\mathcal A}({\mathcal M}_{p\times q})]\subset \hbox{End}(
{\mathcal A})$ denoting the subalgebra of $\hbox{End}(
{\mathcal A})$ generated by the shift-monoid acting on 
the recursive closure ${\mathcal A}=\overline A^{rec}$ of $A$.  

In particular, we have the inequalities
$$\hbox{dim}(\overline{A}^{rec})\leq
\hbox{dim}(\overline{A}^{birec})\leq
\big(\hbox{dim}(\overline{A}^{rec})\big)^2.$$
\end{prop}

\begin{rem} The proof shows in fact that the action of the
  shift-monoid $\rho({\mathcal M}_{p\times q})$ on
  ${\overline A}^{birec}$ is (isomophic to) the obvious linear action by 
right multiplication of $\rho_{\mathcal A}({\mathcal M}_{p\times q})$ on $
{\mathbf K}[\rho_{\mathcal A}({\mathcal M}_{p\times q})]$.

In particular, if $pq\geq 2$ we have generically the equality
$$\hbox{dim}(\overline{A}^{birec})=
\big(\hbox{dim}(\overline{A}^{rec})\big)^2$$
for $A\in\hbox{Rec}_{p\times q}({\mathbb C})$.
\end{rem}

\begin{cor} \label{corbirec} (i) The vector spaces of recurrence matrices
and birecurrence matrices
in ${\mathbf K}^{{\mathcal M}_{p\times q}}$ coincide.

\ \ (ii) The
palindromic automorphism $\iota$ of ${\mathbf K}^{{\mathcal M}_{p\times
      q}}$ restricts to an involutive automorphism of
$\hbox{Rec}_{p\times q}({\mathbf K})$. The
decomposition of a recurrence matrix $A$ into its even and odd
parts $A_+=(A+\iota A)/2,\ A_-=(A-\iota A)/2$ holds thus in
$\hbox{Rec}
({\mathbf K})$ for ${\mathbf K}$ of characteristic $\not= 2$.

\ \ (iii) The palindromic automorphism yields an involutive
automorphic functor of the category $\hbox{Rec}({\mathbf K})$.
\end{cor}

\begin{rem} Defining the {\it palindromic element}
$P_p^\iota\in{\mathbf K}^{{\mathcal M}_{p\times p}}$ by
$P_p^\iota[s_1\dots s_n,s_n\dots s_1]=1$ and
$P_p^\iota[s_1\dots s_n,t_1\dots t_n]=0$ if $s_1\dots s_n\not=t_n\dots
t_1$, the palindromic automorphism of ${\mathbf K}^{{\mathcal
    M}_{p\times q}}$ or $\hbox{Rec}_{p\times q}$ is given by 
$X\longmapsto P_p^\iota XP_q^\iota$. In particular, the palindromic
automorphism is an inner automorphism of the algebra
${\mathbf K}^{{\mathcal M}_{p\times p}}$.

The palindromic element
$P_p^\iota$ is however not a recurrence matrix if $p\geq 2$.
Indeed, we have $(\rho(s_1\dots s_n,t_1\dots t_n)P_p^\iota)[
u_1\dots u_n,w_1\dots w_n]=1$ if $u_1\dots u_n=t_n\dots t_1,
w_1\dots w_n=s_n\dots s_1$ and 
$(\rho(s_1\dots s_n,t_1\dots t_n)P_p^\iota)[
u_1\dots u_n,w_1\dots w_n]=0$ otherwise. This implies 
$\hbox{dim}(\overline{P_p^\iota}^{rec})\geq p^{2n}$ for all 
$n\in{\mathbb N}$ and shows that the palindromic automorphism
of the algebra $\hbox{Rec}_{p\times p}({\mathbf K})$
is not inner for $p\geq 2$. (For $p=1$ the palindromic automorphisms
$P_1^\iota$ is trivial and thus inner in the commutative algebra
$\hbox{Rec}_{1\times 1}({\mathbf K})$.)

Let us finish this remark by adding that the algebras 
$\hbox{Rec}_{p\times p}({\mathbf K})$ admit many more similar
``exterior'' automorphisms. An example is for instance given
by the automorphism $\varphi_k$ (for $k\geq 1$) which acts as $\iota$
on $A[{\mathcal M}_{p\times p}^{nk}]$ and as the identity
on $A[{\mathcal M}_{p\times p}^{n}]$ if $n$ is not divisible by $k$.
It would perhaps be interesting to understand the algebraic
structure of the group of {\it outer automorphism}
which is defined as the quotient of all automorphisms of the algebra
$\hbox{Rec}_{p\times
  p}({\mathbf K})$ by the normal subgroup $\hbox{GL}_{p-rec}({\mathbf
  K})/(\hbox{GL}_{p-rec}({\mathbf K})\cap {\mathcal C}_{p-rec})$
(where ${\mathcal C}_{p-rec}$ denotes the center of $\hbox{Rec}_{p\times
  p}({\mathbf K})$) of inner automorphisms given by conjugations.
\end{rem}

{\bf Proof of Proposition \ref{propbirec}} Choose a presentation
$$(A_1=A,A_2,\dots)[\emptyset,\emptyset],\rho({\mathcal M}_{p\times
  q}^1)\subset\hbox{End}({\mathcal A})$$ of ${\mathcal A}=\overline 
A^{rec}$. The right action $\rho$ of ${\mathcal M}_{p\times q}$ 
on $A$ is obvious and Proposition \ref{propcoeffs} shows that 
$\tilde A=\lambda(s,t)A$ is presented by 
$$(\tilde A_1=\tilde A,\tilde A_2,\dots)[\emptyset,\emptyset]=
\left( (A_1=A,A_2,\dots)[\emptyset,\emptyset]\right)\rho(s,t)$$
with shift-matrices $\rho({\mathcal M}_{p\times q}^1)$ as above.

The right action $\rho({\mathcal M}_{p\times q})$ of the shift-monoid 
on the birecursive closure ${\overline A}^{birec}$ 
is thus (up to conjugacy) given by the right 
multiplication of $\rho_{\mathcal A}({\mathcal M}_{p\times q})$ on the 
algebra ${\mathbf K}[\rho_{\mathcal A}({\mathcal M}_{p\times q})]
\subset \hbox{End}({\mathcal A})$ (where ${\mathcal A}=
{\overline A}^{rec}$). This implies the equality
$$\hbox{dim}(\overline A^{birec})=\hbox{dim}({\mathbf
  K}[\rho_{\mathcal A}({\mathcal M}_{p\times q})])\ .$$
The inequalities are now obvious. \hfill$\Box$

\begin{rem} The content of this chapter can roughly be
paraphrased as follows: Recursively closed subspaces of 
${\mathbf K}^{{\mathcal M}_{p\times q}}$ are direct sums of 
linear vector spaces spanned by orbits of $\rho({\mathcal M}_{p\times
  q})$. Birecursively closed vector spaces are suitable
linear representations of the monoid algebra 
${\mathbf K}[{\mathcal M}_{p\times q}]$. Notice that not every linear
representation corresponds to a birecursively closed subspace
of ${\mathbf K}^{{\mathcal M}_{p\times q}}$; for instance
direct sums containing (up to isomorphism) a common irreducible 
linear subrepresentation of ${\mathcal M}_{p\times q}$ are forbidden.
\end{rem} 

\section{Virtual representations and  birecursively
closed  subspaces}

A finite-dimensional
birecursively closed subspace $\mathcal A\subset \hbox{Rec}_{p\times
  p}(\mathbf K)$ gives rise to a finite-dimensional
representation $\rho_{\mathcal A}:\mathcal M_{p,p}\longrightarrow 
\hbox{End}(\mathcal A)$. Choosing a basis of $\mathcal A$
yields a finite-dimensional linear representation of 
$\rho_{\mathcal A}(\mathcal M_{p\times p})$. 
Considering the equivalence relation on 
representations 
given by conjugation, birecursively closed subspaces
of $\hbox{Rec}_{p\times p}(\mathbf K)$ are in bijection 
with suitable finite-dimensional
matrix-representations of $\mathcal M_{p\times p}$.
We denote by $\hbox{Rep}$ the set of equivalence classes
of all finite-dimensional representations of $\mathcal M_{p\times p}$
and introduce the free vector space $\mathbf K[\hbox{Rep}]$
with basis $\hbox{Rep}$ and elements formal linear
combinations of finite-dimensional representations of 
$\mathcal M_{p\times p}$. Setting $\tau=\rho\sigma$ for
$\rho,\sigma\in \hbox{Rep}$
where
$$\tau(s,t)_{kl,ij}=
\sum_{u=1}^r \rho(s,u)_{k,i}\sigma(u,t)_{l,j}\
,$$
cf. Proposition \ref{presentproduct}, defines an associative
bilinear product on $\mathbf K[\hbox{Rep}]$ and turns it into 
an associative algebra. 

This algebra has two interesting quotients. The first one is 
given by considering the quotient by the ideal generated
by $\sigma=\rho\oplus \tau$ if the representation $\sigma$
is the direct sum of subrepresentations $\rho,\tau$. Its
elements, sometimes called {\it virtual representations}
are thus linear combinations of indecomposable
(but not necessarily irreducible) finite-dimensional
representations of the monoid $\mathcal M_{p\times p}$.

Another quotient of $\mathbf K[\hbox{Rep}]$ is obtained 
by considering the free vector space generated by all
finite-dimensional birecursively closed subspaces in
$\hbox{Rec}_{p\times p}(\mathbf K)$. Given two 
such subspaces $\mathcal A,\mathcal B$ the product
$\mathcal C=\mathcal A\mathcal B$ is the smallest 
birecursively closed subspace containing (and in fact
spanned by) all products $XY$ with $X\in\mathcal A$ and 
$Y\in\mathcal B$.

\begin{rem} The usual tensor product of matrices endows the
vector-space $\mathbf K[\hbox{Rep}]$ with a different,
commutative (and associative) product.
\end{rem}

\begin{rem} All constructions can also be carried out on
the quotient space $\hbox{Rec}_{p\times p}(\mathbf
K)/\mathcal{FS}_{p\times p}$. 

One can of course also consider not necessarily finite-dimensional
representations and subspaces. In this case one might also work
in the full algebra ${\mathbf K}^{\mathcal M_{p\times p}}$.

Let us also mention the obvious algebra structure on the
free vector space spanned by
all recursively (but not necessarily birecursively) closed
finite-dimensional subspaces of $\hbox{Rec}_{p\times p}(\mathbf K)$
(or of  $\hbox{Rec}_{p\times p}(\mathbf K)/\mathcal{FS}_{p\times p}$).

Computations in most algebras discussed in this paragraph can
be done algorithmically and involve only finitely many operations on
finite amounts of data.
\end{rem}


\section{Finite monoids and finite state automata}

This section describes a connection between 
automatic functions or automatic sequences 
(defined by finite state automata) 
and recurrence matrices with finite shift-monoid.

\subsection{Finite monoids}

\begin{prop} \label{finitevalues} Given a recurrence matrix $A\in 
\hbox{Rec}_{p\times q}({\mathbf K})$, 
the following assertions are equivalent:

\ \ (i) The subset $\{A[U,W]\ \vert\ (U,W)\in{\mathcal M}_{p\times q}\}
\subset {\mathbf K}$ of values of $A$ is finite.

\ \ (ii) For all $B\in \overline{A}^{rec}$, 
the subset $\{B[U,W]\ \vert\ (U,W)\in{\mathcal M}_{p\times q}\}
\subset {\mathbf K}$ of values of $B$ is finite.

\ \ (iii) The recursive set-closure $\rho({\mathcal M}_{p\times q})A$
of $A$ is finite.

\ \ (iv) The shift-monoid 
$\rho_{\overline{A}^{rec}}
({\mathcal M}_{p\times q})\in \hbox{End}(\overline{A}^{rec})$ of 
$A$ is finite.
\end{prop}

{\bf Proof of Proposition \ref{finitevalues}} For a recurrence matrix
$A\in\hbox{Rec}_{p\times q}({\mathbf K})$ satisfying assertion (i),
consider the finite set
${\mathcal F}=\{A[U,W]\ \vert\ (U,W)\in {\mathcal M}_{p\times q}\}
\subset
{\mathbf K}$ containing all evaluations of $A$.
The set ${\mathcal F}$ contains thus also all evaluations
$\{X[U,W]\ \vert\ (U,W)\in{\mathcal M}_{p\times q}\}$ for
$X\in \rho({\mathcal M}_{p\times q})A$ in the set-closure of $A$.
Choose $(U_1,W_1),\dots,(U_a,W_a)\in {\mathcal M}_{p\times q}$ such
that $\rho(U_1,W_1)A,\dots,\rho(U_a,W_a)A$ form a basis of 
$\overline{A}^{rec}$.
An arbitrary element $B\in \overline{A}^{rec}$ can thus be written
as a linear combination
$B=\sum_{j=1}^a \beta_j \rho(U_j,W_j)A$ and the set
$\{B[U,W]\ \vert\  (U,W)\in{\mathcal M}_{p\times q}\}$ of all 
its evaluations is a subset of the finite set
$\sum_{j=1}^a \beta_j {\mathcal F}$ containing at most
$1+a\sharp\{{\mathcal F}\setminus \{0\}\}<\infty$
elements. This shows the equivalence of assertions (i) and (ii).

All elements in the set closure $\rho({\mathcal M}_{p\times q})A\subset
\overline{A}^{rec}$ of $A\in\hbox{Rec}_{p\times q}({\mathbf K})$
are determined by their projection onto $\overline{A}^{rec}[
{\mathcal M}_{p\times q}^{\leq N}]$ where $N$ denotes the 
saturation level of
$\overline{A}^{rec}$. Assuming (i), the elements of $
\rho({\mathcal M}_{p\times q})A$
are thus in bijection with a subset of the finite set ${\mathcal F}^
{{\mathcal M}_{p\times q}^{\leq N}}$ where 
${\mathcal F}=\{A[U,W]\ \vert\ (U,W)\in{\mathcal M}_{p\times q}\}
\subset {\mathbf K}$ is the
finite set of all evalutions of $A$. This proves that (i) implies (iii).

The shift-monoid $\rho_{\overline{A}^{rec}}
({\mathcal M}_{p\times q})\subset
\hbox{End}(\overline{A}^{rec})$ of $A$ acts faithfully on the set-closure 
$\rho({\mathcal M}_{p\times q})A$ of $A$. 
Finiteness of $\rho({\mathcal M}_{p\times q})A$
implies thus (iv).

We have $\{A[U,W]\ \vert\ (U,W)\in{\mathcal M}_{p\times q}\}=
\{(\rho(U,W)A)[\emptyset,\emptyset]\ \vert\ 
(U,W)\in{\mathcal M}_{p\times q}\}$ which shows that (iv) implies (i).
\hfill$\Box$

\subsection{Automata} 

An {\it initial automaton} with {\it input alphabet} $X$ and 
{\it output alphabet} $Y$ is a directed graph ${\mathcal A}$
with vertices or {\it states} ${\mathcal V}$ such that
\begin{itemize}
\item ${\mathcal V}$ contains a marked {\it initial state} $v_*$.
\item The set of directed edges originating in a given state 
$v\in{\mathcal V}$ is bijectively labelled by the input alphabet $X$.
\item We have an {\it output function} 
$w:{\mathcal V}\longrightarrow Y$ on the set ${\mathcal V}$ of states.
\end{itemize}

An initial automaton ${\mathcal A}$ is {\it finite state} if its
input alphabet $X$ and its set of states ${\mathcal V}$ are both
finite. 

An initial automaton defines an application
$\alpha:{\mathcal M}_X\longrightarrow Y$ from the free monoid
${\mathcal M}_X$ generated by the input alphabet $X$ into $Y$.
Indeed, an element $x=x_1\dots x_n\in {\mathcal M}_X$ corresponds
to a unique continuous path $\gamma
\subset {\mathcal A}$ of length $n$ starting at $v_*$ and running
through $n$ directed edges labelled $x_1,x_2,\dots,x_n$.
We set $\alpha(x_1\dots x_n)=w(v_x)$ where $v_x$ is the endpoint of
the path 
$\gamma$. Every function $\alpha\in Y^{{\mathcal M}_X}$ can be 
constructed in this way by a suitable initial 
automaton which is generally 
infinite. A function $\alpha\in Y^{{\mathcal M}_X}$ is {\it
automatic} if it is realizable by a 
finite state initial automaton. Automatic functions (with
$Y\subset {\mathbf K}$ a subset of a field) which are converging
elements of ${\mathbf K}^{{\mathcal M}_p}$ are called 
{\it $p-$automatic} or {\it automatic sequences}.

\begin{prop} \label{automatic=recursiv}
Automatic functions with $X={\mathcal M}_{p}$
and $Y={\mathbf K}$ are in bijection with elements of
$\hbox{Rec}_{p}({\mathbf K})$ having finite shift-monoids.
\end{prop}

{\bf Proof} We show first that an automatic function
$\alpha\in{\mathbf K}^{{\mathcal M}_{p}}$ yields a 
recurrence vector with
finite shift-monoid. Choose a finite state initial automaton
${\mathcal A}=(v_*\in{\mathcal V},w\in {\mathbf K}^{\mathcal V})$
realizing $\alpha$ and having $a=\sharp({\mathcal V})$ states. 
Given a word
$S=s_1\dots s_n\in{\mathcal M}_{p}^n$,
we get an application $\gamma_{S}:{\mathcal V}\longrightarrow
{\mathcal V}$ defined by $\gamma_{S}(v)=u$ if the path
of length $n$ starting at state
$v\in {\mathcal V}$ and consisting of directed edges
 with labels $s_1,\dots,s_n,$ ends at state $u$.
Associate to a word $S$ the initial finite state automaton
${\mathcal A}_{S}=(v_*\in {\mathcal V},\tilde w=w\circ \gamma_{S})$
with the same set ${\mathcal V}$ of states and the same initial state
$v_*$ but with output function $\tilde w(v)=w(\gamma_{S}(v))\in
  {\mathbf K}$. We have then by construction
  $\alpha(US)=\alpha_{S}
(U)$ where $\alpha_{S}$ is the automatic function of the
finite state initial automaton ${\mathcal A}_{S}$.
 
Given a fixed finite state initial automaton ${\mathcal A}$,
the automatic function $\alpha_{S}$ depends only on the function
$\gamma_{S}$ which belongs to the finite set 
${\mathcal V}^{\mathcal V}$ of all $a^a$ functions from ${\mathcal V}$
to ${\mathcal V}$. The recursive set-closure 
$$\rho({\mathcal M}_p)\alpha=\{\alpha_{S}\ \vert
\ S\in{\mathcal M}_{p}\}$$
of $\alpha\in{\mathbf K}^{{\mathcal M}_{p}}$ is thus finite.
This proves that $\alpha\in\hbox{Rec}_{p}({\mathbf K})$ and 
Proposition \ref{finitevalues} implies that the shift-monoid 
of $\alpha$ is finite.

Consider now a recurrence matrix $A\in{\mathbf K}^{{\mathcal M}_{p}}$ 
with finite shift-monoid. We have to show that $U\longmapsto
A[U]$ (for $U\in {\mathcal M}_p$) 
is an automatic function. If $A$ is identically zero, there
is nothing to prove. Otherwise, complete $A$ to a basis
$A_1=A,A_2,\dots, A_a$ of its recursive closure $\overline{A}^{rec}$
of dimension $a$ and consider the associated
presentation with shift-matrices
$\rho_{\overline{A}^{rec}}(s)\in {\mathbf K}^{a\times a}$ 
and initial values
$W=(A_1[\emptyset],\dots,A_a[\emptyset])\in
{\mathbf K}^a$. We get thus a matrix-realization 
${\mathcal M}_A=\rho_{\overline{A}^{rec}}({\mathcal M}_{p})\subset 
\hbox{End}({\mathbf K}^a)$ of the finite shift-monoid of $A$. 
Consider now the Cayley graph $\Gamma$ of ${\mathcal M}_A$
with respect to the generators $\rho_{\overline{A}^{rec}}(s)\in 
\rho_{\overline{A}^{rec}}({\mathcal
  M}_{p}^1)$. The graph $\Gamma$ is the finite graph
with vertices indexed by all elements of ${\mathcal M}_A$. 
An oriented (or directed) edge labelled $s$ joins $X$ to $Y$ if
$Y=X\rho_{\overline{A}^{rec}}(s)$. 
Consider the graph $\Gamma$ as an initial
finite state automaton with initial state the identity
$\rho_A(\emptyset)$ and weight function
$w(X)=\big(X^tW^t\big)_1$ where $\big(X^tW^t\big)_1$ is the
first coordinate of the row-vector $X^tW^t\in {\mathbf K}^a$.
By construction, the automatic function $\alpha\in {\mathbf
  K}^{{\mathcal M}_{p}}$ associated to this initial
automaton is given by $\alpha(U)=A_1[U]$ since we have
$$\left(\begin{array}{c}
A_1[u_1\dots u_n]\\A_2[u_1\dots u_n]  
\\ \vdots\\ A_a[u_1\dots u_n]\end{array}
\right)=\rho_{\overline{A}^{rec}}(u_n)^t\cdots 
\rho_{\overline{A}^{rec}}(u_1)^t
\left(\begin{array}{c}
A_1[\emptyset]\\A_2[\emptyset]  \\ \vdots\\ 
A_a[\emptyset]\end{array}\right)$$
(cf. Proposition \ref{propcoeffs}).
 \hfill $\Box$

\begin{rem} When dealing with automatic functions, the output 
alphabet is often a finite field ${\mathbf K}$. 
Shift-monoids of recurrence matrices are then always finite and 
${\mathbf K}-$valued automatic functions on ${\mathcal
  M}_{p}$ are in bijection with $\hbox{Rec}_{p}({\mathbf K})$. Automatic
functions form thus a ring with respect to the Hadamard product 
(Proposition \ref{otherproducts}, assertion (i)) 
or with respect to the polynomial product (see \S 
\ref{sectionpolynomialring}). 
Converging automatic functions (also called automatic sequences,)
form a subring for the Hadamard product (see also
\cite{AS}, Corollary 5.4.5 and Theorem 12.2.6) and a differential
subring (after identification with the corresponding generating
series, see Proposition \ref{propdiffring}) 
for the polynomial product (see also \cite{AS}, Theorem 16.4.1).
\end{rem}


\section{The categories of transducers and finite-state
transducers}

The {\it transducer} of an element $A\in Y^{{\mathcal M}_X}$ 
is the length-preserving
application $\tau_A:{\mathcal M}_X\longrightarrow {\mathcal M}_Y$
defined by $\tau_A(\emptyset)=\emptyset$ and
$$\tau_A(u_1\dots u_n)=A[u_1u_2\dots u_n]A[u_2\dots u_n]A[u_3\dots
u_n]
\dots A[u_{n-1}u_n]A[u_n].$$
A transducer $\tau_A$
is in general not a morphism of monoids (except if 
$A[u_1\dots u_n]=A[u_1]$ for all $u_1\dots u_n\in {\mathcal M}_X$ of
length $\geq 1$) and the identity $\tau_A=\tau_B$ holds
if and only if $A[U,W]=B[U,W]$ for all $(U,W)\in\mathcal M_X\setminus
\emptyset$.

\begin{prop} \label{transducerdescr}
A length-preserving function $\tau:{\mathcal M}_X
\longrightarrow {\mathcal M}_Y$ is a transducer if and only if 
$\tau(uU)$ is of the form $*\tau(U)$ for all $u\in X$ and $U\in
{\mathcal M}_X$.
\end{prop}

\begin{cor} \label{transducercor}
The composition $\tau_B\circ \tau_A:{\mathcal M}_X
\longrightarrow{\mathcal M}_Z$ of two transducers
$\tau_A:{\mathcal M}_X\longrightarrow{\mathcal M}_Y,
\tau_B:{\mathcal M}_Y\longrightarrow{\mathcal M}_Z$ is a transducer.
\end{cor}

{\bf Proof of Proposition \ref{transducerdescr}} Given a 
length-preserving function
$\tau:{\mathcal M}_X\longrightarrow {\mathcal M}_Y$,
set $A[\emptyset]=\alpha$ with $\alpha\in Y$ arbitrary and
$A[uU]=y$ if $\tau(uU)=y\tau(U)$.
The transducer $\tau_A$ associated to $A$ coincides then with
the function $\tau$ if and only if $\tau$ satisfies the conditions
of Proposition \ref{transducerdescr}.\hfill$\Box$

The easy proof of Corollary \ref{transducercor} is left to the reader.
 
Corollary \ref{transducercor} allows to define the 
{\it category of transducers} by considering
free monoids ${\mathcal M}_p$ on $\{0,1,\dots, p-1\}$
as objects with arrows
from ${\mathcal M}_p$ to ${\mathcal M}_q$ given by
transducers $\tau_A:{\mathcal M}_p\longrightarrow{\mathcal M}_q$ 
associated to 
$A\in \{0,\dots,q-1\}^{{\mathcal
    M}_p}$. The set of all transducers from 
${\mathcal M}_p$ to ${\mathcal M}_q$ is in bijection with the set 
of functions $\{0,\dots, q-1\}^{{\mathcal M}_p\setminus\emptyset}$.

Given a transducer 
$\tau:{\mathcal M}_p\longrightarrow{\mathcal
  M}_q$, we define its {\it transducer-matrix}
$M_\tau\in {\mathbf K}^{{\mathcal M}_{q\times p}}$ by
$M_\tau[U,W]=1$ if $\tau(W)=U$ and $M_\tau[U,W]=0$ otherwise.
The set of transducer-matrices is contained in ${\mathbf
  K}^{{\mathcal M}_{q\times p}}$ for any field ${\mathbf K}$.
Since $M_{\tau\circ \tau'}=M_\tau M_{\tau'}$, the application
$\tau\longmapsto M_\tau$ which associates to a transducer
its transducer-matrix is a faithful functor from 
the category of transducers into the category
${\mathbf K}^{\mathcal M}$.

A transducer $\tau_A\in {\mathcal M}_q^{{\mathcal M}_p}$ 
is {\it finite-state} if
$A\in\{0,\dots, q-1\}^{{\mathcal M}_p}$ is automatic. 

\begin{prop} \label{fstransducerprop}
A transducer-matrix $M_\tau\in {\mathbf K}^{
{\mathcal M}_{q\times
      p}}$ is a recurrence matrix if and only if 
$\tau$ is a finite-state transducer.
\end{prop}

\begin{cor} Finite-state transducers form a subcategory 
in the category of transducers. 
\end{cor}

\begin{cor} The category of finite-state transducers can be realized
  as a subcategory of the category $\hbox{Rec}({\mathbf K})$ of
recurrence matrices over any field ${\mathbf K}$.
\end{cor}

{\bf Proof or Proposition \ref{fstransducerprop}} 
Consider a transducer $\tau=\tau_A:
{\mathcal M}_p\longrightarrow{\mathcal M}_q$ associated to
$A\in \{0,\dots,q-1\}^{{\mathcal M}_p}\subset {\mathbb Q}^{{\mathcal
      M}_p}$.
For $T\in {\mathcal M}_p$, we denote by $\tau_T$
the finite state transducer and by $M_T$
the transducer matrix associated to
the function $\rho(T)A$ (defined in the usual way by
$(\rho(T)A)[W]=A[WT]$). We have then
$$\tau(WT)=\tau_T(W)\tau(T)$$
where $\tau=\tau_\emptyset$. Moreover, the formula
$$(\rho(s,t)M_T)[U,W]=M_T[Us,Wt]=M[Us\tau(T),WtT]$$
shows the identities
$\rho(s,t)M_T=M_{tT}$ if $\tau(tT)=s\tau(T)$ (or, equivalently,
if $A[tT]=s$) and 
$\rho(s,t)M_T=0$ otherwise, for all $(s,t)\in {\mathcal
  M}_{q\times p}^1$.

If $A$ is automatic, then $A\in \hbox{Rec}_{p}({\mathbb Q})$
by Proposition \ref{automatic=recursiv}. 
Since $A$ takes all its values in the finite
set $\{0,\dots ,q-1\}$, Proposition \ref{finitevalues} 
implies finiteness of its recursive set-closure 
$\rho({\mathcal M}_{p})A=\{\rho(T)A\ \vert\ T\in {\mathcal
  M}_p\}$. This proves finiteness of the set 
${\mathcal S}=\{M_T\ \vert\ T\in{\mathcal M}_p\}\subset {\mathbf
  K}^{{\mathcal M}_{q\times p}}$. The recursive set-closure 
$\rho({\mathcal M}_{q\times p})M\subset \{0\cup{\mathcal S}\}$ (with
$0\in {\mathbf K}^{{\mathcal M}_{q\times p}}$ denoting the zero
recurrence matrix) is thus also finite and $M=M_\emptyset=M_\tau$ 
is a recurrence matrix.

In the other direction, we consider a transducer matrix $M=M_\tau
\in \hbox{Rec}_{q\times p}({\mathbf K})$
which is a recurrence matrix of a transducer $\tau=\tau_A:
{\mathcal M}_p\longrightarrow {\mathcal M}_q$ 
associated to $A\in \{0,\dots,q-1\}^{{\mathcal M}_p}$. 
The recurrence matrix $\rho(S,T)M$ is then either $0$ or the
transducer matrix of the transducer $\tau_{\rho(T)A}$ 
associated to $\rho(T)A$. Finiteness of $\rho({\mathcal M}_{q\times
  p})M$
(which follows from Proposition \ref{finitevalues}) implies finiteness
of the recursive set-closure $\rho({\mathcal M}_p)A$ and shows that
$A$ is a recurrence matrix with finite shift-monoid. The function $A$ is
thus automatic by Proposition \ref{automatic=recursiv}
and $\tau=\tau_A$ is an automatic transducer.\hfill $\Box$ 

\begin{prop} \label{transd-matrices}
A recurrence matrix $M\in {\mathbf K}^{{\mathcal
      M}_{q\times p}}$ is a transducer-matrix of a finite-state
transducer if and only if it can be given by a presentation with
initial values $(M=M_1,M_2,\dots
,M_d)[\emptyset,\emptyset]=(1,1,\dots,1)$ and shift matrices
$\rho(s,t)
\in\{0,1\}^{d\times d}$ with coefficients in $\{0,1\}$ such that
all column-sums of the $p$ matrices $\sum_{s=0}^{q-1}\rho(s,t)$
(with fixed $t\in\{0,\dots p-1\}$) are $1$.

Moreover, all finite-state transducers $\tau_1,\dots, \tau_d$
defined by $d$ transducer-matrices 
$M_1,\dots M_d$ as above are surjective if and only if all
column-sums of the  $q$ matrices $\sum_{t=0}^{p-1}\rho(s,t)$
are strictly positive.
\end{prop}

\begin{rem} There are $(qd)^{pd}$ presentations of complexity $d$
definining $d$ finite-state transducer-matrices in 
$\hbox{Rec}_{q\times p}$ as in Proposition \ref{transd-matrices}.
Indeed, these presentations are in bijection with matrices 
$\{0,1\}^{qd\times pd}$ having all column-sums equal to $1$ as can be seen
by contemplating the matrix
$$\left(\begin{array}{cccc}
\rho(0,0)&\rho(0,1)&\dots&\rho(0,p-1)\\
\rho(1,0)&\rho(1,1)&\dots&\rho(1,p-1)\\
\vdots&&&\vdots\\
\rho(q-1,0)&\rho(q-1,1)&\dots&\rho(q-1,p-1)\end{array}\right)$$
obtained by gluing all $qp$ shift-matrices $\rho(s,t)$ into
a $qd\times pd$ matrix. Since all $pd$ columns of such a matrix can
be choosen independently with $qd$ possibilities, there are
$(qd)^{pd}$ such matrices.
\end{rem}

{\bf Proof of Proposition \ref{transd-matrices}} Given a 
transducer-matrix $M\in \hbox{Rec}_{q\times p}({\mathbf K})$ 
associated to $A\in \hbox{Rec}_p({\mathbb Q})$,
we have $\rho(s,t)M=0$ if $A[t]\not=s$ and $\rho(s,t)M$ is the
transducer matrix $M_t$ associated to $\rho(t)A\in\hbox{Rec}_p(
{\mathbb Q})$ if $A[t]=s$. The first part of Proposition 
\ref{transd-matrices} follows by considering the presentation 
defined by all distinct transducer matrices in
$\{M_T\ \vert\ T\in{\mathcal M}_p\}=\rho({\mathcal M}_{q\times p})M
\setminus\{0\}$ where
$M_T$ is the transducer matrix associated to $\rho(T)A\in
\hbox{Rec}_p({\mathbb Q})$.

By induction on the length $l$ of words in ${\mathcal M}_p^l$,
the transducers associated to the transducer matrices $M_T$
are all surjective if for each $s\in
\{0,\dots, q-1\}$ there exists an integer $t=t_T,0\leq t<p$
such that $A[tT]=s$. The row corresponding to a transducer 
matrix $M_T$ of the shift-matrix $\rho(s,t_T)\in\{0,A\}^{d\times d}$
contains thus at least one non-zero coefficient and this proves 
the last part of Proposition \ref{transd-matrices} . \hfill $\Box$

\begin{rem} \label{transducerdescription}
Recursive presentation as introduced in chapter \ref{subsrecpres}
are particularly useful for transducer-matrices:
They admit recursive presentations of the form 
$A_i=(\rho_i,R(i)),i=1,\dots,d$ of depth $0$ such that 
all rows of all matrices $R(i)$ contain exactly one non-zero
element which belongs to the set $A_1,\dots A_d$.

As an example we consider the two finite-state
transducer-matrices $M_1,M_2\in \{0,1\}^{{\mathcal M}_{2\times 3}}$
recursively defined by
$$M_1=(1,\left(\begin{array}{ccc}M_1&M_2&0\\0&0&M_1\end{array}\right)),
\quad M_2=(1,
\left(\begin{array}{ccc}0&M_2&M_1\\M_2&0&0\end{array}\right)).$$
They span a recursively closed subspace in $\hbox{Rec}_{2\times 3}
({\mathbb Q})$ with (monoidal) presentation
given  by $(M_1,M_2)[\emptyset,\emptyset]=
(1,1)$ and
$$\begin{array}{l}
\displaystyle
\rho(0,0)=\left(\begin{array}{cc}1&0\\0&0\end{array}\right),
\rho(0,1)=\left(\begin{array}{cc}0&0\\1&1\end{array}\right),
\rho(0,2)=\left(\begin{array}{cc}0&1\\0&0\end{array}\right),\\
\displaystyle
\rho(1,0)=\left(\begin{array}{cc}0&0\\0&1\end{array}\right),
\rho(1,1)=\left(\begin{array}{cc}0&0\\0&0\end{array}\right),
\rho(1,2)=\left(\begin{array}{cc}1&0\\0&0\end{array}\right).\end{array}$$
\end{rem}

\begin{rem} \label{othertransducers}
Many authors define the transducer of a 
function $A\in Y^{{\mathcal M}_X}$ as the length-preserving
application $\tilde\tau_A:{\mathcal M}_X\longrightarrow {\mathcal M}_Y$
defined by $\tau_A(\emptyset)=\emptyset$ and
$$\tilde \tau_A(x_1\dots x_n)=A[x_1]A[x_1x_2]A[x_1x_2x_3]\dots A[x_1\dots
x_n].$$
This definition yields an equivalent theory since 
we have $\tilde \tau_A=\iota\circ \tau_A\circ \iota$ 
with $\iota$ denoting the palindromic 
antiautomorphism $w_1\dots w_n\longmapsto
w_n\dots w_n$ of ${\mathcal M}_X$ or ${\mathcal M}_Y$ considered in
chapter \ref{birecursivity}.
The corresponding 
transducer matrices $M_{\tau_A}$ and $M_{\tilde \tau_A}$ (for
$A\in \big({\mathcal M}_q^1\big)^{{\mathcal M}_p}$ are thus related by
$M_{\tilde \tau}=P_q^\iota M_\tau P_p^\iota$ where 
$P_n^\iota\in {\mathbf K}^{{\mathcal M}_{n\times n}}$ is the palindromic
involution defined in chapter \ref{birecursivity}.
By Proposition \ref{propbirec}, the ``transducer matrix'' of
$\tilde \tau_A$ is thus a recurrence matrix if and only if the transducer
matrix of $\tau_A$ is a recurrence matrix.
\end{rem}

\subsection{Transducers and regular rooted planar trees}

\begin{defin} The {\it $p-$regular rooted planar tree}
is the infinite tree ${\mathbf T}_p$
with vertices in bijection with 
${\mathcal M}_p$ and directed (or oriented) edges labelled 
by $s\in\{0,\dots, p-1\}$
joining a vertex $S\in{\mathcal M}_p$ to the vertex
$sS$. The vertex $\emptyset$ corresponds
to the root. 
\end{defin}

Geometrically, a vertex $s_1\dots s_l\in {\mathcal M}_p$ corresponds to
the endpoint of the 
continuous path of length $l$ starting at the root $\emptyset$
and running through $l$ consecutive oriented edges 
labelled $s_l,s_{l-1},\dots ,s_2,s_1$. 

A transducer $\tau:{\mathcal M}_p\longrightarrow 
{\mathcal M}_q$ induces an application $\tau:{\mathbf T}_p
\longrightarrow{\mathbf T}_q$ by its action on vertices.
In particular, a bijective transducer $\tau:{\mathcal M_p}
\longrightarrow {\mathcal M}_p$ corresponds to an automorphism
of the $p-$regular rooted tree ${\mathbf T}_p$ and the uncountable
group $\Gamma$ of all automorphisms of ${\mathbf T}_p$ is thus contained in
$\{0,1\}^{\mathcal M_{p\times p}}\subset \mathbf K^{\mathcal
  M_{p\times p}}$ where $\{0,1\}^{\mathcal M_{p\times p}}$ denotes the
set of all elements with values in $\{0,1\}$. This group contains
the countable subgroup $\Gamma\cap\hbox{Rec}_{p\times p}(\mathbf K)$
of all automorphisms corresponding to bijective transducers which
are automatic.


\section{Automatic groups}

\begin{prop} The inverse of a bijective transducer $\tau:
{\mathcal M}_p\longrightarrow {\mathcal M}_p$ is a transducer.

Moreover, if a bijective transducer $\tau:
{\mathcal M}_p\longrightarrow {\mathcal M}_p$ is finite-state,
then its inverse transducer $\tau^{-1}$ is also finite-state.
\end{prop}

{\bf Proof} Suppose first that $\tau$ is a finite state transducer
and consider a presentation $M_1,\dots, M_d$ as in Proposition
\ref{transd-matrices} with shift-matrices
$\rho(s,t)\in\{0,1\}^{d\times d}$ of the associated transducer-matrix
$M=M_\tau=M_1$. By Proposition
\ref{transd-matrices}, the sum of all coefficients of the matrix
$\sum_{0\leq s,t<p}\rho(s,t)\in{\mathbb Z}^{d\times d}$ equals $pd$.
The second part of Proposition
\ref{transd-matrices} shows thus that for fixed $s,0\leq s<p$,
each matrix
$\sum_{t=0}^p \rho(s,t)$ has all row-sums equal to $1$. 
Proposition
\ref{transd-matrices} can thus be applied to 
the transposed matrices $M_1^t,\dots M_t^t$
(given by the presentation 
$(M_1^t,\dots,M_d^t)[\emptyset,\emptyset]=(1,\dots,1)$
and shift-matrices $\tilde \rho(s,t)$ defined by $\tilde \rho(s,t)=
\rho(t,s)$) and shows that $M_1^t=M_1^{-1},\dots,M_d^t=M_d^{-1}$ are
transducer-matrices. Consider now the restriction $\tau({\mathcal
  M}_p^{\leq l})$ of an arbitrary transducer
$\tau$ to the set of words of length at most $l$ in ${\mathcal M}_p$.
This restriction coincides with the restriction to ${\mathcal
  M}_p^{\leq l}$ of a suitable finite-state transducer. 
Since the length $l$ can be an arbitrary integer this implies 
the first part of the result. 

The second part has already been proven.\hfill$\Box$

The set of all bijective finite-state transducers from 
${\mathcal M}_p$ to ${\mathcal M}_p$  forms a group which
is isomorphic to the subgroup $\mathcal{TD}_{p-rec}\subset
\hbox{GL}_{p-rec}({\mathbf K})$ formed by all invertible
transducer-matrices in $\hbox{Rec}_{p\times p}({\mathbf K})$. 
We call a subgroup of 
$\mathcal{TD}_{p-rec}$ a {\it $p-$automaton group} or 
{\it automaton group}
for short, also called branched groups by some authors.
For $M\in \mathcal{TD}_{p-rec}$, the matrices 
$M[{\mathcal M}_{p\times p}^l]$
are permutation matrices. We have thus $M^{-1}=M^t$ and the element
$M$ belongs to the orthogonal subgroup (with respect to 
the scalar product given by the identity) of 
$\hbox{GL}_{p-rec}({\mathbf K})$ (see section \ref{subsorthog}).

\begin{rem} There are $(d^p\ p!)^d$ presentations of complexity
$d$ defining $d$ bijective finite-state transducer-matrices 
$M_1,\dots ,M_d\in\hbox{Rec}_{p-rec}(\mathbb Z)$. Indeed, using
the obvious recursive presentation outlined in
Remark \ref{transducerdescription}, 
such a matrix $M_i=(1,R(i))$ is encoded by
a ``coloured'' permutation matrix $R(i)$ of order $p\times p$ with
all coefficients $1$ ``coloured'' (independently) by a colour in
the set $\{M_1,\dots,M_d\}$. For each such matrix $R(i)$, there are
thus $p!\ d^p$ possibilities. 
\end{rem}

\subsection{An example: The first Grigorchuk group $\Gamma$} 
This fascinating 
group appeared first in $\cite{Gri}$.
A few interesting properties (see for instance 
\cite{delaHarpe}, page 211 or \cite{GrPa}) 
of $\Gamma$ are:
\begin{itemize}
\item{} $\Gamma$ has no faithful finite-dimensional representation.
\item{} $\Gamma$ is not finitely presented.
\item{} $\Gamma$ is of intermediate growth.
\item{} $\Gamma$ contains every finite $2-$group.
\end{itemize} 

The group $\Gamma$ is the subgroup generated by four bijections 
$a,b,c,d$ of
the set ${\mathcal S}=\bigcup_{l=0}^\infty \{\pm 1\}^l$. Since $a,b,c,d$
preserve the subsets ${\mathcal S}_l=\{\pm 1\}^l$ we denote
by $a_l,b_l,c_l,d_l$ the restricted bijections induced by $a,b,c,d$ on 
the finite subset ${\mathcal S}_l$. For $l=0$, we have 
${\mathcal S}_0=\emptyset$ with trivial action of the permutations 
$a_0,b_0,c_0,d_0$. For $l>0$ we
write $(\epsilon,x)=(y_1,y_2,\dots,y_l)$ with $\epsilon=y_1\in\{\pm 1\}$
and $x=(y_2,y_3,\dots,y_l)\in\{\pm 1\}^{l-1}$. The action of $a_l,b_l,c_l,d_l$
is then recursively defined by
$$\begin{array}{ll}
\displaystyle a_l(1,x)=(-1,x)\qquad &a_l(-1,x)=(1,x)\\
\displaystyle b_l(1,x)=(1,a_{l-1}(x))\qquad &b_l(-1,x)=(-1,c_{l-1}(x))\\
\displaystyle c_l(1,x)=(1,a_{l-1}(x))\qquad &c_l(-1,x)=(-1,d_{l-1}(x))\\
\displaystyle d_l(1,x)=(1,x)\qquad &d_l(-1,x)=(-1,b_{l-1}(x))
\end{array}$$
(cf VIII.B, pages 
217-218 of \cite{delaHarpe}).

The set ${\mathcal S}$ corresponds to the vertices of the $2-$regular
rooted planar tree ${\mathbf T}_2$ and the four bijections $a,b,c,d$
act as automorphisms on ${\mathbf T}_2$. They correspond thus to
transducers $\tau_a,\tau_b,\tau_c,\tau_d$.
The associated transducer-matrices $M_a,M_b,M_c,M_d$ are recursively
presented by
$$M_a=(1,\left(\begin{array}{cc}0&\hbox{1d}\\\hbox{1d}&0\end{array}
\right)),
M_b=(1,\left(\begin{array}{cc}M_a&0\\0&M_c\end{array}\right)),$$
$$M_c=(1,\left(\begin{array}{cc}M_a&0\\0&M_d\end{array}\right)),
M_d=(1,\left(\begin{array}{cc}\hbox{Id}&0\\0&M_b\end{array}\right))$$
with $\hbox{Id}\in \hbox{GL}_{2-rec}(\mathbb Z)$ denoting the identity 
matrix recursivley presented by $\hbox{Id}=(1,\left(\begin{array}{cc}
\hbox{Id}&0\\0&\hbox{Id}\end{array}\right))$.

A (monoidal) presentation of $M_a,M_b,M_c,M_d$
is given by
$(\hbox{Id},M_a,M_b,M_c,M_d)[\emptyset,\emptyset]=(1,1,1,1,1)$ and
shift-matrices
$$\rho(0,0)=\left(\begin{array}{ccccc}
1&0&0&0&1\\
0&0&1&1&0\\
0&0&0&0&0\\
0&0&0&0&0\\
0&0&0&0&0\end{array}\right),\ 
\rho(0,1)=\left(\begin{array}{ccccc}
0&1&0&0&0\\
0&0&0&0&0\\
0&0&0&0&0\\
0&0&0&0&0\\
0&0&0&0&0\end{array}\right),$$
$$\rho(1,0)=\left(\begin{array}{ccccc}
0&1&0&0&0\\
0&0&0&0&0\\
0&0&0&0&0\\
0&0&0&0&0\\
0&0&0&0&0\end{array}\right),\ 
\rho(1,1)=\left(\begin{array}{ccccc}
1&0&0&0&0\\
0&0&0&0&0\\
0&0&0&0&1\\
0&0&1&0&0\\
0&0&0&1&0\end{array}\right).$$

\begin{rem} The fact that the four generators $a,b,c,d$ of 
the Grigorchuk group are of order $2$ is equivalent to the
identity $\rho(0,1)=\rho(1,0)$ of the shift-matrices described above.
\end{rem}


\section{A few asymptotic problems in $\hbox{Rec}_{p\times p}
({\mathbf K})$}

This section introduces generating series counting 
dimensions associated to recurrence matrices. It contains
mainly definitions and (open) problems.

Given a recurrence matrix $A\in\hbox{Rec}_{p\times p}({\mathbf K})$
one can consider the generating series 
$$1+\hbox{dim}(\overline{A}^{rec})t+\dots=
\sum_{n=0}^\infty \hbox{dim}(\overline{A^n}^{rec})t^n$$
associated to the complexities of its powers. One has the 
inequality $\hbox{dim}(\overline{A^n}^{rec})\leq
\left(\hbox{dim}(\overline{A}^{rec})\right)^n$ which
implies convergency of the series for $\{t\in{\mathbb C}\ \vert
\ \vert t\vert <1/\hbox{dim}(\overline{A}^{rec})\}$.

What can be said about the analytic properties of this generating
series? In a few easy cases (nilpotent or of complexity $1$ for 
instance), it is meromorphic in ${\mathbb C}$. 

If $A\in \hbox{GL}_{p-rec}({\mathbf K})$ is invertible in
$\hbox{Rec}_{p\times p}({\mathbf K})$, one can of course 
also consider the formal sum
$$\sum_{n=-\infty}^\infty \hbox{dim}(\overline{A^n}^{rec})t^n$$
encoding the complexities of all integral powers of $A$.
 
\begin{rem} There are many variations for the generating function(s)
defined above. One can replace the complexity by the
stable complexity or the 
birecursive complexity of powers. One can also consider
the generating function encoding the
dimension $\alpha_n$ of the recursive closure $\overline{\hbox{Id},
A^1,\dots, A^n}^{rec}$ containing all 
powers $A^0=\hbox{Id},A^1,\dots, A^n$ etc.
\end{rem}

Similarly, given a monoid or group generated by a finite
set ${\mathcal G}\subset \hbox{Rec}_{p\times p}({\mathbf K})$ of 
recurrence matrices, one can consider the generating series
$\sum_{n=0}^\infty \hbox{dim}({\mathcal A}_n)t^n$ encoding the 
complexities of the recursive closures ${\mathcal A}_n$
associated to all products of (at most) $n$ elements in ${\mathcal G}$.

A different kind of generating series is given by
considering for $A\in \hbox{Rec}_{p\times p}({\mathbb C})$
the series
$$1+\parallel A\parallel_\infty^\infty\ t+\dots=
\sum_{n=0}^\infty \parallel A^n\parallel_\infty^\infty \ t^n.$$
The inequality $\parallel A^n\parallel_\infty^\infty
\leq p^{2n-1}(\parallel A\parallel_\infty^\infty)^n$
(which can be proven by considering 
a recursive matrix whose coefficients $A[U,W]$
depend only on the length $l$ of $(U,W)\in{\mathcal M}_{p\times p}^l$)
shows again convergency for $t\in{\mathbb C}$ small enough.


\section{A generalization}

Consider a ring ${\mathbf R}$ of functions ${\mathbb N}\longrightarrow
{\mathbf K}$ with values in a commutative field ${\mathbf K}$. An
element $A\in{\mathbf K}^{{\mathcal M}_{p\times q}}$ is an
{\it ${\mathbf R}-$recurrence matrix} (or simply a {\it recurrence matrix}
if the underlying function ring ${\mathbf R}$ is obvious) if there
exists a finite number of elements $A_1=A,A_2,\dots, A_a\in 
{\mathbf K}^{{\mathcal M}_{p\times q}}$ and $pq$ {\it shift-matrices}
$\rho_A(s,t)\in{\mathbf R}^{a\times a}$ with 
coefficients $\rho_A(s,t)_{j,k}\in{\mathbf R}$ (for $1\leq
j,k\leq a$)
in the function ring ${\mathbf R}$ such that
$$\big(\rho(s,t)A_k\big)[{\mathcal M}_{p\times q}^{l+1}]=
\sum_{j=1}^a \rho_A(s,t)_{j,k}(l)\ A_j[{\mathcal M}_{p\times q}^l],\
1\leq k\leq a,$$
for all $0\leq s<p,0\leq t<q$.

We denote by ${\mathbf R}-\hbox{Rec}_{p\times q}({\mathbf K})$
the vector-space of all ${\mathbf
  R}-$recurrence matrices in ${\mathbf K}^{{\mathcal M}_{p\times q}}$.
 
Proposition \ref{propproduct} and its proof can easily 
be modified in order to deal with ${\mathbf
  R}-$recurrence matrices and we get thus a category 
${\mathbf R}-\hbox{Rec}({\mathbf K})$ with ${\mathbf R}-$recurrence 
matrices as morphisms. Moreover, 
the vector-space
${\mathbf R}-\hbox{Rec}_{p\times q}({\mathbf K})$ is invariant under
the action of the shift-monoid if the function ring 
${\mathbf R}$ is preserved by the translation $x\longmapsto x+1$ of
the argument 
(ie. if $\alpha\in{\mathbf R}$ implies that
the function $x\longmapsto \alpha(x+1)$ is also in ${\mathbf R}$).

A few analogies and differences between ${\mathbf R}-$recurrence matrices
and (ordinary) recurrence matrices are:

\begin{itemize}
\item{} Elements of ${\mathbf R}-\hbox{Rec}_{p\times q}({\mathbf K})$
have finite descriptions involving a finite number of elements in the
function ring ${\mathbf R}$.

\item{} The notion of a presentation is in general more involved:
Every presentation of an ordinary recurrence matrix $A$ contains
a basis of $\overline{A}^{rec}$ and can thus be used to 
construct a minimal presentation (defining a basis of 
$\overline{A}^{rec}$). This is no longer true in general for
${\mathbf R}-$recurrence matrices since ${\mathbf R}$ might
contain non-zero elements which have no multiplicative inverse.
One has thus to work with (not necessarily free)
${\mathbf R}-$modules when dealing with presentations.

\item{} Proposition \ref{propcoeffs}, slightly modified,
remains valid. We have
$$\left(\begin{array}{c}A_1[U,W]\\
\vdots\\A_d[U,W]\end{array}\right)=
\rho_{\mathcal A}(u_n,w_n)^t(n-1)\cdots \rho_{\mathcal A}(u_1,u_1)^t(0)
\left(\begin{array}{c}A_1[\emptyset,\emptyset]\\
\vdots\\ A_d[\emptyset,\emptyset]\end{array}\right)$$
where $(U,W)=(u_1\dots u_n,w_1\dots w_n)\in{\mathcal M}_{p\times
  q}^n$. The matrices $\rho_{\mathcal A}(s,t)
\in {\mathbf R}^{d\times d}$
are the obvious shift-matrices of the subspace ${\mathcal A}=\sum
{\mathbf K}A_j$ with respect to the generators $A_1,\dots,A_d$.

\item{} An analogue of the saturation level is no longer
available in general. This makes automated computations impossible
in the general case. 

\end{itemize}

\begin{rem} Non-existence of a saturation level does not necessarily
imply the impossibility of proving a few identities among
${\mathbf R}-$recurrence matrices, since such an identity
can perhaps be proven by induction on the word-length
$l$, after restriction to ${\mathcal M}_{p\times q}^{\leq l}$.
\end{rem}

\subsection{Examples}

The case where the function ring ${\mathbf R}$ consists of
all constant functions ${\mathbb N}\longrightarrow {\mathbf K}$
corresponds of course to the case of (ordinary) recurrence 
matrices studied in this paper.

The case where the function ring ${\mathbf R}$ consists 
of ultimately periodic functions 
${\mathbb N}\longrightarrow {\mathbf K}$ yields again only
ordinary recurrence matrices.

The first interesting new case is given by considering 
the ring ${\mathbf R}={\mathbf K}[x]$ 
of all polynomial functions. The category 
${\mathbb C}[x]-\hbox{Rec}({\mathbb C})$ of ${\mathbb
  C}[x]-$recurrence matrices is indeed larger than
the category $\hbox{Rec}({\mathbb C})$ of ordinary recurrence
matrices since ${\mathbb C}[x]-\hbox{Rec}_1({\mathbb C})$
contains for instance the element $A[0^n]=n!$ which is not in
$\hbox{Rec}_1({\mathbb C})$ by Proposition \ref{majorationcriterion}.

It would be interesting to 
have finiteness results in this case:
Given $A,B,C\in {\mathbf K}[x]-\hbox{Rec}_{p\times q}$ 
defined by presentations of complexity $a,b,c$
(where the complexity is the minimal number of
elements appearing in a finite presentation)
with shift-matrices involving only polynomials
of degree $\leq \alpha,\beta,\gamma$, can one give a bound
$N_+=N_+(a,b,c,\alpha,\beta,\gamma,p,q)$ such that 
the equality $(A+B)[{\mathcal M}_{p\times q}^{\leq N_+}]=
C[{\mathcal M}_{p\times q}^{\leq N_+}]$ 
ensures the equality $A+B=C$ in ${\mathbf K}[x]-\hbox{Rec}_{p\times
  q}$?
Similarly, given $A\in {\mathbf K}[x]-\hbox{Rec}_{p\times r},B
\in {\mathbf K}[x]-\hbox{Rec}_{r\times q},
C\in {\mathbf K}[x]-\hbox{Rec}_{p\times q}$ 
defined by presentations of complexity $a,b,c$
with shift-matrices involving only polynomials
of degree $\leq \alpha,\beta,\gamma$, can one give a bound
$N_\times=N_\times (a,b,c,\alpha,\beta,\gamma,p,q,r)$ such that 
the equality $(AB)[{\mathcal M}_{p\times q}^{\leq N_+}]=
C[{\mathcal M}_{p\times q}^{\leq N_+}]$ 
ensures the equality $AB=C$ in ${\mathbf K}[x]-\hbox{Rec}_{p\times
  q}$? Are there natural and ``interesting'' examples 
of ${\mathbf K}[x]-$recurrence matrices which are not
(ordinary) recurrence matrices?

Another interesting example is given by considering the ring
$\mathbf R$ of all functions ${\mathbb N}\longrightarrow {\mathbf K}$
which are linear combinations involving terms of the form
$$n\longmapsto n^k\lambda^n$$
for $\lambda\in {\mathbf K}$. For $p\geq 2$ and ${\mathbf K}$
algebraically closed, the group of lower triangular
convergent Toeplitz matrices in ${\mathbf
  R}-\hbox{GL}_{p-rec}({\mathbf K})$ contains then the multiplicative
group ${\mathbf K}(x)^*$ of all invertible rational power-series 
(this is not the case for ordinary recurrence matrices over
${\mathbb C}$, see the last part of Remark \ref{remnoinverse} and 
Example \ref{remintnoninvert}).

\subsection{An intermediate category between
$\mathbf K^{\mathcal M}$ and $\hbox{Rec}(\mathbf K)$} 
The category $\mathbf K^{\mathcal M}$ contains many other,
perhaps interesting, subcategories. An example is given by the
set of all elements $A\in\mathbf K^{\mathcal M_{p\times q}}$
(for arbitrary $p,q\in\mathbb N$) such that 
$$\hbox{dim}\left(\overline{A}^{rec}[\mathcal M_{p\times
    q}^l]\right)$$
is bounded by a polynomial in $l$. Proposition \ref{propproduct}
shows easily that this property is preserved by products. 
An example of such an element in $\mathbb Q^{\mathcal M_{2\times 2}}$
is perhaps given by the inverse element of the converging
element with limit the infinite Hadamard matrix
associated to the sequence of coefficients of $\prod_{k=0}^\infty
(1-x^{2^k})$.

\begin{rem} One can of course also consider suitable intermediate
growth classes for definining other subcategories of 
$\mathbf K^{\mathcal M}$.
\end{rem}


\section{Examples of a few Toeplitz matrices in
$\hbox{Rec}_{p\times p}(\mathbf K)$}

This and the next chapters present a few (hopefully) interesting 
recurrence examples of recurrence matrices, mainly elements of
$\hbox{GL}_{p-rec}(\mathbf K)$ for a suitable integer $p\geq 2$ 
and $K$ a subfield or subring of $\mathbb C$ or a finite
field. A few examples contain parameters which can be choosen
in suitable, easily specified rings or fields. 

Most of the 
computations are straightforward but tedious and omitted.

\subsection{}
Consider the matrix $A(n)$ of square size $n\times n$ with
coefficients  $A(n)_{i,j},0\leq i,j<n$ given by
$$A(n)_{i,j}=\left\lbrace\begin{array}{ll}
\displaystyle 4&\displaystyle i=j,\\
\displaystyle 3&\displaystyle i<j,\\
\displaystyle 3+3(i-j)\quad&\displaystyle i>j.\end{array}\right.$$
The infinite matrix $A(\infty)$ defines thus a converging element
(still denoted) $A\in\hbox{Rec}_{p\times p}(\mathbb Z)$ for all $p\in
\mathbb N$. We leave it to the reader to prove that $A=LU$
(where $L$ is lower triangular unipotent and $U$ is upper triangular)
with $L,U\in\hbox{GL}_{p-rec}(\mathbb Q)$ for all $p\in\mathbb N$.
An inspection of $U$ proves 
$$\det(A(n))=\left\lbrace\begin{array}{ll}
\displaystyle 1&\displaystyle n\equiv 0\pmod 3,\\
\displaystyle 4&\displaystyle n\equiv 1\pmod 3,\\
\displaystyle -2&\displaystyle n\equiv 2\pmod 3.\end{array}\right.$$
It follows that $A\in \hbox{GL}_{p-rec}(\mathbb Z)$ if 
$p\equiv 0\pmod 3$.

\begin{rem} More generally, one can consider the $n\times n$
matrix $A(n)$ with coefficients
$$A(n)_{i,j}=\left\lbrace\begin{array}{ll}
\displaystyle x+1&\displaystyle i=j,\\
\displaystyle x&\displaystyle i<j,\\
\displaystyle x+x(i-j)\quad&\displaystyle i>j\end{array}\right.$$
for $0\leq i,j<n$. It follows then for instance from 
\cite{ZaPe} that the sequence $\det(A(1)),\det(A(2)),\det(A(3)),\dots$
satisfies a linear recursion. More precisely we have
$$\det(A(n))+(x-3)\det(A(n-1))+(3-x)\det(A(n-2))-\det(A(n-3))=0$$
which implies
$$\det(A(n))=1-(-1)^n\sum_{k=0}^n (-1)^k{2n-1-k\choose k}x^{n-k}\ .$$
Particularly interesting are the evaluations $x\in \{1,2,3\}$
where the case $x=3$ has been considered above.
\end{rem}

\subsection{} A similar example is given by the 
$n\times n$ matrix $A(n)$ with coefficients
$$A(n)_{i,j}=\left\lbrace\begin{array}{ll}
\displaystyle 1&\displaystyle i=j,\\
\displaystyle 1&\displaystyle i<j,\\
\displaystyle 1+i-j\quad&\displaystyle i>j.\end{array}\right.$$
Its characteristic polynomial is 
$$t^n-\sum_{k=1}^n{n-1+k\choose n-k}t^{n-k}$$
and $A(n)$ is thus invertible over $\mathbb Z$ for all $n$
since $\det(A(n))=-(-1)^n$ for $n\geq 1$.

Setting $A[\mathcal M^l]=A(p^l)$ we get for all $p\in \mathbb N$ 
a converging element (still denoted) $A\in\hbox{GL}_{p-rec}(\mathbb Z)$
which satisfies
$A=LU$ with $L\in \hbox{GL}_{p-rec}(\mathbb Z)$ converging lower
triangular unipotent and $U\in \hbox{GL}_{p-rec}(\mathbb Z)$ 
converging upper triangular.

\subsection{}
The square matrix $A(2n)$ of even size $2n\times 2n$
with coefficients $A(2n)_{i,j},0\leq i,j<2n$ given by
$A(2n)_{i,j}=1$ if $(i-j)^2=1$
and $0$ otherwise has determinant $(-1)^n$. The inverse matrix 
$B(2n)$ of $A(2n)$ has coefficients $B(2n)_{i,j}=0$ 
if $\hbox{min}(i,j)\equiv 1\pmod 2$
and $B(2n)_{i,j}=[x^{\vert i-j\vert}]\frac{x}{1+x^2}$ otherwise (for $0\leq
i,j<2n$). For even $p\in 2\mathbb N$, this yields 
thus inverse matrices (still denoted) $A,B=A^{-1}
\in\hbox{Rec}_{p\times
  p}(\mathbb Z)$ defined by $A[\emptyset,\emptyset]=
B[\emptyset,\emptyset]=1$ and $A[\mathcal M_{p\times
  p}^l]=A(p^l)$, $B[\mathcal M_{p\times  p}^l]=B(p^l)$.
\subsubsection{}
Over $\mathbb F_2$ there is a very similar $2-$recursive example 
with coefficients $A_{i,j}=1$ except for $i=j$.

\subsection{}
Define a symmetric $n\times n$ Toeplitz matrix $A(n)$ by
$A(n)_{i,j}=\alpha_{\vert i-j\vert},0\leq i,j$ where
$$\alpha_0=0,\alpha_1=1,\alpha_2=-\frac{1}{2},\alpha_3=\frac{5}{4},
\alpha_4=-\frac{9}{8},\dots,\alpha_n=-\frac{1}{2}\alpha_{n-1}+
 \alpha_{n-2},\dots$$
are the coefficients of the rational series
$\sum_{n=1}^\infty \alpha_nx^n=\frac{x}{1+x/2-x^2}.$
For $n\geq 4$, the coefficients $B(n)_{i,j},0\leq i,j<n$ of
the inverse matrix $B(n)=A(n)^{-1}$ are then given by
$$\left\lbrace\begin{array}{cl}
\displaystyle 1&\displaystyle \hbox{if }i=j\in\{0,n-1\},\\
\displaystyle 5/4&\displaystyle \hbox{if }i=j\in\{1,n-2\},\\
\displaystyle 9/4&\displaystyle \hbox{if }i=j\in\{2,3,\dots,n-3\},\\
\displaystyle 1/2&\displaystyle \hbox{if }(i,j)\in\{(0,1),(1,0),(n-2,n-1),(n-1,n-2)\},\\
\displaystyle -1&\displaystyle \hbox{if }\vert i-j\vert=2,\\
\displaystyle 0&\displaystyle \hbox{otherwise}.\end{array}\right.$$
Choosing an odd prime $\wp$ and a natural number $p\in\mathbb N$, we
get thus an element (still denoted) $A\in\hbox{GL}_{p-rec}(\mathbb
F_\wp)$ defined by $A[\mathcal M^l]=A(p^l)\pmod \mathbb F_\wp$.

\begin{rem} We have The matrices $B(n)=A(n)^{-1}$ are converging 
for $n\rightarrow \infty$ in the sense that the all coefficients 
with fixed indices are ultimately constant. 
The limit-matrix $B(\infty)$ has an
$LU-$decomposition given by
$$LL^t=\left(\begin{array}{ccccccccc}
1&1/2&-1&0&0&0&0&\cdots\\
1/2&5/4&0&-1&0&0&0&\\
-1&0&9/4&0&-1&0&0\\
0&-1&0&9/4&0&-1&0&\\
\vdots&&\ddots&\ddots&\ddots&\ddots&\ddots\end{array}\right)$$
for 
$$L=\left(\begin{array}{cccccc}
1\\
1/2&1\\
-1&1/2&1\\
0&-1&1/2&1\\
&&\ddots&\ddots&\ddots\end{array}\right)$$
the lower triangular Toeplitz matrix defined by $L_{i,i}=1,
L_{i+1,i}=1/2,L_{i+2,i}=-1$ and $L_{i,j}=0$ otherwise.
The reduction of $L$ modulo an odd prime $\wp$ defines a
converging element $L\in\hbox{GL}_{p-rec}(\mathbb F_\wp)$.
Proposition \ref{majorationcriterion} shows however that
$L\in \hbox{Rec}_{p\times p}(\mathbb Q)$ is not invertible in
$\hbox{Rec}_{p\times p}(\mathbb Q)$.

Let us also mention a curious experimental fact concerning
the characteristic polynomial of $-A(2n)$.

A symmetric Toeplitz matrix of size $n\times n$ preserves the
eigenspaces of the involution
$$\iota:(x_0,\dots,x_{n-1})\longmapsto (x_{n-1},x_{n-2},\dots,x_1,x_0)$$
and its characteristic polynomial $\chi$ factorizes thus 
as $\chi=\chi_+\chi_-$ where $\chi_+$ of degree $\lceil n/2\rceil$ 
corresponds to the trivial eigenspace (formed by eigenvectors
$(x_0,\dots,x_{n-1})$ of eigenvalue $1$ satisfying $x_i=x_{n-1-i}$ for
$i=0,\dots,\lfloor n/2\rfloor$) and
where $\chi_-$ of degree $\lfloor n/2\rfloor$
corresponds to the eigenspace of eigenvalue $-1$
formed by all vectors $(x_0,\dots,x_{n-1})$ such that $x_i=-x_{n-1-i}$ for
$i=0,\dots,\lfloor n/2\rfloor$.

For the characteristic polynomial of $-A(2n)$ the corresponding
factorisation seems to be of the form
$$\det(t\hbox{ Id}+A(2n))=(t^n+\sum_{k=0}^{n-1}\gamma_{n,k}t^k)
  (t^n-\sum_{k=0}^{n-1}\gamma_{n,k}t^k)$$
with $\gamma_{n,0},\dots,\gamma_{n,n-1}$ strictly positive rational
numbers. In particular, the bilinear product 
defined by $A(2n)$ endows seemingly the trivial eigenspace 
(associated to the eigenvalue $1$) of
$\iota$ with an Euclidean scalar product.
Normalized in order to have leading term $4^{n-1}t^n$, the first few
divisors $\tilde \chi_+(n)$ of $\det(t\hbox{ Id}+A(2n))$
associated to the trivial eigenspace for $\iota$ are given by
$$\begin{array}{l}
\displaystyle \tilde \chi_+(1)=t+1\\
\displaystyle \tilde \chi_+(2)=4t^2+9t+4\\
\displaystyle \tilde \chi_+(3)=16t^3+65t^2+72t+16\\
\displaystyle \tilde \chi_+(4)=64t^4+441t^3+844t^2+432t+64\\
\displaystyle \tilde
\chi_+(5)=256t^5+2929t^4+8208t^3+7008t^2+2304t+256
\end{array}$$
and satisfy the recurrence relation
$$\tilde \chi_+(n)-(13t+4)\tilde\chi_+(n-1)+4t(13t+4)\tilde\chi_+(n-2)-
64t^3\tilde\chi_+(n-3)=0$$
(the same recurrence relation is also satisfied by the
corresponding complementary divisors, or equivalently, by the geometric
progression $1,(4t),(4t)^2,(4t)^3,\dots$).
Equivalently, for $k\geq 1$, 
the coefficient $t^{n-k}$ of $\chi_+(n)$ is given by the
coefficient of $x^n$ in the rational series $4^{k-1}(x/(1-9x+16x^2))^k$.
\end{rem}

\subsubsection{}
A similar example (defining also elements in $\hbox{GL}_{p-rec}
(\mathbb F_\wp)$ for $\wp$ an odd prime) is obtained by
considering
$A(n)_{i,j}=f_{\vert i-j\vert}, 0\leq i,j$ where
$$f_0=0,f_1=1,f_2=1,f_3=2,f_4=3,\dots, f_n=f_{n-1}+f_{n-2}$$
defined by $\sum_{n=1}^\infty f_nx^n=\frac{x}{1-x-x^2}$
is the Fibonacci sequence. We have then $\det(A(n))=-(-1)^{n-2}$ for
$n\geq 2$.

\subsubsection{}
Another similar example (satisfying
$\det(A(n))=-(-1)^n(\lambda^2-1)^{n-1}$ for $n\geq 1$ and giving 
non-trivial elements in $\hbox{GL}_{p-rec}(\mathbb F_\wp)$
for $\wp$ a suitable prime depending on $\lambda\in\mathbb Q$)
is given by $A(n)_{i,j}=\lambda^{\vert i-j\vert}$.

\subsection{}
Given an even integer $n$ such that $n\equiv 0\pmod 4$, consider the
symmetric 
Toeplitz matrix $A(n)$ of square size $n\times n$ 
with coefficients given by
$$A(n)_{i,j}=\left\lbrace\begin{array}{ll}
\displaystyle (2-n)/2+a&\displaystyle\hbox{if }i=j\\
\displaystyle -n/2+\vert i-j\vert+a&\displaystyle\hbox{if
}i\not=j\end{array}\right.$$
The matrix $A(n)$ is then invertible with inverse matrix
having coefficients $\left(A(n)^{-1}\right)_{i,j},0\leq i,j<n$
given by
$$\left\lbrace\begin{array}{ll}
\displaystyle 1-a&\displaystyle\hbox{if }i=j,\\
\displaystyle -(-1)^{(i-j)/2}a&\displaystyle\hbox{if
}i\not=j\hbox{ and }i\equiv j\pmod 2,\\
\displaystyle -(-1)^{\lfloor\vert i-j\vert/2\rfloor}
+(-1)^{(i+j-1)/2}a&\displaystyle\hbox{if
}i\not\equiv j\pmod 2.\end{array}\right.$$
For $p\geq 2$ an even natural number, one gets thus 
an element (still denoted) $A\in\hbox{GL}_{p-rec}(\mathbb C)$ 
(respectively $A\in\hbox{GL}_{p-rec}(\mathbb Z)$) 
by choosing $a\in\mathbb C$ (respectively $a\in\mathbb Z$)
and setting $A[\mathcal M^l]=
A(p^l)$ if $l\geq 2$ (after appropriate choices for $A[\mathcal M^0],
A[\mathcal M^1]$).

\begin{rem}
The matrix $A(n)$ is singular if $n\equiv 2\pmod 4$ since it 
contains then the vector $(1,-1,1,-1,\dots,1,-1)$ in its kernel.
\end{rem}

\subsubsection{}
There are many possible variations on the above example. 
One can for instance consider the symmetric 
Toeplitz matrix $A(n)$ with coefficients
$$A(n)_{i,j}=\left\lbrace\begin{array}{ll}
\displaystyle a+1&\displaystyle\hbox{if }i=j,\\
\displaystyle a+\vert i-j\vert&\displaystyle\hbox{if
}i\not=j\end{array}.\right.$$
for $0\leq i,j<n$. For $n\equiv 0\pmod 4$, the inverse matrix 
$A(n)^{-1}$ has coefficients $(A(n)^{-1})_{i,j},0\leq i,j<n$
given by
$$\left\lbrace\begin{array}{ll}
\displaystyle (2-n)/2-a&\displaystyle\hbox{if }i=j,\\
\displaystyle -(-1)^{(i-j)/2}(n/2+a)&\displaystyle\hbox{if
}i\not=j\hbox{ and }i\equiv j\pmod 2,\\
\displaystyle (-1)^{\lfloor\vert i-j\vert/2+j\rfloor}(n/2+a-(-1)^j)&\displaystyle\hbox{if
}i\not\equiv j\pmod 2\hbox{ and }i> j,\\
\displaystyle (-1)^{\lfloor\vert i-j\vert/2+i\rfloor}(n/2+a-(-1)^i)&\displaystyle\hbox{if
}i\not\equiv j\pmod 2\hbox{ and }i< j.\end{array}\right.$$
For $p$ even and $a\in\mathbb Z$, we get thus an element (still denoted) 
$A\in\hbox{GL}_{p-rec}(\mathbb Z)$
by setting $A[\mathcal M^l]=A(p^l)$ for $l\geq 2$
(and by defining $A[\mathcal M^0],A[\mathcal M^1]$ suitably).

\subsubsection{}
Similarly, consider a natural odd integer $n$ and define
$\epsilon(n)\in\{\pm 1\}$ such that $n\equiv\epsilon(n)\pmod 4$. 
Let $A(n)$ denote the $n\times n$
matrix with coefficients given by
$$\left(A(n)\right)_{i,j}=
\left\lbrace\begin{array}{ll}
\displaystyle \frac{2-\epsilon(n)-n}{2}+a&\displaystyle \hbox{if }i=j\\
\displaystyle \frac{-\epsilon(n)-n}{2}+\vert i-j\vert+a&\displaystyle \hbox{if }
i\not=j\end{array}\right.$$
For invertible $a$, the matrix $A(n)$ is then invertible;
If $n\equiv 1\pmod 4$ the
inverse matrix has coefficients $\left(A(n)^{-1}\right)_{i,j},0\leq
i,j<n$ given
by
$$\left\lbrace\begin{array}{ll}
\displaystyle 1/a&\displaystyle \hbox{if }i=j\equiv 0\pmod 2,\\
\displaystyle 2&\displaystyle \hbox{if }i=j\equiv 1\pmod 2,\\
\displaystyle (-1)^{(i-j)/2}(1-a)/a&\displaystyle 
\hbox{if }i\not=j\hbox{ and }i\equiv j\equiv 0\pmod 2,\\
\displaystyle (-1)^{(i-j)/2}&\displaystyle 
\hbox{if }i\not=j \hbox{ and }i\equiv j\equiv 1\pmod 2,\\
\displaystyle -(-1)^{\lfloor\vert i-j\vert/2\rfloor}&\displaystyle 
\hbox{if }i\not\equiv j\pmod 2.\end{array}\right.$$

If $n\equiv 3\pmod 4$ the
inverse matrix has coefficients $\left(A(n)^{-1}\right)_{i,j},
0\leq i,j<n$ given by
$$\left\lbrace\begin{array}{ll}
\displaystyle 0&\displaystyle \hbox{if }i=j\equiv 0\pmod 2,\\
\displaystyle 2+1/a&\displaystyle \hbox{if }i=j\equiv 1\pmod 2,\\
\displaystyle -(-1)^{(i-j)/2}&\displaystyle 
\hbox{if }i\not=j\hbox{ and }i\equiv j\equiv 0\pmod 2,\\
\displaystyle (-1)^{(i-j)/2}(a+1)/a&\displaystyle 
\hbox{if }i\not=j \hbox{ and }i\equiv j\equiv 1\pmod 2,\\
\displaystyle -(-1)^{\lfloor\vert i-j\vert/2\rfloor}&\displaystyle 
\hbox{if }i\not\equiv j \pmod 2.\end{array}\right.$$

For $p\geq 3$ an odd natural integer and $a\in\mathbb C^*$, 
one gets thus an element $A\in\hbox{GL}_{p-rec}(\mathbb C)$
by setting $A[\mathcal M^l]=A(p^l)$. The inverse element 
$A^{-1}$ is converging if $p\equiv 1\pmod 4$. Moreover,
$A\in\hbox{GL}_{p-rec}(\mathbb Z)$ for $a=1$ and $a=-1$.

\subsection{}
Another nice example is given by considering the matrix 
$A(n)$ of square size $n\times n$ with coefficients
$$(A(n))_{i,j}=\left\lbrace\begin{array}{ll}
\displaystyle (-1)^{(j-i)/2}&\displaystyle\hbox{if }i\leq j
\hbox{ and }i\equiv j\pmod 2,\\
\displaystyle 0&\displaystyle\hbox{if }i\leq j
\hbox{ and }i\not\equiv j\pmod 2,\\
\displaystyle (-1)^{i-j}&\displaystyle\hbox{if }i\geq j.\end{array}
\right.$$
There are nice formulae for $A(n)^{-1}$ (having all its coefficients
in the finite set $\{0,1,2,3,4\}$) which the reader can easily
write down inspecting the matrices $A(6)^{-1}$ and $A(7)^{-1}$:
$$\left(\begin{array}{cccccc}1&1&2&2&2&1\\
1&2&3&4&4&2\\0&1&2&3&4&2\\
0&0&1&2&3&2\\0&0&0&1&2&1\\0&0&0&0&1&1\end{array}\right),
\left(\begin{array}{ccccccc}1&1&2&2&2&2&1\\
1&2&3&4&4&4&2\\0&1&2&3&4&4&2\\
0&0&1&2&3&4&2\\0&0&0&1&2&3&1\\0&0&0&0&1&2&1\\0&0&0&0&0&1&1
\end{array}\right)\ .$$
Setting $A[\mathcal M^l]=A(p^l)$ for $p\geq 2$
yields thus an element (still denoted) 
$A\in\hbox{GL}_{p-rec}(\mathbb Z)$.
Moreover, $A$ and $A^{-1}$ have very simple
$LU$ decompositions in 
$\hbox{Rec}_{p\times p}(\mathbb Z)$ with $L$ lower triangular
and $U$ upper triangular both unipotent.

\subsection{}\label{exple012012}

Consider the $n\times n$ Toeplitz matrix 
$$T(n)=\left(\begin{array}{rrrrrrrrrrrrrr}
1&-1&1&0&-1&1&0&-1&\dots\\
-2&1&-1&1&0&-1&1&0\\
2&-2&1&-1&1&0\\
0&2&-2&1&-1\\
-2&0&2&\\
\vdots\end{array}\right)$$
with coefficients $T(n)_{i,i}=1,T(n)_{i,j}\equiv i-j\pmod 3$ with
values in $\{-1,0,1\}$ for $i<j$ and with values in $\{-2,0,2\}$
for $i>j$.

For $n\equiv 0\pmod 3$ we have $T(n)\in\hbox{GL}_3(\mathbb Z)$ with
inverse matrix having coefficients $(T(n)^{-1})_{i,j},0\leq i,j<n$
given by
$$(T(n)^{-1})_{i,j}=
\left\lbrace\begin{array}{ll}
\displaystyle -2\cdot n/3+1&\displaystyle \hbox{if }i=j\\
\displaystyle -2\cdot n/3+j-i&\displaystyle \hbox{if }i<j\\
\displaystyle -2\cdot n/3+2(i-j)&\displaystyle \hbox{if
}i>j\end{array}\right.$$
for $0\leq i,j<n$.

For $p\equiv 0\pmod 3$ we get thus a converging element (still
denoted) $T\in\hbox{GL}_{p-rec}(\mathbb Z)$ by
setting $T[\mathcal M^l]=T(p^l)$.

\begin{rem} More generally, one can consider the $n\times n$ 
Toeplitz matrix 
with coefficients $T(n)_{i,i}=a,T(n)_{i,j}\equiv i-j\pmod 3$ with
values in $\{-1,0,1\}$ for $i<j$ and with values in $\{-2,0,2\}$
for $i>j$. The example above corresponds to $a=1$.
Other interesting values are $a=0,a=3$ and $a=-3$.
For $a=3$ for instance, the values of $3^{-n}\det(T(n))$ 
display (experimentally) the following interesting 
$12-$periodic pattern (which one can probably prove
adapting the methods of \cite{ZaPe})
$$\begin{array}{ll}
n\equiv 0,\pm 1\pmod {12}\qquad& 1,\\
n\equiv \pm 2\pmod {12}\qquad&7/9,\\
n\equiv \pm 3\pmod {12}\qquad&5/9,\\
n\equiv \pm 4\pmod {12}\qquad&1/3,\\
n\equiv \pm 5,6\pmod {12}\qquad&1/9.\end{array}$$
\end{rem}

\begin{rem} Another variation on the above example is 
given by considering the matrix $\tilde T(n)$ with coefficients
$T(n)_{i,i}=a+1,T(n)_{i,j}=a+j-i$ if $i<j$ and $T(n)_{i,j}=a+2i-2j$
if $i>j$ (obtained by adding the constant
$a+2n/3$ to all coefficients of the inverse matrix $T(n)^{-1}$
defined above for $n\equiv 0\pmod 3$).

A last variation is given by the matrix defined by $T(n)_{i,i}=x,
T(n)_{i,j}=1+j-i$ if $i<j$ and $T(n)_{i,j}=1+2i-2j$ if $i>j$
(corresponding, up to addition of a constant to all diagonal terms
to the evaluation $a=1$ of the matrix $\tilde T(n)$ considered above). 
Interesting evaluations are: $x=2$ (yielding seemingly unimodular
matrices for $n\equiv 0\pmod 3$, this can probably been proven 
using the ideas of \cite{ZaPe}), $x=4$ (yielding seemingly matrices
of determinant $\det(T(n))=3^{n}$ if $n\equiv 0\pmod 6$ and 
$\det(T(n))=3^{n-2}$ if $n\equiv 3\pmod 6$) and $x=5/2$ yielding seemingly matrices
of determinant $\det(T(n))=(3/2)^{n}$ if $n\equiv 0\pmod 4$ and 
$\det(T(n))=3^{n-2}/2^n$ if $n\equiv 2\pmod 4$).
\end{rem}

\subsection{Digression}
This small digression ``explains'' formulae for the determinants
of several examples treated above. 

Consider the Toeplitz matrix $T(n)$ with coefficients
$T(n)_{i,j},0\leq i,j<n$
given by $T(n)_{i,j}=a+b(j-i)$
if $j>i$, $T(n)_{i,i}=x$, $T(n)_{i,j}=\alpha+\beta(i-j)$ if $i>j$.
We have thus for example
$$T(4)=\left(\begin{array}{cccc}
x&a+b&a+2b&a+3b\\
\alpha+\beta&x&a+b&a+2b\\
\alpha+2\beta&\alpha+\beta&x&a+b\\
\alpha=3\beta&\alpha+2\beta&\alpha+\beta&x\end{array}\right).$$

The generating function $\sum_{n=0}^\infty \det(T(n)t^n$ 
of the associated determinants is then a rational function $P/Q$
(the proof is essentially the same as in \cite{ZaPe}) given by 
$P=\sum_{n=0}^5P_nt^n,Q=\sum_{n=0}^6Q_nt^n$
with
$$\begin{array}{l}
\displaystyle P_0=1\\
\displaystyle P_1=3(a+\alpha)+2(b+\beta)-5x\\
\displaystyle
P_2=3(a+\alpha)^2+2a\alpha+(b^2+b\beta+\beta^2)+3(a+\alpha)(b+\beta)+(ab+\alpha\beta)\\
\displaystyle \qquad -(12(a+\alpha)+6(b+\beta))x+10x^2\\
\displaystyle
P_3=(a+\alpha)^3+4(a^2\alpha+a\alpha^2)+(a+\alpha)^2(b+\beta)+ab(a+b)+\alpha\beta(\alpha+\beta)\\
\displaystyle \qquad
+3a\alpha(b+\beta)+ab\beta+b\alpha\beta\\
\displaystyle \qquad
-(9(a+\alpha)^2+6a\alpha+(b+\beta)^2+6(a+\alpha)(b+\beta)+2(ab+\alpha\beta))x\\
\displaystyle \qquad
+(18(a+\alpha)+6(b+\beta))x^2-10x^3\\
\displaystyle
P_4=(a-x)(\alpha-x)\Big(2(a+\alpha)^2+a\alpha+(a+\alpha)(b+\beta)+(ab+\alpha\beta)+b\beta\\
\displaystyle \qquad
-(7(a+\alpha)+2(b+\beta))x+5x^2\Big)\\
\displaystyle P_5=(a-x)^2(\alpha-x)^2(a+\alpha-x)
\end{array}$$
and
$$\begin{array}{l}
\displaystyle Q_0=1\\
\displaystyle Q_1=3(a+\alpha)+2(b+\beta)-6x\\
\displaystyle Q_2=
3(a+\alpha)^2+3a\alpha+(b+\beta)^2+4(a+\alpha)(b+\beta)\\
\displaystyle \qquad -(15(a+\alpha)+8(b+\beta))x+15x^2,\\
\displaystyle
Q_3=(a+\alpha)^3+6(a\alpha^2+a^2\alpha)+(a+\alpha)(b+\beta)^2+2(a+\alpha)^2(b+\beta)\\
\displaystyle \qquad
+4a\alpha(b+\beta)-(12(a+\alpha)^2+12a\alpha+2(b+\beta)^2+12(a+\alpha)(b+\beta))x\\
\displaystyle \qquad+(30(a+\alpha)+12(b+\beta))x^2-20x^3\\
\displaystyle Q_4=(a-x)(\alpha-x)Q_2\\
\displaystyle Q_5=(a-x)^2(\alpha-x)^2Q_1\\
\displaystyle Q_6=(a-x)^3(\alpha-x)^3\ .\end{array}$$

\subsection{An example given by a symmetric Toeplitz matrix}

We associate to a formal power series $s=\sum_{j=0}^\infty s_jx^j\in
\mathbf K[[x]]$ the symmetric Toeplitz matrix $T(n)$
of size $n\times n$ with coefficients $T(n)_{i,j},0\leq i,j<n$ given by
$$T(n)_{i,j}=[x^{\vert i-j\vert}]s=s_{\vert i-j\vert}\ .$$

Recall that a finite integral square-matrix $A$ is unimodular if
$\det(A)\in\{\pm 1\}$.

\begin{prop} \label{invertsymToep} 
Let $k\geq 1$ be a natural integer. For 
$s=(1-x-x^2)/(1-x-x^2-x^k+x^{k+2})$ and $n\geq 2k+5$,
the matrix $T(n)$ is unimodular. For $j=k+2,\dots,n-k-3$, 
the generating function $g_j=g_j(x)$ of coefficients for the 
$j-$th row (with rows indexed from $0$ to $n-1$) of
$T(n)^{-1}$ is given by 
$$g_j=x^{j-k-2}(1-x-x^2-x^k+x^{k+2})(1-x^2-x^k-x^{k+1}+x^{k+2})\ .$$
The corresponding generating series $g_0,\dots,g_{k+1}$ for
the first $k+2$ rows are given by
$$g_j=(1-x-x^2-x^k+x^{k+2})\left(\sum_{k=0}^jx^k[x^k]\left(
x^{j-k-2}(1-x^2-x^k-x^{k+1}+x^{k+2})\right)\right)$$
where $$\left(\sum_{k=0}^jx^k[x^k]\left(
x^{j-k-2}(1-x^2-x^k-x^{k+1}+x^{k+2})\right)\right)$$
denotes the non-singular part of the 
Laurent-polynomial $$x^{j-k-2}-x^{j-k}-x^{j-2}-x^{j-1}+x^j.$$
The easy identity $g_{n-j}(x)=x^{n-1}
g_j\left(\frac{1}{x}\right)$ determines now the last $k-2$
rows of $T(n)^{-1}$.
\end{prop}

\begin{cor} Choosing an integer $p\geq 2$ and 
a prime $\wp$, reduction $\pmod \wp$ of the matrices $T(p^l)$ 
yields an element (still denoted)  $T\in\hbox{Rec}_{p\times p}(\mathbb
  F_{\wp})$ which is invertible, up to modifications of the
  first few matrices  $T[\mathcal M^l]$ with $p^l<2k+5$.
\end{cor}

\begin{rem} We have $P(1)=-1$ and 
$\hbox{lim}_{x\rightarrow\infty}P(x)=\infty$ for the polynomial 
$P(x)=1-x-x^2-x^k-x^{k+2}$. The polynomial $P(x)$ has thus a real root 
$\rho>1$. For $p\geq 2$, it follows thus from Proposition 
\ref{majorationcriterion} that
the element of $\prod_{l=0}^\infty \mathbb Z^{p^l\times p^l}$ defined by 
$T[\mathcal M^l]=T(p^l)$ is not in $\hbox{Rec}_{p\times
  p}(\mathbb Z)$.
\end{rem}

{\bf Sketch of proof for proposition \ref{invertsymToep}} The identity
$g_{n-j}(x)=x^{n-1}g_j\left(\frac{1}{x}\right)$ for the 
generating functions of the rows of 
$T(n)^{-1}$ follows easily from the fact that $T(n)$ is a symmetric
Toeplitz matrix.

The 
generating series $\tilde f_j=\tilde f_j(x)$ of the $j-$th row of 
$T(n)$ is associated to the coefficients $[x^{-j}]f,\dots,[x^{-j+n-1}]f$
where 
$$f=\frac{1-x-x^2}{1-x-x^2-x^k+x^{k+2}}+\frac{1-x^{-1}-x^{-2}}{
1-x^{-1}-x^{-2}-x^{-k}+x^{-k-2}}-1\in\mathbb Z[[x,x^{-1}]]\ .$$

The proposition follows now from the identities $[x^0]\left(\tilde f_i
\left(\frac{1}{x}\right)g_j(x)\right)=1$ if $i=j$ and 
$[x^0]\left(\tilde f_i\left(\frac{1}{x}\right)g_j(x)\right)=0$ otherwise.
\hfill $\Box$

\begin{rem} Computations suggest that the sequence
$d_n=\det(T(n))$ of determinants 
is given by
$d_n=1$ if $n\leq k$, $d_{k+1}=0$ and $d_n=-1$ if $n>k+1$.

A analogous example is given by the matrices
associated to the function $s=-(1-x^2)/(1-x^2-x^k-x^{k+1}-x^{k+2})$.
The sequence of  determinants $d_n=\det(T(n))$ is
seemingly given by 
$d_n=(-1)^n$ for $n\leq k$, $d_{k+1}=0$ and $d_n=-(-1)^k$ for 
$n>k+1$.

Let us also mention the matrices
$T(n)$ associated to
$s=(1-ax-x^2)/(1-ax-x^2-x^k-(1-a)x^{k+1}+x^{k+2})$. 
Experimentally,
the determinant $d_n=\det(T(n))$ seems to be given by 
$d_n=1$ if $n\leq k$,
$d_{k+1}=0$ and $d_n=-(2a-1)^{n-k-2}$ if $n>k+1$.
\end{rem}

\begin{rem} The formal ``infinite inverse'' matrix $T(\infty)^{-1}$
associated to $T(\infty)$ for $s=(1-x-x^2)/(1-x-x^2-x^k+x^{k+2})$ 
with rows given by $g_j$ as in Proposition 
\ref{invertsymToep} for $n>j+k$ is also interesting.
All submatrices defined by its first $n$ rows and columns are
unimodular and it defines thus a converging unimodular
matrix in $\hbox{Rec}_{p\times p}(\mathbb Z)$ which is not invertible
(for $p>1$) over $\mathbb Z$ but has an invertible 
reduction in $\hbox{Rec}_p(\mathbb F_\wp)$ for $\wp $ an arbitrary 
prime.
\end{rem}

\begin{rem} There are other similar examples, eg. by considering
$s=(1-x^4)/(1-x^2+x^6)),s=(1+x^2-x^4)/(1-x^4+x^6),
s=(1-x^2-x^4)/(1-x^2-x^3-x^4+x^7)$ or
$s=(1-2x+x^2+2x^3-x^4)/(1-2x+4x^3-x^4-2x^5+x^6)$.
It would perhaps be interesting to classify all rational fractions
giving rise to series $s\in\mathbb Z[[x]]$ such that the associated 
symmetric Toeplitz matrices $T(n)$ are
unimodular for almost all $n\in\mathbb N$.

\end{rem}

\subsection{Two more symmetric Toeplitz matrices}

For $\epsilon\in\{1,-1\}$, consider the infinite symmetric 
Toeplitz matrix $T(n)$ with coefficients
$$T(n)_{i,j}=
[x^{\vert i-j\vert}]\left(\epsilon+\frac{x}{1+x^2}\right),
\ 0\leq i,j$$
associated to the series $s=\epsilon+\frac{x}{1+x^2}=\epsilon+x-x^3+x^5-x^7+\dots$.

For $n\equiv 0\pmod 4$ we have $T(n)\in\hbox{GL}_n(\mathbb Z)$ 
with coefficients $(T(n)^{-1})_{i,j},0\leq i,j$ of the inverse
matrix given by
$$
\left\lbrace\begin{array}{ll}
\displaystyle -(n/2-1)\epsilon&\displaystyle \hbox{if }i=j,\\
\displaystyle (n/2-\vert i-j\vert)(-\epsilon)^{i-j+1}
&\displaystyle \hbox{otherwise}.\end{array}
\right.$$ 
For even $p\in\mathbb N$, we get thus a converging element (still denoted)
$T\in\hbox{Rec}_{p\times p}(\mathbb Z)$
by setting $T[\mathcal M^l]=T(p^l)$ for $l\geq 2$ (and defining 
$T[\mathcal M^0],T[\mathcal M^1]$ appropriately).

\begin{rem} For $\epsilon=1$ one gets also an interesting larger
  family by adding a constant $a$ to all coefficients of $T(n)$.
\end{rem}

\section{The Baobab-example related to powers of $2$}

The infinite symmetric Hankel matrix 
$$R=\left(\begin{array}{ccccccccccccc}1&1&0&1&0&0&0&1&\dots\\
1&0&1&0&0\\
0&1&0&\dots\\
\vdots\end{array}\right)$$
with coefficients $R_{i,j}\in\{0,1\},0\leq i,j$ given by
$R_{i,j}=1$ if $(i+j+1)$ is a power of $2$ and $R_{i,j}=0$ 
otherwise defines a converging element (still called) 
$R\in\hbox{Rec}_{2\times 2}(\mathbb Z)$ recursively presented by
$$R=(1,\left(\begin{array}{cc}R&A\\A&0\end{array}\right)),
A=(1,\left(\begin{array}{cc}0&A\\A&0\end{array}\right))\ .$$
Shift matrices with respect to the basis $(R,A)$ (having
initial values $(R,A)[\emptyset,\emptyset]=(1,1)$) of
$\overline{R}^{rec}$ are given by
$$\rho(0,0)=\left(\begin{array}{cc}1&0\\0&0\end{array}\right),
\rho(0,1)=\rho(1,0)=\left(\begin{array}{cc}0&0\\1&1\end{array}\right),
\rho(1,1)=\left(\begin{array}{cc}0&0\\0&0\end{array}\right).$$

The matrix $R$ has an $LU$ decomposition with
$L=L_1\in\hbox{Rec}_2(\mathbf K)$ the lower triangular matrix presented by 
$(L_1,L_2,L_3,L_4)[\emptyset,\emptyset]=(1,1,1,-1)$ and shift matrices
$$\rho(0,0)=\left(\begin{array}{cccc}1&0&1&0\\0&0&0&0\\
0&0&0&0\\0&0&0&0\end{array}\right), 
\rho(0,1)=\left(\begin{array}{cccc}0&0&0&0\\
0&0&0&0\\0&0&0&0\\0&1&0&1\end{array}\right),$$
$$\rho(0,0)=\left(\begin{array}{cccc}0&0&0&0\\1&1&-1&1\\
0&0&0&0\\0&0&0&0\end{array}\right), 
\rho(0,1)=\left(\begin{array}{cccc}0&0&0&0\\
0&0&0&0\\1&1&1&-1\\0&0&0&0\end{array}\right).$$
The upper triangular matrix $U$ is given by $DL^t$ where $D$ is the 
converging diagonal matrix with diagonal limit the $2-$periodic 
sequence $1,-1,1,-1,1,-1,\dots$ recursively presented by
$(S_1,S_2)[\emptyset]=(1,-1)$, $S_1=S_2=(S_1\ S_2)$.

The recurrence matrix $L$ is invertible in $\hbox{Rec}_2(\mathbb Z)$
with inverse $M_1=M=L^{-1}$ presented by
$(M_1,M_2,M_3,M_4)=(1,-1,1,1)$ and 
$$\rho(0,0)=\left(\begin{array}{cccc}1&-1&1&1\\0&0&0&0\\
0&0&0&0\\0&0&0&0\end{array}\right), 
\rho(0,1)=\left(\begin{array}{cccc}0&0&0&0\\
0&0&0&0\\0&0&0&0\\0&1&0&1\end{array}\right),$$
$$\rho(0,0)=\left(\begin{array}{cccc}0&0&0&0\\1&1&-1&1\\
0&0&0&0\\0&0&0&0\end{array}\right), 
\rho(0,1)=\left(\begin{array}{cccc}0&0&0&0\\
0&0&0&0\\1&0&1&0\\0&0&0&0\end{array}\right).$$

Consider the converging diagonal matrix $D'=D'_1\in\hbox{Rec}_{2\times 2}(
\mathbb Z)$ presented by $(D'_1,D'_2)[\emptyset,\emptyset]=(1,1)$
and $$\rho(0,0)=\left(\begin{array}{cc}1&-1\\0&0\end{array}\right),
\rho(0,1)=\rho(1,0)=\left(\begin{array}{cc}0&0\\0&0\end{array}\right),
\rho(1,1)=\left(\begin{array}{cc}0&0\\1&1\end{array}\right).$$
All diagonal coefficients of $D'_1$ (or $D'_2$) are in $\{\pm 1\}$ and
$D'_1$ is thus an element of order $2$ in $\hbox{GL}_{2-rec}(\mathbb
Z)$. The $2-$automatic sequence formed by all diagonal coefficients of
$D'_1$ starts as $1,1,-1,1,-1,-1,-1,1,\dots$ and can be constructed
in the following way: Given two words $W_n,W_n'$ of length $2^n$
in the alphabet $\{\pm 1\}$ construct $W_{n+1}=W_nW_n'$ and
$W_{n+1}'=(-W_n)W_n'$ by concatenating $W_n$ (respectively $-W_n$)
with $W_n'$. The diagonal sequence of $D'_1$ is the limit word
$W_\infty$ obtained from $W_0=W_0'=1$. The converging element $W$ can
thus be recursively presented by $(W,W')[\emptyset]=(1,1)$ and
substition matrices $W=(W\ W'),\ W'=(-W\ W')$.

One can show the following result (see \cite{Bacpapfldng}):

\begin{prop} \label{propbaobab}
(i) The limit of the converging recurrence
matrix $D'LD'\in\hbox{GL}_{2-rec}(\mathbb Z)$ is the infinite
lower triangular matrix $\tilde L$ with coefficients in $\{0,1\}$
defined by $\tilde L_{i,j}\equiv {2i+1\choose i-j}\pmod 2$ 
for $0\leq i,j$.

\ \ (ii) The limit of the converging recurrence matrix 
$DD'MD'D\in\hbox{GL}_{2-rec}(\mathbb Z)$ is the infinite
lower triangular matrix $\tilde L$ with coefficients in $\{0,1\}$
defined by $\tilde M_{i,j}\equiv {i+j\choose 2j}\pmod 2$ for $0\leq i,j$.
\end{prop}

\begin{rem} The recurrence matrices (still denoted) $\tilde L,\tilde M\in
  \hbox{Rec}_{2\times 2}(\mathbb Z)$ defined by Proposition 
\ref{propbaobab} display the following self-similar structure:
Writing
$$\tilde L[\mathcal M_{2\times 2}^l]=\left(\begin{array}{cc}
A\\B&A\end{array}\right),\tilde M[\mathcal M_{2\times 2}^l]=
\left(\begin{array}{cc}A'\\B'&A'\end{array}\right)$$
for $l\geq 1$ (where $A=(\rho(0,0)\tilde L)[\mathcal M_{2\times 2}^{l-1}]=
(\rho(1,1)\tilde L)[\mathcal M_{2\times 2}^{l-1}]$, etc)
we have
$$\tilde L[\mathcal M_{2\times 2}^{l+1}]=\left(\begin{array}{cccc}
A\\B&A\\0&B&A\\B&A&B&A\end{array}\right),
\tilde M[\mathcal M_{2\times 2}^{l+1}]=
\left(\begin{array}{cccc}A'\\B'&A'\\A'&B'&A'\\B'&0&B'&A'
\end{array}\right).$$
The name of this example is motivated by the word
A{\bf BA0BAB}ABA defining the matrix $\tilde L$.
The matrices $\tilde L=A$ and $\tilde M=A'$ have recursive
presentations given by
$$A=(1,\left(\begin{array}{cc}A&0\\B&A\end{array}\right)),
B=(1,\left(\begin{array}{cc}0&B\\B&A\end{array}\right))$$
and
$$A'=(1,\left(\begin{array}{cc}A'&0\\B'&A'\end{array}\right)),
B'=(1,\left(\begin{array}{cc}A'&B'\\B'&0\end{array}\right)).$$
\end{rem}

\section{The Prouhet-Thue-Morse example}

This example has been described in \cite{Bac} and was 
``guiding principle'' and main motivation.

We denote by $\tau(n)$ the integer-valued Prouhet-Thue-Morse
sequence defined by $\tau(n)=\tau\left(\sum_{j=0}^l \epsilon_j2^j\right)=
\sum_{j=0}^l \epsilon_j$ for $n$ a binary integer with binary digits
$\epsilon_0,\dots,\epsilon_l$. The infinite symmetric 
Hankel matrix of the sequence
$i^{\tau(0)},i^{\tau(1)},i^{\tau(2)},\dots$ (where $i^2=-1$)
with generating series 
$\prod_{k=0}^\infty (1+ix^{2^k})$ defines then a converging element
(still denoted) $H\in\hbox{Rec}_{2\times 2}(\mathbb Z[i])$ presented
by $(H=H_1,H_2)[\emptyset,\emptyset]=(1,i)$ and shift matrices
$$\rho(0,0)=\left(\begin{array}{cc}1&i\\0&0\end{array}\right),\
\rho(0,1)=\rho(1,0)=\left(\begin{array}{cc}0&-i\\1&1+i\end{array}\right),
\ 
\rho(1,1)=\left(\begin{array}{cc}i&i\\0&0\end{array}\right)$$

The element $H$ has an $LU$ decomposition with
$L\in\hbox{Rec}_{2\times 2}(\mathbb Z[i])$ unipotent lower triangular
presented by
$(L=L_1,L_2,L_3,L_4)[\emptyset,\emptyset]=(1,i,1,0)$ and shift
matrices
$$\begin{array}{ll}
\displaystyle \rho(0,0)=\left(\begin{array}{cccc}
1&i&1&0\\0&0&0&0\\0&0&0&0\\0&0&0&0\end{array}\right),&
\displaystyle \rho(0,1)=\left(\begin{array}{cccc}
0&0&0&0\\0&0&0&0\\0&0&0&0\\0&1&0&1\end{array}\right),\\
\displaystyle \rho(1,0)=\left(\begin{array}{cccc}
0&-i&-1+i&-i\\1&1+i&-i&1\\0&0&0&0\\0&0&0&0\end{array}\right),&
\displaystyle \rho(1,1)=\left(\begin{array}{cccc}
0&0&0&0\\0&0&0&0\\1&1+i&1&i\\0&i&0&i\end{array}\right).\end{array}$$

We have $U=DL^t\in\hbox{Rec}_{2\times 2}
(\mathbb Z[i])$ for the upper triangular matrix $U\in\hbox{Rec}_{2\times 2}
(\mathbb Z[i])$ where $D\in \hbox{Rec}_{2\times 2}(\mathbb Z[i])$ 
is diagonal with diagonal coefficients given by presentation
$(D=D_1,D_2,D_3)[\emptyset]=(1,1+i,1+i)$ and shift matrices
$$\rho(0,0)=\left(\begin{array}{ccc}1&0&0\\0&0&0\\0&1&1\end{array}\right),
\rho(0,1)=\rho(1,0)=\left(\begin{array}{ccc}0&0&0\\0&0&0\\0&0&0
\end{array}\right),
\rho(1,1)=\left(\begin{array}{ccc}0&2&0\\1&1&1\\0&-2&0\end{array}\right)\
.$$

The inverse $M=L^{-1}$ has presentation $(M=M_1,M_2,M_3,M_4,M_5)
[\emptyset,\emptyset]=(1,-i,1,-i,-1+i)$ and shift-matrices
$$\begin{array}{ll}
\displaystyle \rho(0,0)=\left(\begin{array}{ccccc}
1&0&1&-i&-1\\0&0&0&0&0\\0&0&0&-i&1\\
0&1&0&-i&0\\0&0&0&-i&1\end{array}\right)&
\displaystyle \rho(0,1)=\left(\begin{array}{ccccc}
0&0&0&0&0\\0&0&0&0&0\\0&i&0&1&0\\0&i&0&1&0
\\0&i&0&1&0\end{array}\right)\\
\displaystyle \rho(1,0)=\left(\begin{array}{ccccc}
0&0&0&0&0\\1&0&-i&-i&i\\0&0&1-i&2&-1+i\\
0&0&0&0&0\\0&1&0&1&0\end{array}\right)&
\displaystyle \rho(1,1)=\left(\begin{array}{ccccc}
0&0&0&0&0\\0&0&0&0&0\\1&0&1&i&-1\\
0&0&0&0&0\\0&i&0&i&0\end{array}\right)\end{array}$$

The product $\prod_{k=0}^\infty (1+ix^{2^k})$ has thus 
a continuous fraction of Jacobi-type with coefficients forming
converging $2-$recursive sequences of $\hbox{Rec}_{2}(\mathbb Z[i])$.

\begin{rem} The somewhat similar Hankel matrices with generating
  series $(1-x)\prod_{k=1}^\infty (1+ix^{2^k})$ or 
$(1+x)\prod_{k=1}^\infty (1+ix^{2^k})$ have very similar properties.
\end{rem}

\subsection{The Hankel matrix of
  $\frac{1+x}{1+i}\prod_{k=0}^\infty(1+ix^{2^k})$}

Consider the Hankel matrix $H$ associated to the sequence
$$\beta_1,\beta_2,\dots=1,1+i,i,i,i,-1+i,-1,i,-1+i,-1,\dots$$
defined by
$$\sum_{n=0}^\infty \beta_nx^n=\frac{1+x}{1+i}\prod_{k=0}^\infty
\left(1+ix^{2^k}\right)$$
where we drop the constant term $\beta_0=\frac{1-i}{2}$.

\begin{thm} \label{thmHankelvarianteun}
The determinant $d(n)=\det(H(n))$ of the $n\times n$
Hankel matrix with coefficients $H_{i,j}=[x^{i+j+1}]
\frac{1+x}{1+i}\prod_{k=0}^\infty\left(1+ix^{2^k}\right)$ is
given by $d(n)=(-i)^{\lfloor n/2\rfloor}$.
\end{thm}

The converging Hankel matrix $H\in\hbox{Rec}_{2\times 2}(\mathbb
Z[i])$ is presented by $(H=H_1,H_2)[\emptyset,\emptyset]=(1,1+i)$
with the same shift-matrices
$$\rho(0,0)=\left(\begin{array}{cc}1&i\\0&0\end{array}\right),
\rho(0,1)=\rho(1,0)=\left(\begin{array}{cc}0&-i\\1&1+i\end{array}\right),
\rho(1,1)=\left(\begin{array}{cc}i&i\\0&0\end{array}\right)$$
as for the previous example.

$L$ of $LU$ decomposition: 
$(L=L_1,L_2,L_3,L_4)[\emptyset,\emptyset]=
(1,1+i,1,i)$ and the shift-matrices
$$\begin{array}{ll}
\displaystyle \rho(0,0)=\left(\begin{array}{cccc}
1&i&1&0\\0&0&0&0\\0&0&0&0\\0&0&0&0\end{array}\right)&
\displaystyle \rho(0,1)=\left(\begin{array}{cccc}
0&0&0&0\\0&0&0&0\\0&0&0&0\\0&1&0&1\end{array}\right)\\
\displaystyle \rho(1,0)=\left(\begin{array}{cccc}
0&-i&-1+i&-i\\1&1+i&-i&1\\0&0&0&0\\0&0&0&0\end{array}\right)&
\displaystyle \rho(1,1)=\left(\begin{array}{cccc}
0&0&0&0\\0&0&0&0\\1&1+i&1&i\\0&i&0&i\end{array}\right)
\end{array}$$
are again as in the previous case.

Inverse $M=L^{-1}$ is presented by
$(M=M_1,M_2,M_3,M_4)[\emptyset,\emptyset]=(1,-1-i,1,-1)$ with
shift-matrices
$$\begin{array}{ll}
\displaystyle \rho(0,0)=\left(\begin{array}{cccc}
1&0&1&-i\\0&0&0&0\\0&0&0&0\\0&1&0&-i\end{array}\right)&
\displaystyle \rho(0,1)=\left(\begin{array}{cccc}
0&0&0&0\\0&0&0&0\\0&0&0&0\\0&i&0&1\end{array}\right)\\
\displaystyle \rho(1,0)=\left(\begin{array}{cccc}
0&0&0&0\\1&0&-i&-i\\0&-1&1-i&1\\0&0&0&0\end{array}\right)&
\displaystyle \rho(1,1)=\left(\begin{array}{cccc}
0&0&0&0\\0&0&0&0\\1&-i&1&0\\0&0&0&0\end{array}\right)
\end{array}$$

The diagonal matrix $D$ such that $H=LDL^t$ is the converging
matrix with $2-$periodic diagonal entries $1,-i,1,-i,1,-i,\dots$.

{\bf Proof of theorem \ref{thmHankelvarianteun}} The result 
follows at once from the form of the diagonal matrix $D$ involved in
$H=LDL^t$.\hfill$\Box$ 

\begin{rem} The converging Hankel matrix $H$ presented by 
$(H_1=H,H_2)[\emptyset,\emptyset]=(1,1)$ with shift matrices 
$$\rho(0,0)=\left(\begin{array}{cc}1&i\\0&0\end{array}\right),
\rho(0,1)=\rho(1,0)=\left(\begin{array}{cc}0&-i\\1&1+i\end{array}\right),
\rho(1,1)=\left(\begin{array}{cc}i&i\\0&0\end{array}\right)$$
as above is also in $\hbox{GL}_{2-rec}(\mathbb Q[i])$ with inverse
$R$ (of Hankel type) presented by $(R_1=R,R_2,R_3)[\emptyset,
\emptyset]=(1,\frac{1+i}{2},\frac{-1-i}{2})$ having shift-matrices
$\rho(0,0),\rho(0,1)=\rho(1,0),\rho(1,1)$ given by
$$\left(\begin{array}{ccc}
0&-i&1+i\\0&-1+i&-2i\\i&-1+i&1-2i\end{array}\right),
\left(\begin{array}{ccc}
0&1&-1+i\\1&-1-i&2-i\\0&-1-i&2\end{array}\right),
\left(\begin{array}{ccc}
0&1+i&-1\\0&-2i&1+i\\1&1-2i&1+i\end{array}\right).$$
The Hankel matrix $H$ has however no $LU$ decomposition since (for
instance) the
submatrix consisting of its first $9$ rows and columns is singular.
\end{rem}

\subsection{The Hankel matrix of
  $\frac{x^2-1}{1-i}\prod_{k=0}^\infty(1+ix^{2^k})$}

Let $H$ denote the infinite Hankel matrix associated to the sequence
$$\gamma_2,\gamma_3,\gamma_4,\dots=1,i,0,0,i,-1,-i,1,i,-1,0,0,-1,\dots$$
of coefficients of $\frac{x^2-1}{1-i}\prod_{k=0}^\infty(1+ix^{2^k})$
with $\gamma_0=\frac{-1-i}{2},\gamma_1=\frac{1-i}{2}$ dropped.

The convergent Hankel matrix $H\in\hbox{Rec}_{2\times 2}(\mathbb
Z[i])$ is presented by $(H=H_1,H_2,H_3)[\emptyset,\emptyset]=(1,i,0)$
with shift-matrices $\rho(0,0),\rho(0,1)=\rho(1,0),\rho(1,1)$
as follows
$$\left(\begin{array}{ccc}1&0&i\\0&0&0\\0&1&0\end{array}\right),\ 
\left(\begin{array}{ccc}0&1-i&0\\1&1+i&i\\0&i&0\end{array}\right),\ 
\left(\begin{array}{ccc}0&i&-i\\0&0&0\\1&0&1+i\end{array}\right).$$

The convergent lower triangular unipotent matrix $L$ involved in the 
$LU$ decomposition of $H$ is presented by $(L=L_1,L_2,L_3)
[\emptyset,\emptyset]=(1,i,0)$ with shift-matrices
$$\begin{array}{ll}
\displaystyle \rho(0,0)=\left(\begin{array}{ccc}
1&0&i\\0&0&0\\0&1&0\end{array}\right)&
\displaystyle \rho(0,1)=\left(\begin{array}{ccc}
0&-1&0\\0&-i&0\\0&0&0\end{array}\right)\\
\displaystyle \rho(1,0)=\left(\begin{array}{ccc}
0&1-i&0\\1&1+i&i\\0&i&0\end{array}\right)&
\displaystyle \rho(1,1)=\left(\begin{array}{ccc}
1&0&0\\0&1&0\\0&-i&1\end{array}\right)\end{array}$$

The inverse matrix $M=L^{-1}$ is presented by
$(M=M_1,M_2,M_3,M_4)[\emptyset,\emptyset]=(1,-i,1,0)$ with
shift-matrices
$$\begin{array}{ll}
\displaystyle \rho(0,0)=\left(\begin{array}{cccc}
1&0&1&-i\\0&0&0&0\\0&0&0&0\\0&1&0&-i\end{array}\right)&
\displaystyle \rho(0,1)=\left(\begin{array}{cccc}
0&0&0&0\\0&0&0&0\\0&0&0&0\\0&i&0&1\end{array}\right)\\
\displaystyle \rho(1,0)=\left(\begin{array}{cccc}
0&0&0&0\\1&0&-i&-i\\0&-1&1-i&1\\0&0&0&0\end{array}\right)&
\displaystyle \rho(1,1)=\left(\begin{array}{cccc}
0&0&0&0\\0&0&0&0\\1&-i&1&0\\0&0&0&0\end{array}\right)
\end{array}$$

The diagonal matrix has diagonal coefficients defining
$D\in\hbox{Rec}_2(\mathbf [i])$ presented by
$(D=D_1,D_2,D_3)[\emptyset,\emptyset]=(1,1,i)$ 
with non-zero shift-matrices given by
$$\rho(0,0)=\left(\begin{array}{ccc}
1&0&1\\0&0&0\\0&1&0\end{array}\right),\ 
\rho(1,1)=\left(\begin{array}{ccc}
0&0&-1\\1&1&1\\0&0&1\end{array}\right).$$

\begin{rem} The converging Hankel matrix of the sequence 
$n\longmapsto (-1)^{\tau(n)}$ with generating series $
\prod_{k=0}^\infty (1-x^{2^k})$ seems to be a non-invertible element
of $\hbox{Rec}_{2\times 2}(\mathbb Q)$. The matrix 
$H[\mathcal M^k]$ (with coefficients 
$(-1)^{\tau(i+j)},0\leq i,j<2^k$), presented by
$(H_1=H,H_2)[\emptyset,
\emptyset]=(1,-1)$ and shift-matrices
$$\rho(0,0)=\left(\begin{array}{cc}1&-1\\0&0\end{array}\right),
\rho(0,1)=\rho(1,0)=\left(\begin{array}{cc}0&1\\1&0\end{array}\right),
\rho(1,1)=\left(\begin{array}{cc}-1&-1\\0&0\end{array}\right),$$
seems however to
have a nice and interesting inverse matrix $(H[\mathcal M^k])^{-1}$
given by the Hankel matrix associated to the
sequence $\frac{1}{2}\sigma(k),-\frac{1}{2}\sigma(k)$ (with
last term unused) where
$\sigma(k)$ is of length $2^k$ and is recursively defined by
$\sigma(0)_1=2$, $\sigma(k+1)_1=\frac{1}{2}\sum_{i=1}^{2^k}\sigma(k)_i$,
$\sigma(k+1)_{2i+1}=\sigma(k+1)_{2i}$ for $i>0$ and
$\sigma(k+1)_{2i}=\sigma(k+1)_{2i-1}-\sigma(k)_i$.

The first sequences $\sigma(k)$ are:

$\sigma(0)=2$

$\sigma(1)=(1,-1)$

$\sigma(2)=(0,-1,-1,0)$

$\sigma(3)=(-1,-1,-1,0,0,1,1,1)$

$\sigma(4)=(0,1,1,2,2,3,3,3,3,3,3,2,2,1,1,0)$.
\end{rem}


I would like to thank Jean-Paul Allouche, Laurent Bartholdi,
Michel Brion, Rostislav Grigorchuk, Pierre de la
Harpe, Jeffrey Shallit and many other people for interesting
discussions and remarks related to this paper.




\noindent Roland BACHER

\noindent INSTITUT FOURIER

\noindent Laboratoire de Math\'ematiques

\noindent UMR 5582 (UJF-CNRS)

\noindent BP 74

\noindent 38402 St Martin d'H\`eres Cedex (France)
\medskip

\noindent e-mail: Roland.Bacher@ujf-grenoble.fr


\begin{thebibliography}{99}

\bibitem{AS} J.-P. Allouche, J. Shallit, Automatic Sequences.
Theory, Applications, Generalizations, Cambridge
University Press (2003).

\bibitem{Bac} R. Bacher, {\it La suite de Thue-Morse et la cat\'egorie 
$\hbox{Rec}$}, C. R. Acad. Sci. Paris, Ser. I.

\bibitem{Bacpapfldng} R. Bacher, {\it Paperfolding and Catalan
    numbers}, arXiv: math.CO/0406340.

\bibitem{BGN} L. Bartholdi, R. Grigorchuk, V. Nekrashevych, {\it From 
Fractal Groups to Fractal Sets}, arXiv:math.GR/0202001 v4 (2001).

\bibitem{delaHarpe} P. De la Harpe, Topics in Geometric Group Theory,
Chicago Lectures in Mathematics (2000).

\bibitem{Flajolet} P. Flajolet, {\it Combinatorial aspects of
    continued fractions}, Discr. Math., {\bf 32} (1980), 125-161.



\bibitem{Gri} R.I. Grigorchuk, {\it Burnside's problem on
periodic groups}, Funct. Anal. Appl. {\bf 14} (1980), 41-43.

\bibitem{GrPa} R.I. Grigorchuk, I. Pak, {\it Groups of Intermediate
    Growth: an Introduction for Beginners}, arXiv: math.GR/0607384v1

\bibitem{SaMaLa} Mac Lane, Saunders,
Categories for the working mathematician. Springer (1998).

\bibitem{Nek} V. Nekrashevych, Self-similar groups, book

\bibitem{ZaPe}
H. Zakraj\v sek, M. Petkov\v sek, {\it
Pascal-like determinants are recursive},
Adv. in Appl. Math. {\bf 33} (2004), no. 3, 431-450. 

\end{thebibliography}
\end{document}